\documentclass{amsart} 
\usepackage{calrsfs}
\usepackage{fullpage}

\newcommand{\ord}{\mbox{ord}}
\newcommand{\Gal}{\mbox{Gal}}

\newcommand{\mte}[1]{\mbox{\te {#1}}}

\newcommand{\be}{\begin{enumerate}}
\newcommand{\ee}{\end{enumerate}}
\newcommand{\id}{\mbox{id}}

\newcommand{\rank}{\mbox{rank}}

\newcommand{\dg}{\leq_{Dioph}}

\newcommand{\calA}{{\mathcal A}}
\newcommand{\calB}{{\mathcal B}}
\newcommand{\calC}{{\mathcal C}}
\newcommand{\calD}{{\mathcal D}}
\newcommand{\calE}{{\mathcal E}}

\newcommand{\calL}{{\mathcal L}}
\newcommand{\calM}{{\mathcal M}}
\newcommand{\calN}{{\mathcal N}}

\newcommand{\calP}{{\mathcal P}}

\newcommand{\calR}{{\mathcal R}}
\newcommand{\calS}{{\mathcal S}}
\newcommand{\calT}{{\mathcal T}}
\newcommand{\calU}{{\mathcal U}}
\newcommand{\calV}{{\mathcal V}}
\newcommand{\calW}{{\mathcal W}}
\newcommand{\calX}{{\mathcal X}}
\newcommand{\calY}{{\mathcal Y}}
\newcommand{\calZ}{{\mathcal Z}}

\newcommand{\C}{{\mathbb C}}
\newcommand{\F}{{\mathbb F}}

\newcommand{\N}{{\mathbb N}}
\newcommand{\Q}{{\mathbb Q}}
\newcommand{\R}{{\mathbb R}}
\newcommand{\Z}{{\mathbb Z}}

\newcommand{\pp}{{\mathfrak p}}
\newcommand{\Pp}{{\mathfrak P}}
\newcommand{\qq}{{\mathfrak q}}

\newcommand{\ttt}{{\mathfrak t}}

\newtheorem{theorem}{Theorem}[section]
\newtheorem{lemma}[theorem]{Lemma}
\newtheorem{corollary}[theorem]{Corollary}
\newtheorem{proposition}[theorem]{Proposition}

\theoremstyle{definition}
\newtheorem{definition}[theorem]{Definition}

\theoremstyle{remark}
\newtheorem{remark}[theorem]{Remark}

\newtheorem{notation}[theorem]{Notation}
\newtheorem{notationassumption}[theorem]{Notation and Assumptions}
\newtheorem{assumption}[theorem]{Assumptions}
\font\te=eufm10
\newcommand{\su}{\subsection*}

\begin{document}

\bibliographystyle{plain}%
  \title{Diophantine Definability and Decidability in the Extensions of Degree 2 of Totally Real Fields}
\date   {June 02, 2006}
\author{Alexandra Shlapentokh}

\thanks{The research for this paper has been partially supported by NSF grant DMS-0354907 and ECU
Faculty Senate Summer 2005 Research Grant.}
\address{Department of Mathematics \\ East Carolina University \\ Greenville, NC 27858}
\email{shlapentokha@ecu.edu }
\urladdr{www.personal.ecu.edu/shlapentokha}
\maketitle

\begin{abstract}%
We investigate Diophantine definability and decidability over some subrings of algebraic numbers contained in
quadratic  extensions of totally real algebraic extensions of  $\mathbb Q$.  Among other results we prove the
following.  The big subring definability and undecidability results previously shown by the author to hold over
totally complex extensions of degree 2 of totally real number fields, are shown to hold for {\it all} extensions
of degree 2 of totally real number fields.  The definability and undecidability results for integral closures of
``small'' and ``big'' subrings of number fields in the infinite algebraic extensions of $\mathbb Q$, previously
shown by the author to hold for totally real fields, are extended to a large class of extensions of degree 2 of
totally real fields.  This class includes infinite cyclotomics  and abelian extensions with finitely many ramified
rational primes.
\end{abstract}%
\maketitle

\section{\bf Introduction}
The interest in the questions of existential definability and decidability over rings goes back to a
question that was posed by Hilbert: given an arbitrary polynomial equation in several variables over
$\Z$, is there a uniform algorithm to determine whether such an equation has solutions in
$\Z$? This question, otherwise known as Hilbert's Tenth Problem, has been answered negatively
in the work of M. Davis, H. Putnam, J. Robinson and Yu. Matiyasevich. (See \cite{Da1} and
\cite{Da2}.)   Since the time when this result was obtained, similar questions have been raised for
other fields and rings. In other words, let $R$ be a recursive ring.  Then, given an arbitrary
polynomial equation in several variables over $R$, is there a uniform algorithm to determine
whether such an equation has solutions in $R$?
    One way to resolve the question of Diophantine decidability negatively over a ring of
characteristic 0 is to construct a Diophantine definition of $\Z$ over such a ring.  This
notion is defined below.

\begin{definition}
Let $R$ be a ring and let $A \subset R$.  Then we say that $A$ has a Diophantine definition over
$R$ if there exists a polynomial $f(t,x_1,\ldots,x_n) \in R[t,x_1,\ldots,x_n]$ such that for any $t
\in R$,
\[
\exists x_1,\ldots,x_n \in R, f(t,x_1,...,x_n) = 0 \Longleftrightarrow t \in A.
\]

If the quotient field of $R$ is not algebraically closed,   we can allow
a Diophantine definition to consist of several polynomials without changing the nature of the
relation. (See \cite{Da2} for more details.)
\end{definition}

The usefulness  of Diophantine definitions stems from the following easy lemma.
\begin{lemma}
Let $R_1 \subset R_2$ be two recursive rings such that the quotient field of $R_2$ is not
algebraically closed. Assume that  Hilbert's Tenth Problem (abbreviated as ``HTP'' in the future)   is
undecidable over $R_1$, and $R_1$ has a Diophantine definition over $R_2$. Then HTP is
undecidable over $R_2$.%
\end{lemma}

Using norm equations, Diophantine definitions have been obtained for $\Z$ over the rings of algebraic integers of
some number fields. Jan Denef constructed a Diophantine definition of $\Z$ for the finite degree totally real
extensions of $\Q$. Jan Denef and Leonard Lipshitz extended Denef's results to all the extensions of degree 2 of
the finite degree totally real fields. Thanases Pheidas and the author of this paper have independently
constructed Diophantine definitions of $\Z$ for number fields with exactly one pair of non-real conjugate
embeddings. Finally Harold N. Shapiro and the author of this paper showed that the subfields of all the fields
mentioned above ``inherited'' the Diophantine definitions of $\Z$. (These subfields include all the abelian
extensions.) The proofs of the results listed above can be found in \cite{Den1}, \cite{Den2}, \cite{Den3},
\cite{Ph1}, \cite{Sha-Sh}, and \cite{Sh2}.  The author also showed that the  totally real fields which are
non-trivial extensions of $\Q$, and the totally complex extensions of degree 2 of the totally real fields contain
``big'' rings, i.e. rings consisting of  algebraic numbers with infinitely many primes allowed in the
denominators of divisors,  where $\Z$  is definable and Hilbert's Tenth Problem has no solution.  The details of
these results can be found in \cite{Sh1}, \cite{Sh5} and \cite{Sh3}.  Subsequently these results were extended by
the author to the integral closures of some  rings of $\calS$-integers where $\calS$ is finite ( in the future
referred to as ``small'' rings)   and to the integral closures of some ``big'' rings in  a class of infinite
totally real extensions of $\Q$ (see \cite{Sh17} and \cite{Sh26}.)

The investigation of ``big'' rings was prompted by difficulties of resolving the status of Hilbert's Tenth
Problem over $\Q$. These difficulties, in part, gave rise to a series of conjectures by Barry Mazur which can be
found in \cite{M1}, \cite{M2}, \cite{M3}, and \cite{M4}. The strongest of the conjectures in \cite{M1} was refuted
in \cite{CTSDS}, but the status of the other conjectures as well as the Diophantine status of $\Q$ and other number
fields is still unknown. Among other things, Mazur's conjectures  implied that infinite discrete (in Archimedean or
$\mte{p}$-adic topology) sets are not existentially definable over $\Q$ (or other number fields) and thus $\Z$ is
not Diophantine over $\Q$. Cornelissen and Zahidi showed that one of the conjectures also implied that another
possible method of proof of Diophantine undecidability of $\Q$ was not viable, i.e. they showed that one of
Mazur's conjectures implied that  there was no Diophantine model of $\Z$ over $\Q$. (See \cite{CZ} for more
information on this result.)

In arguably the most important development in the subject since the solution of the original problem, Poonen showed
  that there exist ``really big'' recursive subrings  of $\Q$ (that is recursive subrings with the natural density
of primes allowed in the denominators equal to 1), where one could construct an infinite, existentially definable
over the ring discrete in archimedean topology set, which is also a Diophantine model of $\Z$.  Thus, Poonen
showed that a ring version of a  Mazur's conjecture failed for this ring and Hilbert's Tenth Problem was
undecidable over the ring.  (See  \cite{Po2} for more details.)   In a paper joint with the author (see
\cite{PS}), this result was lifted to any number field which has a rank one elliptic curve. Poonen and the author
also showed in \cite{PS} that some ``really big'' subrings of  number fields  with a rank one elliptic curve ($\Q$
being one of these fields) had Diophantine sets which were simultaneously discrete in the usual archimedean and
every non-archimedean topology of the field.  Additionally, the paper contained examples of ``big'' rings with
$\mte{p}$-adically discrete Diophantine sets  contained in totally real number fields and their totally complex
extensions of degree 2.

Elliptic curves have also been used to show undecidability of rings of algebraic integers.  The first use of
elliptic curves for this purpose is due to Denef  who proved the following proposition in \cite{Den2}.%
\begin{theorem}%
Let $K_{\infty}$ be a totally real algebraic possibly infinite extension of $\Q$.  If there exists an elliptic
curve $\calE$  over $\Q$ such that $[\calE(K):\calE(\Q)] < \infty$, then $\Z$ has a Diophantine definition over
the rings of algebraic integers of $K$.
\end{theorem}%

Extending ideas of Denef, Bjorn Poonen has shown  the following in \cite{Po}.%
\begin{theorem}%
\label{thm:po}%
Let $M/K$ be a number field extension with an elliptic curve $\calE$ defined over $K$, of rank one over $K$, such
that the rank of $\calE$ over $M$ is also one.  Then $O_K$ (the ring of integers of $K$) is Diophantine over $O_M$.\\%
\end{theorem}%

In a recent paper (see \cite{CPZ}), Cornelissen, Pheidas and Zahidi weakened somewhat assumptions of Poonen's theorem. Instead
of requiring a rank 1 curve retaining its rank in the extension, they require existence of a rank 1 elliptic curve over the
bigger field and an abelian variety over the smaller field retaining its rank in the extension. Further, Poonen and the author
have independently shown that the conditions of Theorem \ref{thm:po} can be weakened to remove the assumption that rank is one
and require only that the rank in the extension is the same (see \cite{Sh33} and \cite{Po3}).  In \cite{Sh33}, the author also
showed that the elliptic curve technique extends to ``big rings''.\\

In this paper, using norm equations, we extend the ``big'' ring results to \emph{all} extensions of degree two of
totally real number fields and some totally real infinite extensions of $\Q$. In the case of infinite extensions,
we will all also obtain new results for rings of $\calS$-integers but not for rings of integers. (Corresponding
results for rings of integers can be obtained via elliptic curve methods as in \cite{Sh33}. We intend to do this
in the future.) En route to the results above we obtain  improvements relative to \cite{Sh26} for the results
concerning  ``big'' and ``small'' rings of some totally real infinite extensions of $\Q$, as well as results on
integrality at finitely many ``primes'' in infinite extensions. The main theorems of the paper are stated below.

\su{Theorem}
Let $G$ be an extension of degree 2 of a totally real number field. Then for any $\varepsilon
>0$ there exists a set $\calV_G$ of primes of $G$ whose natural density is bigger than $1-1/[G:\Q] - \varepsilon$
and such that $\Z$ has a Diophantine definition over $O_{G,\calV_G}$.\\

\su{Theorem}%
Let $A_{\infty}$ be an abelian (possibly infinite) extension of $\Q$ with finitely many ramified primes. Then the
following statements are true.
\begin{itemize}%
\item If the ramification degree of 2 is finite, then for any number field $A \subset A_{\infty}$ there exists an
infinite set of $A$-primes $\calW_A$ such that $\Z$ is existentially definable in the integral closure of
$O_{A,\calW_A}$ of $A_{\infty}$.%
\item For any number field $A \subset A_{\infty}$ and any finite non-empty set $\calS_A$ of its
primes, we have that $\Z$ is existentially definable in the integral closure of $O_{A,\calS_A}$ in $A_{\infty}$.
\end{itemize}%

\su{Theorem}%
Let $q$ be a rational prime. Let $L$ be an algebraic, possibly infinite extension of $\Q$. Let $\Pp_L$ be a prime
of $L$ (a prime ideal of $O_L$ -- the ring of algebraic integers of $L$) such that it is  relatively prime to $q$
(meaning the ideal does not contain $q$), the residue field of $\Pp_L$ has an extension of degree $q$, and for any
number field $M \subset L$, it is the case that any $M$-prime $\pp_M$ lying below $\Pp_L$ is unramified over
$\Q$. Then for any number field $M \subset L$, there exists  a subset $\calX$ of $L$ satisfying the following
conditions:
\begin{itemize}%
\item If $x \in \calX$, then $x$ is integral with respect to $\pp_M$, an  $M$-prime below $\Pp_L$.%
\item If $x \in M$ and $x$ is integral at $\pp_M$, then $x \in \calX$.%
\item $\calX$ is Diophantine over $L$.%
\end{itemize}%

 \section{\bf Preliminary Results.} %
 \label{sec:prelim}
In this section we state some definitions and a few well-known technical propositions which will be used in the
proofs. We start with a definition of ``big'' rings and
integrality at finitely many primes in a number field.
\begin{definition}%
\label{def:rings}%
Let $K$ be a number field.  Let $\calW_K$ be a set of its non-archimedean primes.  Then define
\[%
O_{K,\calW_K}=\{x \in K: \ord_{\pp}x\geq 0, \, \forall \pp \not \in \calW_K\}.
\]%
If $\calW_K$ is empty the ring $O_{K,\calW_K}=O_K$ is the ring of integers of $K$. If $\calW_K$ is finite, then the
ring $O_{K,\calW_K}$ is called the ring of $\calW_K$-integers or a ``small'' ring.  If $\calW_K$ is infinite, we
will call  the ring $O_{K,\calW_K}$ a ``big'' ring.
\end{definition}%
\begin{proposition}%
\label{prop:finmany}%
Let $K$ be a number field. Let $\calW_K$ be any set of primes of $K$. Let $\calS_K \subseteq
\calW_K$ be a finite set.  Let $\calV_K = \calW_K \setminus \calS_K$.  Then
$O_{K,\calV_K}$ has a Diophantine definition over $O_{K,\calW_K}$.  (See, for example, \cite{Sh5}.)\\
\end{proposition}%

Next we state another technical proposition which is also quite important for the proofs in this paper.%
\begin{proposition}%
\label{prop:non-zero}%
Let $K$ be a number field. Let ${\mathcal W}_K$ be any set of primes of $K$. Then the set of non-zero elements of
$O_{K,\calW_K}$ has a Diophantine definition over $O_{K,\calW_K}$. Further, let $K_{\infty}$ be an
algebraic extension of $K$ and let $O_{K_{\infty},\calW_{K_{\infty}}}$ be the integral closure of $O_{K,\calW _K}$
in $K_{\infty}$. Then the set of non-zero elements of $O_{K_{\infty},\calW_{K_{\infty}}}$ has a Diophantine
definition over $O_{K_{\infty},\calW_{K_{\infty}}}$ (See, for example, \cite{Sh5} and \cite{Sh26}.)\\
\end{proposition}%

We will also need the following easy proposition:
\begin{proposition}%
\label{prop:relprime}
Let $F$ be an algebraic extension of $\Q$.  Let $O_F$ be the ring of integers of $F$.  Then the following set of
pairs of elements of $O_F$ is Diophantine over $O_F$: $\{(a,b) \in (O_F)^2: (a,b)=1\}$.
\end{proposition}%
\begin{proof}%
It is enough to note that $(a,b) =1 \Leftrightarrow (\exists A, B \in O_F)(Aa+Bb)=1$.
\end{proof}%

The following proposition will allow us to set up bounds for real valuations.

\begin{proposition}%
\label{prop:positive}%
Let $F$ be an algebraic extension of $\Q$. Let $P=\{x \in F| \mbox{ For all embeddings } \sigma: F \longrightarrow \R,
\sigma(x) \geq 0\}$.  Then $P$ is Diophantine over $F$. (See \cite{Den2}, Lemma 10.)
\end{proposition}%

The next proposition deals with rewriting equations using variables from finite degree subfields.

\begin{proposition}%
\label{prop:rewriting}%
Let $F_2/F_1$ be a finite field extension. Let $R_1 \subseteq F_1$ be a ring whose fraction field is $F_1$ and let
$R_2$ be the integral closure of $R_1$ in $F_2$.  Assume further that the set of non-zero elements is Diophantine
over $R_1$.  Let
\begin{equation}%
\label{sysr:1}
P(X_1,\ldots,  X_n, z_1,\ldots,z_m)=0
\end{equation}%
be an equation with coefficients in $F_2$. Then for some positive integers $l$ and $r$, there exists a system of equations
\begin{equation}%
\label{sysr:2}
\{Q_i(t_1,\ldots,t_r, z_1,\ldots,z_m)=0, i=1,\ldots, l\}
\end{equation}%
over $R_1$ such that (\ref{sysr:2}) has solutions $t_1,\ldots,t_r,z_1,\ldots,z_m \in R_1$ if and only if (\ref{sysr:1}) has
solutions $X_1,\ldots,X_n$ in $R_2$, $z_1,\ldots,z_m \in R_1$.
\end{proposition}%
\begin{proof}%
The rewriting proceeds in several steps.  First of all let $\{\omega_1,\ldots, \omega_k\}$ be a basis of $F_2$
over $F_1$.  Let
\begin{equation}%
\label{sysr:3} P(X_1,\ldots,  X_n, z_1,\ldots,z_m)=\sum_{i_1,\ldots, i_n,e_1,\ldots,e_m}A_{i_1,\ldots, i_n,
e_1,\ldots,e_m}X_1^{i_1}\ldots X_n^{i_n}z_1^{e_1}\ldots z_m^{e_m}=0.
\end{equation}%
Next note that $A_{i_1,\ldots,i_n} = \sum_{j=1}^kb_{i_1,\ldots,i_n,e_1,\ldots,e_m,j}\omega_j$, where
$b_{i_1,\ldots,i_n,j} \in F_1$.   Now we replace (\ref{sysr:3}) by

\begin{equation}%
\label{sysr:4}%
 \sum_{i_1\ldots i_n,e_1\ldots e_m}\left(\sum_{j=1}^kb_{i_1,\ldots,i_n,e_1,\ldots, e_m, j}\omega_j\right )\left
(\sum_{j=1}^kx_{1,j}\omega_j\right )^{i_1}\ldots \left (\sum_{j=1}^kx_{n,j}\omega_j\right )^{i_n}z_1^{e_1}\ldots z_m^{e_m}=0.%
\end{equation}%
The next step is to make sure that for each $i$ we have that $\sum_{j=1}^kx_{i,j}\omega_j$ is in the integral
closure of $R_1$ and all the variables range over $R_1$. To reach the latter goal we will replace each $x_{i,j}$
by a ratio $u_{i,j}/v_{i,j}$ where $u_{i,j}, v_{i,j}$ will range over $R_1$. To insure that
$\sum_{j=1}^k(u_{i,j}/v_{i,j})\omega_j$ is in $R_2$, the integral closure of $R_1$ in $F_2$,  we add an equation
requiring that the all the symmetric functions of the conjugates of the sum over $F_1$ are in $R_1$ or
alternatively that the coefficients of the characteristic polynomial are all in $R_1$. We also add equations
stating that $v_{i,j} \not = 0$ as elements of $R_1$.

The third step is multiply all the factors out treating $z_j, u_{i,j}/v_{i,j}$'s as elements of  $F_1$ and to
replace products of the elements of the basis by their linear combinations over $F_1$.  This operation will
produce a linear combination of $\omega$'s with coefficients which are polynomials in $u_{i,j}/v_{i,j}$.  The last
step is then to multiply by appropriate powers of $v_{i,j}$'s and denominators of $b_{i_1,\ldots,i_n,e_1,\ldots,
e_m, j}$ with respect to $R_1$ to clear all the denominators.
\end{proof}%
Now we state a slightly different version of the rewriting proposition whose proof is completely analogous to the
proof above.
\begin{proposition}%
\label{prop:rewriting3}%
Let $F_2/F_1$ be a finite field extension. Let $R_1 \subseteq F_1$ be a ring whose fraction field is $F_1$ and let
$R_2=R_1[\nu]$, where $\nu$ generates $F_2$ over $F_1$ and $\nu$ is integral over $R_1$.  Let
\begin{equation}%
P(X_1,\ldots,  X_n, z_1,\ldots,z_m)=0
\end{equation}%
be an equation with coefficients in $F_2$. Then for some positive integers $l$ and $r$, there exists a system of equations
\begin{equation}%
\{Q_i(t_1,\ldots,t_r, z_1,\ldots,z_m)=0, i=1,\ldots, l\}
\end{equation}%
over $R_1$ such that (\ref{sysr:2}) has solutions $t_1,\ldots,t_r,z_1,\ldots,z_m \in R_1$ if and only if (\ref{sysr:1}) has
solutions $X_1,\ldots,X_n$ in $R_2$, $z_1,\ldots,z_m \in R_1$.
\end{proposition}%
 Using similar reasoning, one can also prove the following easy proposition.
\begin{proposition}%
\label{prop:rewriting2}%
Let $F_1 \subset F_2 \subset F_3$ be a finite extension of number fields. Let $R_1 \subset F_1$ be an integrally
closed subring of $F_1$ such that its fraction field is $F_1$. Let $\nu_2 \in F_2$ be a generator of $F_2$ over
$F_1$ such that it is  integral over $R_1$.  Similarly, let $\nu_3 \in F_3$ be a generator of $F_3$ over $F_2$
such that it is integral over $R_1[\nu_2]$. Let $R_2=R_1[\nu_2], R_3=R_1[\nu_2,\nu_3]$. Then for some positive
integers $l$ and $r$ there exists a system of polynomial equations
\begin{equation}%
\label{eq:replace}%
\{Q_i(t_1,\ldots,t_l)=0, i=1,\ldots,r\}
\end{equation}%
with coefficients in  $R_1$, such that ${\mathbf N}_{F_3/F_2}(\varepsilon)=1$ has a solution $\varepsilon \in R_3$
if and only if (\ref{eq:replace}) has a solution $(t_1,\ldots,t_l) \in R_1^l$.
\end{proposition}%

We finish this section with a notational convention which we will observe throughout the paper.
\begin{notation}
\label{not:algclosure}%
Let $\tilde{\Q}$ be
the algebraic closure of $\Q$ inside $\C$. All the algebraic extensions of $\Q$ discussed in the paper will be
assumed to be subfields of $\tilde{\Q}$. Further, given a finite set of fields $F_1,\ldots, F_k \subset
\tilde{\Q}$, we will interpret $F_1\ldots F_k$ to mean the smallest subfield of $\tilde{\Q}$ containing
$F_1,\ldots, F_k$.
\end{notation}%

\section{\bf Integrality at a Prime in Infinite Extensions.}%
This section contains some technical material necessary for the proofs for the infinite extension cases.  However,
the result may be of independent interest.  More specifically we will discuss existential definability of
integrality at finitely many primes in infinite extensions.  We will rely heavily on Theorem 2.1 of \cite{Sh23}
which is the technical version of the assertion that integrality at finitely many primes is existentially
definable over number fields.
\begin{notationassumption}%
\label{not:prime}
In this section we will use the following notation and assumptions.
\begin{itemize}%
\item Let $K$ be a number field. %
\item Let $F$ be an algebraic (possibly infinite) extension of $K$.%
\item Let $q>2$ be a rational prime number. Let $\xi_q$ be a $q$-th primitive root of unity.%
\item Assume $\xi_q \in K$.%
\item Let $b \in O_K$ and assume that $X^q -b$ is irreducible over $F$.%
\item Let $\pp_K$ be a prime of $K$ satisfying the following conditions.
\begin{itemize}%
\item $\pp_K$ is not a factor of $q$.%
\item Let $\pp_M$ be any factor  of $\pp_K$ in some finite extension $M$ of $K$ such that $M \subset F$. Then  $X^q
-b$ is irreducible in the residue field of $\pp_M$ and the ramification degree of $\pp_M$  over $\pp_K$ is not
divisible by $q$.  We will separate out the case of ramification degree being 1 for all $M$ and all $\pp_M$, and
will refer to this case as the ``unramified" case.  Note also that the irreducibility assumption implies that
$\ord_{\pp_M}b=0$.
\end{itemize}%
\item Let $\displaystyle \prod \mte{a}_i^{r_i}$ be the $K$-divisor of $b$. For each $i$, let $A_i$ be the rational
prime below $\mte{a}_i$.%
\item Let $s \geq \max_i\{3qr_i[K:\Q]\}$ be a natural number not divisible by $q$.%
\item  Let $g \in K$ satisfy the following conditions.%
\begin{itemize}%
\item $\ord_{\pp_K}g=s\not \cong 0 \mod q$.%
\item $g \cong 1 \mod bq^3$.%
\item The divisor of $g$ is of the form $\displaystyle \frac{\pp_K^s}{\prod\qq_i^{n_i}}$, where for all $i$ we have
that $\qq_i$ is a prime of $K$ such that $\qq_i \not = \pp$ and $n_i \in \Z_{>0}$. (Such a $g \in K$ exists by
the Strong Approximation Theorem.)
\end{itemize}%

\item Let $\displaystyle r=q^{(3q[K:\Q])}(q^{(q[K:\Q])!}-1)\left(\prod (A_i^{(q[K:\Q])!}-1)\right )$%

\item For $x \in F$, let $h =(q^3b)^r(g^{-1}x^{r(s+1)} + g^{-q})+1$. %
\item Let $z \in O_K, z \not \not \cong 0 \mod \pp$.%
\item Let $\beta(x) \in \tilde{\Q}$ be a root of $T^q - (h^{-1} + z^q)$. \item Let $\beta \in \tilde{\Q}$ be a root
of the polynomial $X^q-b$.%
\item Let $\displaystyle N(a_0,\ldots,a_{q-1})=\prod_{j=1}^{q}\left (\sum_{i=0}^{q-1}a_i\xi_q^{ij}\beta^i\right)$.%
\item If $y \in M$, where $M$ is a finite extension of $K$, then we will say that ``$y$ is integral at $\pp_K$'' if $y$
is integral at every factor of $\pp_K$ in $M$.
\end{itemize}%
\end{notationassumption}%

We are now ready to state and prove the main technical proposition of this section.\\
\begin{proposition}%
\label{prop:main}
 Let $x \in F$.  Let $M \subset F$ be a number field containing $K(x)$. Then
\begin{equation}%
\label{eq:n}%
 N(a_0,\ldots,a_{q-1})=h
\end{equation}%
has solutions $a_0,\ldots,a_{q-1} \in L_x=F(\beta(x))$ only if for all fields $M$ as above, for all  factors
$\pp_M$ of $\pp_K$ in $M$, we have that
\begin{equation}%
\label{eq:orderbound}
\ord_{\pp_M}x > \frac{(-q+1)\ord_{\pp_M}(g)}{r(s+1)}.
\end{equation}%
In the unramified case we can make a stronger statement: equation (\ref{eq:n}) has solutions $a_0,\ldots,a_{q-1}
\in L_x=F(\beta_x)$ only if $x$ is integral at every factor of $\pp_K$. (Note that (\ref{eq:orderbound}) is
automatically satisfied if $x$ is integral at all the factors of $\pp_K$ in $F$.) Finally, if $x \in K$ and is
integral at $\pp_K$, then $(\ref{eq:n})$ has solutions $a_0,\ldots, a_{q-1} \in K(\beta_x)$.
\end{proposition}%
\begin{proof}%
We start with the first part of the proposition concerning the necessary conditions for the existence of
$a_0,\ldots,a_{q-1}$, that is suppose (\ref{eq:n}) has solutions as described in the statement of the proposition.
We will first show that this part of the proposition holds for a particular class of fields $M$. We will consider
two cases: $\beta(x) \in F$ and $\beta(x) \not \in F$. In the first case let $M \supseteq K(x, \beta(x))$. If
$\beta(x) \not \in F$, then it is of degree $q$ over $F$ by Lemma \ref{le:extq}, and for each $i$ we can write
\[%
a_i =A_{i,0} + A_{i,1}\beta(x) + \ldots + A_{i,q-1}\beta(x)^{q-1},
\]%
where $A_{i,0}, \ldots, A_{i,q-1} \in F$.   In this case, let $M \supseteq K(x, A_{0,0}, \ldots,A_{q-1,q-1})$. Then
each $a_i \in M(\beta(x))$ and $[M(\beta(x)):M]=q$.

Now assume that either we are in the unramified case and for some factor $\pp_M$ of $\pp_K$, we have that $x$ is
not integral at $\pp_M$, or (\ref{eq:orderbound}) does not hold.

We begin with the unramified case. Let $\pp_M$ be an $M$-factor of $\pp_K$ such that $\ord_{\pp_M}x <0$. Then since
in the case under consideration we have that $\ord_{\pp_M}g= \ord_{\pp_K}g$, we conclude that  $\ord_{\pp_M}h <0$
and $\ord_{\pp_M}h \not \cong 0 \mod q$. Thus, $h^{-1} +z^q \cong z^q \mod \pp_M$. Hence if the extension
$M(\beta(x))/M$ is non-trivial, $\pp_M$ will split completely in this extension. Let $\pp_{M(\beta(x))}$ be any
factor of $\pp_M$ in $M(\beta(x))$. Given the arguments above, whether $M(\beta(x))= M$ or is a non-trivial
extension of $M$, we have that $X^q-b$ is irreducible over the residue field of $\pp_{M(\beta(x))}$. Since
$\pp_{M(\beta(x))}$ is prime to $q$ and does not occur in the divisor of $b$, the irreducibility of
$X^q-b$ over the residue field of $\pp_{M(\beta(x))}$ implies that $\pp_{M(\beta(x))}$ does not split in the
extension $M(\beta(x),\beta)/M(\beta(x))$. Further,
\[%
\ord_{\pp_{M(\beta(x))}}h=\ord_{\pp_M}h=\ord_{\pp_K}h \not \cong 0 \mod q.
\]%
Therefore, $h$ cannot be a $M(\beta(x))$-norm of an element from $M(\beta, \beta(x))$. But (\ref{eq:n}) asserts
precisely that. Consequently, we have a contradiction and conclude that in the unramified case, if for some $x \in
F$ it is the case that (\ref{eq:n}) has solutions $a_0,\ldots,a_{q-1} \in M(\beta(x))$, then $x$ is integral at
$\pp_K$.

We will now drop the assumption that $\pp_K$ has no ramified factors in any finite extension of $K$ contained
in $F$.  Define a  field $M$ as above and assume that  (\ref{eq:orderbound}) does not hold.  Then, given
our definition of $h$ and our assumption on the ramification degrees, we still have  $\ord_{\pp_M}h <0$ and
$\ord_{\pp_M}h \not \cong 0 \mod q$ as above.  Therefore from this point on we can proceed as before.\\

Now assume $M'$ is an arbitrary subfield of $F$ containing $K$. Let $x \in M'$ and suppose that equation
$(\ref{eq:n})$ has solutions $a_0,\ldots,a_{p-1} \in L_x=F(\beta(x))$ as described above. Let $M = M'(x,\beta(x))$
or $M=M'(x, A_{0,0}, \ldots,A_{q-1,q-1})$ as above depending on whether $\beta(x) \in F$. Let $\pp_{M'}$ be the
prime of $M'$ below $\pp_M$. Then by the arguments above, depending on whether we are in the ramified
or the unramified case, we have that either inequality (\ref{eq:orderbound}) holds or $\ord_{\pp_M}x >0$. Let $e$
the ramification degree of $\pp_M$ over $\pp_{M'}$. Then we also have that either
\[%
\ord_{\pp_M}x > \frac{(-q+1)\ord_{\pp_M}(g)}{r(s+1)} \Rightarrow \frac{1}{e}\ord_{\pp_M}x >
\frac{(-q+1)\ord_{\pp_M}(g)}{er(s+1)} \Rightarrow \ord_{\pp_{M'}}x > \frac{(-q+1)\ord_{\pp_{M'}}(g)}{r(s+1)},
\]%
or
\[%
\ord_{\pp_M}x >0 \Rightarrow \ord_{\pp_{M'}}x > 0.
\]%
Now the second assertion of the proposition follows directly from Theorem 2.1 of \cite{Sh23}.
\end{proof}%
Our next task is to reduce the number of assumptions on the field $F$ necessary for the definability of integrality
at a prime.  This will be accomplished in the two lemmas below: one for the general case and one for the unramified
case.  We treat the general case first.%

\begin{lemma}%
\label{le:exb}%
Let $L$ be an algebraic, possibly infinite extension of $\Q$. Let $Z$ be a number field contained in $L$ such that
$L$ is normal over $Z$.  Let $\pp_Z$ be a prime of $Z$ and assume that the following conditions are satisfied.
\begin{itemize}%
\item  There exists a non-negative integer $m_f$ such that any prime  lying  above $\pp_Z$ in a number field
contained in $L$ has a relative degree $f$ over $Z$ with $\ord_qf \leq m_f$.
\item There exists a non-negative integer $m_e$ such that any prime  lying above $\pp_Z$ in a number field
contained in $L$ has a ramification degree $e$ over $Z$ with $\ord_qe \leq m_e$.%
\end{itemize}%
Then there exists a finite extension $K$ of $Z$ such that $K$ and $F=KL$ satisfy the assumptions in
\ref{not:prime} for the general case with respect to all factors  $\pp_K$ of $\pp_Z$ in $K$.
\end{lemma}%
\begin{proof}%
  Let ${\mathbf M}_e$ be the set of all exponents $m$ such that $e=e_0q^m$ with $(e_0,q)=1$ is a ramification
degree over $Z$ for a number field prime lying above $\pp_Z$.  This is a set of non-negative integers bounded from
above and therefore must have a maximal element $\bar m$.  Let $U_0$ be a number field with a prime $\pp_{U_0}$
above $\pp_Z$ with  the  ramification degree over $Z$ divisible by $q^{\bar m}$.   Let $U_e$ be the Galois closure
of $U_0$ over $Z$ and observe that $U$ must also have  a prime $\pp_{U_e}$ above $\pp_Z$ with  the  ramification
degree over $Z$ divisible by $q^{\bar m}$, since $e(\pp_{U_e}/\pp_Z) = e(\pp_{U_e}/\pp_{U_0})e(\pp_{U_0}/\pp_Z)$.
Further, since $U_e$ is Galois over $Z$, it is the case that for any $U_e$-prime $\qq_{U_e}$ above $\pp_Z$ we have
that $\ord_{q}e(\qq_{U_e}/\pp_Z)=\bar{m}$.  Next we note that if $\bar{U}/U_e$ is a finite extension of number
fields with $\bar{U} \subset L$ and $\pp_{\bar{U}}$ is a prime above $\pp_{U_{e}}$, then $\ord_qe(\pp_{\bar
U}/\pp_{U_e})=0$.

Similarly we can find a field $U_f$, Galois over $Z$, so that for any finite extension $\hat{U}$ of $U_f$ and any
prime $\pp_{\hat{U}}$ lying above $\pp_Z$ in $\hat{U}$ we have that $\ord_qf(\pp_{\hat  U}/\pp_{U_f})=0$.     Let
$U=U_eU_f$ (observe that $U/Z$ is Galois),  let $K=U(\xi_q)$ and let $F=KL=L(\xi_q)$.  Note that $F/Z$ is also a
normal extension.   Let $N$ be a number field such that $K \subset N \subset F$.  Let $\pp_U$ be any prime lying
above above $\pp_Z$ in $U$, and let $\pp_N$ be any prime above $\pp_U$ in $N$.  Let $B \in R(\pp_Z)$ (the residue
field of $\pp_H$) be such that $B$ is not a $q$-th power.  We claim that $B$ is not a $q$-th power in $R(\pp_N)$
and $e(\pp_N/\pp_U) \not \equiv 0 \mod q$.

Indeed, since $N \subset L(\xi_q)$, for some field $T$ such that $U \subset T \subset L$ we have that $N \subset
T(\xi_q)$. Without loss of generality we can assume that $N = T(\xi_q)$. Let $\pp_T$ be a prime above $\pp_Z$ in
$T$. Then, by construction of $U$, we know that $e(\pp_T/\pp_U) \not \equiv 0 \mod q$ and  $f(\pp_T/\pp_U) \not
\equiv 0 \mod q$.  Next we observe that $e(\pp_N/\pp_T)$ and   $f(\pp_N/\pp_T)$ are both divisors  of $q-1$ and
therefore are not divisible by $q$, so that  $e(\pp_N/\pp_U) \not \equiv 0 \mod q$ and  $f(\pp_N/\pp_U) \not
\equiv 0 \mod q$. Consequently, we also have $e(\pp_N/\pp_K) \not \equiv 0 \mod q$ and $f(\pp_N/\pp_K) \not \equiv
0 \mod q$.  Further, since $([R(\pp_N):R(\pp_U)],q)=1$, we also conclude that $B$ is not a $q$-th power in
$R(\pp_N)$.
\end{proof}%

We now consider the unramified case.

\begin{lemma}%
\label{le:exb2}%
Let $L$ be an algebraic, possibly infinite extension of $\Q$. Let $Z$ be a number field contained in $L$ such that
$L$ is normal over $Z$.  Let $\pp_Z$ be a prime of $Z$ and assume that the following conditions are satisfied.
\begin{itemize}%
\item  There exists a non-negative integer $m_f$ such that any prime  lying  above $\pp_Z$ in a number field
contained in $L$ has a relative degree $f$ over $Z$ with $\ord_qf \leq m_f$.
\item There exists a non-negative integer $m_e$ such that any prime  lying above $\pp_Z$ in a number field
contained in $L$ has a ramification degree $e$ over $Z$ with $e \leq m_e$.%
\end{itemize}%
Then there exists a finite extension $K$ of $Z$ such that $K$ and $F=KL$ satisfy the assumptions in
\ref{not:prime} for the unramified case with respect to all factors  $\pp_K$ of $\pp_Z$ in $K$.
\end{lemma}%

\begin{proof}%
The proof of this lemma is almost identical to the proof of Lemma \ref{le:exb}. The only difference will come in
the way the field $U_e$ is selected. First we will let ${\mathbf M}_e$ be the set of all non-negative integers
$e$ such that $e$ is a ramification degree for a prime lying above $\pp_Z$ in some number field contained in $L$.
The set ${\mathbf M}_e$,  as in Lemma \ref{le:exb}, will have a maximal element $\bar e$. Let $U_e$ be a finite
Galois extension of $Z$  contained in $L$ such that for some $U_e$-prime $\pp_U$ lying above $\pp_Z$ the
ramification degree over $Z$ is $\bar{e}$. From this point on the proof proceeds as in Lemma \ref{le:exb}.
\end{proof}%

\begin{remark}%
\label{rem:also}%
 From the proof of the lemmas it is clear that if a field $Z \subset L$ and a $Z$-prime $\pp_Z$ satisfy the assumptions
of Lemmas \ref{le:exb} or \ref{le:exb2}, then any finite extension $T$ of $Z$ and any $T$-prime lying above $\pp_Z$ will
also satisfy the requirements of  Lemmas \ref{le:exb} or \ref{le:exb2} respectively.
\end{remark} %
Finally we state the main results of this section.  As above we separate out the general and the unramified case
for readability.  We start with a general case again.

\begin{theorem}%
\label{thm:ifmp}%
Let $L, Z, \pp_Z, q$ be as in Lemma \ref{le:exb}.  Then  there exist  an element $u \in L$  and  a subset $\calX$
of $L$ satisfying the following conditions:
\begin{itemize}%
\item If $x \in \calX$, then $ux$ is integral with respect to $\pp_Z$.%
\item If $x \in Z$ and $x$ is integral at $\pp_Z$, then $x \in \calX$.%
\item $\calX$ is Diophantine over $L$.%
\end{itemize}%
\end{theorem}%
\begin{proof}%
Let $F$ and $K$ be as in Lemma \ref{le:exb} . For each factor $\pp_K$ of $\pp_Z$ in $K$, let $g=g(\pp_K) \in K$ be
defined as in Notation and Assumptions \ref{not:prime} with respect to $\pp_K$. Then by Proposition \ref{prop:main},
and since intersection of Diophantine sets is Diophantine, there exists a set $\calY \subset F$ satisfying the
following conditions.
\begin{itemize}%
\item If $y \in \calY$, then $\left (\prod_{\pp_K| \pp_T}g(\pp_K)^q\right )y$ is integral with respect to $\pp_Z$.%
\item If $y \in K$ and $y$ is integral at every $\pp_K$ above $\pp_Z$, then $y \in \calY$.%
\item $\calY$ is Diophantine over $F$.%
\end{itemize}%
Let $\calX=\calY \cap L$. Then $\calX$ is Diophantine over $L$ by Proposition \ref{prop:rewriting}. Next  let $u
\in T$ be such that $\ord_{\pp_Z}u \geq \ord_{\pp_K}g(\pp_K)^q$ for all factors of $\pp_z$ in $K$.  Now observe
that if $x \in \calX$ then $x \in \calY$ and  $ux$ is  integral at $\pp_Z$. Conversely, if $x \in Z$ and is
integral at $\pp_Z$, then $x \in K$ and is integral at every factor of  $\pp_Z$ in $K$.  Thus $x \in \calY \cap
L=\calX$.
\end{proof}%

We now proceed to the unramified case.

\begin{theorem}%
\label{thm:ifmpu}%
Let $L, Z, \pp_Z, q$ be as in Lemma \ref{le:exb2}. Then there exists a subset $\calX$ of $L$ satisfying the following
conditions:
\begin{itemize}%
\item If $x \in \calX$, then $x$ is integral with respect to $\pp_Z$.%
\item If $x \in Z$ and $x$ is integral at $\pp_Z$, then $x \in \calX$.%
\item $\calX$ is Diophantine over $L$.%
\end{itemize}%
\end{theorem}%
\begin{proof}%
The proof is completely analogous to the proof of  Theorem \ref{thm:ifmp} but relies on Lemma \ref{le:exb2} and the
unramified case of Proposition \ref{prop:main}.
\end{proof}%
We finish this section with a corollary which we will use for the cases of infinite cyclotomic and abelian
equations.
\begin{corollary}%
\label{cor:deg}
Suppose $L$ is a normal, algebraic, possibly infinite extension of $\Q$ such that for some rational prime $q$ for
every number field $M$ contained in $L$ we have that $[M:\Q] \not \equiv 0 \mod q$. Then any number field $M
\subset L$ and any $M$-prime $\pp_M$ satisfy the assumptions of Theorem \ref{thm:ifmp} or Theorem \ref{thm:ifmpu}.
\end{corollary}%

\begin{proof}%
Let $M \subset L$ be a number field.  Let $M'$ be a finite extension of $M$.  Without loss of generality we can assume
that $M'$ is Galois over $\Q$.   Then for any prime $\qq_{M'}$ of $M'$ we have that
$e(\qq_{M'}/\qq_{\Q})f(\qq_{M'}/\qq_{\Q}) | [M':\Q]$, where $\qq_{\Q}$ is the rational prime below $\qq_{M'}$.  Thus,
$e(\qq_{M'}/\qq_{\Q})f(\qq_{M'}/\qq_{\Q}) \not \equiv 0\mod q$.  Let $\pp_{M}$ be any prime of $M$ and let $\pp_{M'}$ be
a prime of $M'$ above it.  Then
\[%
e(\pp_{M'}/\pp_{M})= \frac{e(\pp_{M'}/\pp_{\Q})}{e(\pp_{M}/\pp_{\Q})} \not \equiv 0 \mod q.
\]%
Similarly,
\[%
f(\pp_{M'}/\pp_{M})= \frac{f(\pp_{M'}/\pp_{\Q})}{f(\pp_{M}/\pp_{\Q})} \not \equiv 0 \mod q.
\]%
\end{proof}%

We now introduce the following notation.
\begin{notation}%
\label{not:I}%
Let $K$ be a number field, let $F$ be an algebraic possibly infinite extension of $K$, and let  $\calC_K$ be a
finite set of primes  of $K$.  Then let $I_{\calC_K/F}(x, t_1,\ldots,t_k) \in K[x,t_1,\ldots,t_k]$ be such that
\begin{equation}
\label{eq:I}
I_{\calC_K/F}(x, t_1,\ldots,t_k) =0
\end{equation}%
has solutions in $F$ only if $u_{\calC_K/F}x$ is integral at all the primes of $\calC_K$, where $u_{\calC_K/F} \in
\Z_{>0}$ is fixed and depends only $\calC_K$.  (As we have seen   $u_{\calC_K/F}$  can be equal to 1 in some
cases.)  Conversely, if $x \in K$ and is integral at all the primes of $\calC_K$, then (\ref{eq:I}) has solutions
in $K$.
\end{notation}%

\section{\bf Norm Equations over Totally Real Fields: an Update.}
\setcounter{equation}{0}
In this section we generalize the results from Section 3 of \cite{Sh26}.   The main reason for this generalization
is to allow for the treatment of arbitrary rings of $\calS$-integers of infinite abelian extensions of $\Q$. We
will point out the nature of the generalization below.
\begin{notationassumption}
\label{not:new}%
We start with a new set of assumptions and notation.
\begin{itemize}%

\item Let $M$ be a totally real number field of degree $n$ over $\Q$.%
\item Let $K$ be a subfield of $M$.
\item Let $L$ be a totally complex extension of degree 2 of $K$ such that $LM$ has no non-real roots of unity.%
\item Let $E_1/K, E_2/K$ be  totally real cyclic extensions of odd prime degrees $p_1$ and $p_2$
respectively with $(p_i,[M_G:\Q])=1$ for $i=1,2$, where $M_G$ is the Galois closure of $M$ over $\Q$.%
\item If $F$ is any number field, then let $U_F$ denote the group of its integral units and let $\calP(F)$ denote
the set of all non-archimedean primes of $F$.%
\item Let $N$ be any finite extension of a  number field $U$. Let $\calT_U$ (or $\calV_U$, $\calW_U$, $\calS_U$,
$\calE_U$, $\calN_U$, $\calL_U$, $\calR_U$, \ldots) be any set of primes of $U$. Then let $\calT_N$  (or
$\calV_N$, $\calW_N$, $\calS_N$, $\calE_N$, $\calN_N$, $\calL_N$, $\calR_N$, \ldots) be the set of primes of $N$
lying above the primes of $\calT_U$. Also let $\overline{\calT_N}$ be the closure of $\calT_N$ in $\calP(N)$ with
respect to conjugation over $\Q$.
\item Let $\calV_K$ be a set of primes of $K$ not splitting in either of the extensions $E_i/K$, $i=1,2$.%
\item Let $\calS_K=\{\pp_1,\ldots, \pp_s\} \subset \calP(K)$ be a set of primes not splitting in the extensions
$M/K$ and $E_2/K$.  Assume that at least one prime in $\calS_K$ splits completely in the extension $E_1L/K$.
\item Let $\calW_K=\calV_K \cup \calS_K$.%
\item Let $h_{LE_1}$ be a class number of $LE_1$.
\end{itemize}%
\end{notationassumption}%
Below we prove a generalization of Lemma 3.1 of \cite{Sh26}.  The main difference between the old and the new
versions of the lemma is that we replace a single $K$-prime $\pp$ splitting in the extension $LE_1/K$ and
remaining prime in the extension $M/K$ with a finite set of primes satisfying the same conditions.

\begin{lemma}%
\label{le:norms}%
Let $x \in O_{MLE_1, \calW_{MLE_1}}$ be a solution to the following system of equations.
\begin{equation}
\label{norm:1}
\left \{
\begin{array}{c}
{\bf N}_{MLE_1/ML}(x) =1\\
{\bf N}_{MLE_1/E_1M}(x) = 1
\end{array}
\right .
\end{equation}
Then $x^{2h_{LE_1}} \in  E_1L$.  Further such a solution exists.
\end{lemma}%
\begin{proof}%
First of all note that no prime of $\calW_M \setminus \calS_M$,  splits in the extension $ME_1/M$ by Lemmas
\ref{le:notsplit} and \ref{le:split}. Further, we note that given our assumptions on $\calS_K$, by Lemma
\ref{le:split}, all primes of $\calS_M$ split completely in the extension $MLE_1/M$.

Suppose now that $x \in LE_1M$ is a solution to the system of norm equations. Then the divisor of $x$ must be
composed of the primes lying above primes of $E_1M$ and $LM$ splitting in the extensions $LE_1M/E_1M$ and
$LE_1M/LM$ respectively. Given the fact that both extensions are cyclic of distinct prime degrees, we can conclude
that $LE_1M$-primes occurring in the divisor of $x$ lie above $M$-primes splitting completely in the extension
$LE_1M/M$. Thus, if $x \in O_{MLE_1,\calW_{MLE_1}}$ is a solution to the norm system, its divisor consists of
$MLE_1$-factors of primes in $\calS_M$ only. Further, since $LE_1M/E_1M$ is a totally complex extension of degree
2 of a totally real field, all the integral solutions to the second equation have to be roots of unity. Since
$ME_1L$ does not have any complex roots of unity, we can conclude the following. Let $x_1, x_2$ be two solutions
to the second norm equation above such that $x_1$ and $x_2$ have the same divisor.  Then on the one hand, $x_1 =
\pm x_2$. On the other hand, since primes of $\calS_K$ do not split in the extension $M/K$ and primes of $\calS_K$
split completely in the extension $LE_1/K$, we know that $LE_1$-factors of primes in $\calS_K$ do not split in the
extension $MLE_1/LE_1$ by Lemma \ref{le:notsplitabove}. Thus, there exists $y \in O_{LE_1,\calW_{LE_1}}$
such that $y$ has the same divisor as $x^{h_{LE_1}}$. Therefore, $y = \nu x^{h_{E_1L}},$ where $\nu$ is an
integral unit of $MLE_1$, and ${\bf N}_{ME_1L/E_1M}(y) = \mu$ is an integral unit of $E_1M$. On the other hand,
since $[E_1L:E_1]=[LE_1M:E_1M]=2$, we have that  ${\bf N}_{ME_1L/E_1M}(y) = {\bf N}_{E_1L/E_1}(y)$ and therefore
$\mu$ is an integral unit of $E_1$. Let $\bar{y} = \mu^{-1}y^{2}$. Then ${\bf N}_{ME_1L/E_1M}(\bar{y}) = {\bf
N}_{E_1L/E_1}(\bar{y}) = 1$. The divisors of $\bar{y}$ and $x^{2h_{E_1L}}$ are the same and therefore
$x^{2h_{E_1L}} \in E_1L$.

The proof of the fact that the system always has solutions in $O_{LE_1,\calW_{LE_1}}$ which are not roots of unity
remains the same as in Lemma 3.1 of \cite{Sh26}.

\end{proof}%

\section{\bf Norm Equations and Extensions of Degree 2 Real Fields.}%
\label{sec:norms}%
\setcounter{equation}{0}%
In this section we revisit a result of Denef and Lipshitz from \cite{Den2} and show that under some assumptions it
has a version that holds in the infinite extensions too. We will also prove a version of this result for ``large''
rings. We start with a notation and assumptions list and some facts about the primes and norm equations in the
extensions under consideration.
\begin{notationassumption}%
\label{not:norms}%
Below we  use Notation and Assumptions \ref{not:new} as well as the following  notation and assumptions.
\begin{itemize}%
\item Let $G$ be an extension of degree 2 of $K$ generated by $\alpha \in O_G$ with $\alpha^2=a \in O_K$.%
\item Let $H$ be an extension of degree 2 of $K$ generated by $\delta \in O_H$ such that $\delta^2=d \in O_K$.%
\item For all embeddings $\sigma : K \longrightarrow \C$, assume that $\sigma(d) >0$ if and only if $\sigma(a)
<0$. Further, assume that   $|\sigma(a)| >1$ for all $\sigma$'s  as above.%
\item For any number field $N$, let  $r_{N}$ be the number of real embeddings of $N$ into $\tilde{\Q}$, and let
$2s_{N}$ be the number of non-real embeddings of $N$ into $\tilde{\Q}$.%
\item Assume $s_G >0$.%
\item Let $\calZ_K$ be a set of primes of $K$ not splitting in the extension $E_2/K$.  We will assume that $\calZ_K
\supseteq \calW_K$.
\end{itemize}%
\end{notationassumption}

\begin{lemma}%
\label{le:rank1}
There exists $\varepsilon \in O_{HGE_2}$ such that $\varepsilon$ is not a root of unity, and
\begin{equation}%
\label{sys:norm}
\left \{%
\begin{array}{c}
{\mathbf N}_{HGE_2M/E_2GM}(\varepsilon)=1,\\%
{\mathbf N}_{GHE_2M/HGM}(\varepsilon)=1.%
\end{array}
\right .
\end{equation}
\end{lemma}%
\begin{proof}%

Since $[M:\Q]=n$, we have that  $[GM:\Q]=2n$, $[GHM:\Q]=4n$, $[HGME_2:\Q]=4p_2n$, and $[E_2GM:\Q]=2p_2n$. Next we
note that given a number field extension $U/F$, the integral solutions to the equation ${\mathbf N}_{U/F}(x)=1$ in
$U$ form a multiplicative group whose rank is equal to the difference of ranks between the integral unit groups of
$U$ and $F$. Let
\[%
A_{GM}=\{x \in O_{HGME_2}: {\mathbf N}_{HGME_2/GM}(x)=1\},
\]%
\[%
A_{GMH}=\{x \in O_{HGME_2}: {\mathbf N}_{HGME_2/GMH}(x)=1\},
\]%
\[%
A_{E_2GM}=\{x \in O_{HGME_2}: {\mathbf N}_{HGME_2/E_2GM}(x)=1\},
\]%
It is clear that $A_{GMH} \cup A_{E_2GM} \subseteq A_{GM}$.  By computing the ranks of the integral unit groups
involved we will show that
\begin{equation}%
\label{ineq:ranks}%
\rank A_{GMH} + \rank A_{E_2GM} > \rank A_{GM}.
\end{equation}%
This inequality implies that $A_{GMH} \cap A_{E_2GM}$ contains an element of infinite order. First of all, we have
that
\[%
r_{GM}+2s_{GM}=2n.
\]%
Further by Dirichlet Unit Theorem we know that
\[%
\rank U_{GM} = r_{GM}+s_{GM}-1.
\]%
Given our assumptions (see Notation and Assumptions \ref{not:norms}) on $HM$ and $GM$, every real embedding of $GM$
will extend to two non-real embeddings of $GMH$. (Non-real embeddings obviously always extend to non-real
embeddings.) Thus, $GMH$ will
have no real embeddings and $4n$ non-real embeddings into $\tilde{\Q}$. Therefore,%
\[%
\rank U_{GMH}=2n-1=r_{GM}+2s_{GM}-1.%
\]%
Since $E_2$ is a totally real field, $GME_2$ has $p_2r_{GM}$ real embeddings and $2p_2s_{GM}$ non-real embeddings
with
\[%
\rank U_{E_2GM}=p_2r_{GM} + p_2s_{GM} -1.
\]%
Finally, adjoining $E_2$ to $GMH$ will result in the field with $4p_2n$ non-real embeddings so that
\[%
\rank U_{GMHE_2} = 2p_2n-1= p_2r_{GM} + 2p_2s_{GM} -1.
\]%
To show that (\ref{ineq:ranks}) holds note the following:
\[%
\rank A_{GM} = \rank U_{HGME_2} - \rank U_{GM}= p_2r_{GM}+2p_2s_{GM} - r_{GM}-s_{GM},
\]%
\[%
\rank A_{GMH} = \rank U_{HGME_2} - \rank U_{GMH}= p_2r_{GM}+2p_2s_{GM} -r_{GM}-2s_{GM},
\]%
\[%
\rank A_{E_2GM} = \rank U_{HGME_2} - \rank U_{E_2GM}= p_2r_{GM}+2p_2s_{GM} - p_2r_{GM}-p_2s_{GM}=p_2s_{GM},
\]%
Thus,
\[%
\rank A_{GMH} + \rank A_{E_2GM} = p_2r_{GM} + 3p_2s_{GM} -r_{GM}-2s_{GM} >p_2r_{GM} +2p_2s_{GM} - r_{GM}- s_{GM} = \rank
A_{GM},
\]%
as long as $p_2s_{GM} >s_{GM}$.  This last inequality obviously  holds for any $p_2 >1$ since we assumed that
$s_{GM} \geq 1$.
\end{proof}%

The next lemma will state an easy result which will be crucial in eliminating the unwanted primes in the denominator.%

\begin{lemma}%
\label{le:nosplit}
The primes in $\calZ_{GMH}$ do not split in the extension $GMHE_2/GMH$.
\end{lemma}%
\begin{proof}%
This lemma follows from the fact that $([E_2:K],[GMH:K])=1$ and Lemma \ref{le:notsplit}.
\end{proof}%

\begin{corollary}%
\label{cor:m}%
Let $m$ be divisible by the size of the group of roots of unity in $GMHE_2$.  Then for any $\varepsilon \in
O_{GMHE_2,\calZ_{GMHE_2}}$ such that it  is a solution to (\ref{sys:norm}) we have that $\varepsilon^m \in
O_{MHE_2}$.  Further, if we assume that all the roots of unity in $GMHE_2$ are already in $GM$, we can replace $m$
by 2.
\end{corollary}%
\begin{proof}%
First, since $\varepsilon \in  O_{GMHE_2,\calW_{GMHE_2}}$, the only primes which can occur in the denominator of the
divisor of $\varepsilon$ are primes from $\calW_{GMHE_2}$.  Secondly, since ${\mathbf
N}_{GME_2H/HGM}(\varepsilon)=1$, the only primes which can occur in the numerator of the divisor of $\varepsilon$
are the primes that have a distinct conjugate over $GMH$ which is allowed in the denominator of the divisors of
the elements of  $O_{GMHE_2,\calW_{GMHE_2}}$.    But by Lemma \ref{le:nosplit}, no primes of $\calW_{GMHE_2} $ has a
distinct conjugate over $GMH$.  Consequently, $\varepsilon$ has a a trivial divisor and therefore is an integral
unit.

Second, let
\[%
A_{E_2M}=\{x \in O_{HE_2M}: {\mathbf N}_{HE_2M/E_2M}(x)=1\},
\]%
and note that  $A_{E_2M} \subseteq A_{E_2GM}$ since $GM$ and $HM$ are linearly disjoint  over $E_2M$ by assumptions
in \ref{not:norms}. Using notation from Lemma \ref{le:rank1}, we have that $\rank A_{E_2GM}=p_2s_{GM}$.  To
compute the
rank of $A_{E_2M}$ we need to compute the ranks of integral unit groups of $HE_2M$ and $E_2M$.  It is easy to see that%
\[%
\rank U_{E_2M}=p_2n-1= p_2\frac{r_{GM}+2s_{GM}}{2}-1=\frac{pr_{GM}}{2}+ps_{GM}-1.%
\]%
To compute, $U_{HE_2M}$ we can look at the number of real and non-real embeddings of $HM$ first and then multiply
these numbers by $p_2$ to get the analogous information for $E_2HM$. From the assumptions in \ref{not:norms} we
have that  $HM$ has $2s_{GM}$ real and $r_{GM}$ non-real embeddings into $\tilde{\Q}$. Thus $E_2HM$ has
$2p_2s_{GM}$ real and
$p_2r_{GM}$ non-real embeddings. Therefore,%
\[%
\rank U_{HE_2M}= 2p_2s_{GM}+\frac{p_2r_{GM}}{2}-1.%
\]%
Hence, $\rank A_{E_2M}=p_2s_{GM}=\rank A_{E_2GM}$. Now suppose $\varepsilon \in A_{E_2GM}$. Since the ranks are the
same and $A_{E_2M} \subseteq A_{E_2GM}$, we conclude that for some positive $l \in \N$ we have that $\varepsilon^l
\in O_{HE_2M}$. Let $\varepsilon'$ be the conjugate of $\varepsilon$ over $E_2HM$. Then,
$(\varepsilon/\varepsilon')^l=1$ or, in other words, $\varepsilon/\varepsilon'=\xi$ -- a root of unity in $GMHE_2$
with ${\mathbf N}_{GMHE_2/GME_2}(\xi)=1$. Thus $l|m$. If we now assume that $\xi \in GM$, we conclude that ${\mathbf
N}_{GMHE_2/GME_2}(\xi)=\xi^2=1$. Thus, $\varepsilon'=\pm\varepsilon$ and therefore $\varepsilon^2 \in HE_2M$.
\end{proof}%

\section{\bf Bounds for Extensions of Degree 2.}%
\label{sec:bound2}
\setcounter{equation}{0}%

\begin{notation}
\label{not:bounds}%
We now add to the list of Notation and Assumptions \ref{not:new} and \ref{not:norms}.
\begin{itemize}%
\item Assume that  an integral unit $\mu$ generates  $E_2$ over $K$. Denote the monic
irreducible polynomial of $\mu$ over $K$ by $P(X)$, and assume that $P(X) \in \Z[X]$. (This assumption implies
that $\mu$ is of degree $p_2$ over $\Q$, $K$, $M$ and $GM$.)%
\item Let $\calM_{GM} \subset \calZ_{GM}$ be a set of $GM$-primes lying above $M$-primes not splitting in the
extension $GM/M$.%
\item Let $\calU_{GM} =\calM_{GM} \cup \calS_{GM}$.%
\item Let $\calE_K$ contain all the primes $\pp$ of $\calZ_K$ such that $\pp$ divides the discriminant of $P(X)$
and all the primes of $\calS_K$.   Given our definition of $\calZ_K$, we have that $\calE_K$ is a finite set of
$K$-primes. (In the future we might add primes to this set, but it will always remain finite.)  Let $N_{\calE}$ be
a positive integer divisible by all the primes of $\calE_K$. \item Let $Q(X) = P(N_{\calE_K}X)$.
\end{itemize}%
\end{notation}%
Eventually we will use $P(X)$ to get away from the factors of primes in $\calZ_K$ in the denominator.  We know that
for all primes  $\pp \in \calZ_K$ not dividing the discriminant of $P(X)$, for all $x \in K$, it is the case that
$\ord_{\pp}P(x) \leq 0$.  However we need to take care of the finitely many extra primes in $\calE_K$ possibly
dividing the discriminant of $P(X)$ or inconvenient in some other ways (as will be explained later).  To that
effect we adjust $P(X)$.
\begin{lemma}%
\label{le:allprimes}%
Let $x \in GM$ and assume that $x$ is integral at all the primes of $\overline{\calE}_{GM}$.   Then for all $\pp
\in \overline{\calZ}_{GM}$ we have that $\ord_{\pp}P(N_{\calE}x) \leq 0$.  Further, for all $\pp \in
\overline{\calE}_{GM}$ we have that $\ord_{\pp}P(N_{\calE}x) = 0$.
\end{lemma}%

\begin{proof}%
Let $C$ be the Galois closure of $GM$ over $\Q$.  It is enough to show that the lemma holds for $C$ in place of
$GM$ and for $\calZ_C$  and $\calE_C$ in place of $\calZ_{GM}$ and $\calE_{GM}$ respectively.  First  observe that
$[C:\Q]= [M_G:\Q]2^j$, for some $j \in \Z_{>0}$.  Thus, given our assumption that $p_2$ is odd and
$([M_G:\Q],p_2)=1$, we conclude that $([C:K], p_2)=1$.   Second, by Lemma \ref{le:notsplitabove}, no prime of
$\calZ_{C}$ will split in the extension $CE_2/C$.  Suppose for some $\pp \in \calZ_{C} \setminus \calE_{C}$, some
$x \in C$  we have that $\ord_{\pp}Q(X) >0$.  Then $P(X)$ has a root modulo $\pp$.  But then $\pp$ has a relative
degree one factor in the extension $CE_2/C$ (see \cite{L}, page 25), contradicting our arguments above.  Suppose
now that $\pp \in \calE_C$, $x \in C$ is integral at $\pp$ and $Q(x)= P(N_{\calE}x) \cong 0 \mod \pp$.  But given
our assumption on $x$ and $N_{\calE}$ we have that $N_{\calE}x \cong 0 \mod \pp$ and the free term of $P(X)$ is an
integral unit.  Hence we have a contradiction. Finally, suppose  that $\bar{\qq}$ is a conjugate of $\qq \in
\calZ_{C}$ over $\Q$, and for some $\bar{x} \in C$ we have that $\ord_{\bar{\qq}}Q(\bar{x}) >0$.  Then, since
$Q(T) \in \Z[T]$, we have that $\ord_{\qq}Q(x) >0$.

Suppose now that $\pp \in \overline{\calE}_C$.  By the argument above we have that   $\ord_{\pp}P(N_{\calE}x) \leq
0$.  On the other hand, $ N_{\calE}x$ is, by assumption, integral at $\pp$, and $P(X) \in \Z[X]$.  Thus,
$P(N_{\calE}x)$ is also integral at $\pp$.  Consequently, $P(N_{\calE}x)$ is a unit at $\pp$.
\end{proof}%

\begin{notationassumption}%
\label{not:Q}
Here we make additions to our notation set.
\begin{itemize}%

\item Let $\beta_1,\ldots,\beta_{p_2}$ be all the roots of the polynomial $P(X)$.%
\end{itemize}
\end{notationassumption}
We continue with a series of lemmas often used to obtain bounds on non-archimedean valuations.
\begin{lemma}%
\label{le:primesofG}%
Primes of $\calM_{GME_2}$ lie above primes of $M$ not splitting in the extension $E_2GM/M$.
\end{lemma}%
\begin{proof}%
Since $([M:K], [GE_2:K])=1$, the assertion of the lemma follows from Lemma \ref{le:notsplit}.
\end{proof}%
\begin{lemma}%
\label{le:order}%
Let $x \in O_{GME_2,\calU_{GME_2}}$, $x =y_0+y_1\alpha \equiv 0 \mod \mte{Z}$ in $O_{GME_2,\calU_{GME_2}}$, where
$y_0,y_1 \in E_2M$, $\mte{Z}$ is an integral divisor of $E_2M$ without any factors in $\calU_{GME_2}$.   (We
remind the reader that $\alpha$ which has been defined in Notation and Assumptions \ref{not:norms} is an integral
generator of $GM$ over $M$.) Assume additionally that for any $\mte{t} \in \overline{\calU}_{GME_2}$ we have that
$\ord_{\mte{t}}x \leq 0$ and $x$ is a unit at all the primes of $\overline{\calS}_{GME_2}$.  Let ${\mathbf
N}_{E_2GM/\Q}(x)=\frac{X}{Y}$, where $X, Y \in \Z$ and $(X,Y)=1$ in $\Z$. Let $Z= {\mathbf
N}_{E_2GM/\Q}(\mte{Z})$. Then $\frac{Y}{Z}{\mathbf N}_{GME_2/\Q}(2\alpha y_1)$ is an integer.
\end{lemma}%

\begin{proof}%
Let $\bar{x} =y_0-y_1\alpha$ be the conjugate of $x$ over $E_2M$. By assumption on $\calS_{GME_2}$ we have that
$\ord_{\pp}x = 0$ for all $\pp \in \calS_{GME_2}$. Also, by assumption on $\calM_{GME_2}$, for every $\mte{t} \in
\calM_{GME_2}$ we have that
\[%
\ord_{\mte{t}}x=\ord_{\mte{t}}\bar{x}.
\]%
Thus,  we have that
\[%
\ord_{\mte{t}}2\alpha y_1 =\ord_{\mte{t}}(x-\bar{x}) \geq \ord_{\mte{t}}x.
\]%
Since $x$ does not have positive order at any prime of $\overline{\calU}_{GME_2}$, the last inequality also asserts
that
\[%
|\ord_{\mte{t}}2\alpha y_1| \leq |\ord_{\mte{t}}x|.
\]%

Let $\frac{\mathfrak B}{\mathfrak C}$ be the $GME_2$-divisor of $2\alpha y_1$ with $\mte{B}, \mte{C}$ being
relatively prime integral divisors such that all the factors of  $\mte{C}$ are in $\calM_{GME}$.   Let ${\mathbf
N}_{GME_2/\Q}(2\alpha y_1)=\frac{B}{C}$, where $B, C \in \Z$ and are relatively prime in $\Z$.  Then  $C \Big{|}
{\mathbf N}_{GME_2/\Q}(\mte{C})$ in $\Z$.  Next, let $\frac{\mte{X}}{\mte{Y}}$, where $\mte{X}, \mte{Y}$ are
relatively prime integral divisors with all the factors of $\mte{Y}$ in $\calM_{GME_2}$ and no factor of $\mte{X}$
is in $\overline{\calM}_{GME_2}$,  be the $GME_2$-divisor of $x$.    Given our assumptions on $\mte{X}$ and
$\mte{Y}$, we can conclude that ${\mathbf N}_{GME/\Q}(\mte{X})$ and ${\mathbf N}_{GME/\Q}(\mte{Y})$ are relatively
prime and therefore the divisor of $Y$ is ${\mathbf N}_{GME_2/\Q}(\mte{Y})$.  We claim that ${\mathbf
N}_{GME_2/\Q}(\mte{C}) \Big{|} {\mathbf N}_{GME_2/\Q}(\mte{Y})$.  Let
\begin{equation}
\label{eq:divisors}
\mte{C}=\prod_{\mte{p} \in \calU_{GME_2}}\mte{p}^{a(\mte{p})},
\end{equation}%
\begin{equation}%
\label{eq:divisors2}
\mte{Y}=\prod_{\mte{p} \in \calU_{GME_2}}\mte{p}^{b(\mte{p})},
\end{equation}%
where $a(\mte{p})=b(\mte{p})=0$ for all but finitely many $\mte{p}$.  Further, by the argument above, we also have
that $b(\mte{p}) \geq a(\mte{p})$ for all $\mte{p} \in \calU_{GME_2}$.  Using (\ref{eq:divisors}) and
(\ref{eq:divisors2}),  we can write
\[%
{\mathbf N}_{GME_2/\Q}(\mte{C})=\prod_{\mte{p} \in \calU_{GME_2}}p(\mte{p})^{a(\mte{p})f(\mte{p})},
\]%
\[%
Y={\mathbf N}_{GME_2/\Q}(\mte{Y})=\prod_{\mte{p} \in \calU_{GME_2}}p(\mte{p})^{b(\mte{p})f(\mte{p})},
\]%
where $p(\mte{p})$ is the rational prime below a $GME_2$-prime $\mte{p}$, and $f(\mte{p})$ is the relative degree of
$\mte{p}$ over $\Q$.  Now assertion follows from the fact that $b(\mte{p}) \geq a(\mte{p})$.

Next we note that for any $\mte{q}$ such that $\ord_{\mte{q}}\mte{Z} >0$ we have that $\ord_{\mte{q}}2\alpha y_1
\geq\ord_{\mte{q}}\mte{Z}$. This follows from the fact that $x \cong \bar{x} \equiv 0 \mod \mte{Z}$, when we
consider $\mte{Z}$ as an ideal of $O_{GME_2,\calM_{GME_2}}$. Since no prime factor of $\mte{Z}$ is allowed in the
denominator of the elements of our ring, we have that  $Z \Big{|} B$ in $\Z$, and
the lemma holds.
\end{proof}%
In the next lemma we remove the assumption that the primes allowed in the denominator stay prime in the extension
$GME_2/GM$.
\begin{lemma}%
\label{le:order2}%
Let $x \in O_{GME_2,\calZ_{GME_2}}, x =y_0+y_1\alpha \equiv 0 \mod \mte{Z}$ in $O_{GME_2,\calW_{GME_2}}$, where
$\mte{Z}$ is an integral divisor, $\ord_{\mte{t}}\mte{Z}=0$ for all $\mte{t} \in \calZ_{GME_2}$.  Assume
additionally that for any $\ttt \in \overline{\calZ}_{GME_2}$ we have that $\ord_{\mte{t}}x \leq 0$.  Let ${\mathbf
N}_{E_2GM/\Q}(x)=\frac{X}{Y}$, where $X, Y \in \N$ and $(X,Y)=1$ in $\Z$. Let $Z= {\mathbf N}_{E_2GM/\Q}(\mte{Z})$.
Then $\frac{Y^2}{Z}{\mathbf N}_{GME_2/\Q}(2\alpha y_1)$ is an integer.
\end{lemma}%
\begin{proof}%

First of all, as above, we have that  $2\alpha y_1 = x-\bar{x}$, where $\bar{x}$ is the conjugate of $x$ over
$E_2$. Therefore, given our assumptions, for any prime $\mte{t}\in \calZ_{GME_2}$,  if $\ord_{\mte{t}}(2\alpha y_1)
<0$, then $\ord_{\mte{t}}x <0$ or $\ord_{\mte{t}}\bar{x} <0$.  Thus, we need to consider three cases:
\[%
r(\ttt)=\ord_{\ttt}x=\ord_{\mte{t}}\bar{x},
\]%
\[%
r_1(\ttt)=\ord_{\mte{t}}x <r_2(\ttt)=\ord_{\mte{t}}\bar{x},
\]%
\[%
r_2(\ttt)=\ord_{\mte{t}}\bar{x} <r_1(\ttt)=\ord_{\mte{t}}x.
\]%
In the first case, $\ord_{\mte{t}}2\alpha y_1\geq r(\ttt)$.   In the second case,
\[%
\ord_{\mte{t}}2\alpha y_1=\ord_{\mte{t}}x = r_1,
\]%
and in the third case
\[%
\ord_{\mte{t}}2\alpha y_1=\ord_{\mte{t}}\bar{x} = r_2.
\]%
Next let
\[%
\calY_0= \{ \mte{p} \in \calP(GME_2): \ord_{\mte{p}}x = \ord_{\mte{p}}\bar{x} <0\}
\]%
\[%
\calY_1= \{ \mte{p} \in \calP(GME_2): \ord_{\mte{p}}x  < \ord_{\mte{p}}\bar{x}\leq 0\}
\]%
\[%
\calY_2= \{ \mte{p} \in \calP(GME_2): \ord_{\mte{p}}\bar{x}  < \ord_{\mte{p}}x\leq 0\}
\]%
Observe that since $\ord_{\pp}x= \ord_{\bar{\pp}}\bar{x}$, where $\pp$ and $\bar{\pp}$ are primes conjugate over
$E_2M$,  we have that $\calY_2$ consists of the  conjugates over $E_2$ of primes in  $\calY_1$, and
$r_1(\pp)=r_2(\bar{\pp})$. Further $\calY_0$ is closed under conjugation. As above, let $\frac{\mte{X}}{\mte{Y}}$,
$\frac{\overline{\mte{X}}}{\overline{\mte{Y}}}$, where $\mte{X}, \mte{Y}, \overline{\mte{X}}, \overline{\mte{Y}}$
are integral divisors and $(\mte{X}, \mte{Y})=(\overline{\mte{X}},\overline{\mte{Y}})=1$, be the $GME_2$-divisors
of $x$ and $\bar x$ respectively, and let $\frac{\mte{B}}{\mte{C}}$, where $\mte{B}, \mte{C}$ are integral
relatively prime divisors, be the $GME_2$-divisor of $2\alpha y_1$.  Then we can write
\[%
\mte{Y}=\prod_{\pp \in \calY_0}\mte{p}^{r(\pp)}   \prod_{\pp \in \calY_1 \cup \calY_2}\mte{p}^{r_1(\pp)},
\]%
and
\[%
\bar{\mte{Y}}=\prod_{\bar{\pp} \in \calY_0}\bar{\mte{p}}^{r(\bar{\pp})}   \prod_{\bar{\pp} \in \calY_1 \cup
\calY_2}\bar{\mte{p}}^{r_1(\bar{\pp})}= \prod_{\pp \in \calY_0}\mte{p}^{r(\pp)}   \prod_{\pp \in \calY_1 \cup
\calY_2}\mte{p}^{r_2(\pp)}.
\]%

Consequently,
\[%
{\mathbf N}_{GME_2/\Q}(\mte{Y}) =   {\mathbf N}_{E_2M/\Q}({\mathbf N}_{GME_2/E_2M}\mte{Y})= {\mathbf
N}_{E_2M/\Q}(\mte{Y}\overline{\mte{Y}})={\mathbf N}_{E_2M/\Q}\left (\prod_{\pp \in \calY_0}\mte{p}^{2r(\pp)}
\prod_{\pp \in \calY_1 \cup \calY_2}\mte{p}^{r_1(\pp)+r_2(\pp)} \right ).
\]%
Next we note that
\[%
\mte{C}=\prod_{\pp \in \calY_0}\pp^{r_0(\pp)}\prod_{\pp \in \calY_1}\pp^{r_1(\pp)}\prod_{\pp \in
\calY_2}\pp^{r_2(\pp)},
\]%
where $r_0(\pp) \leq r(\pp)$.  Since the conjugate of $2\alpha y_1$ over $ME_2$ is $-2\alpha y_1$, we have that
$\bar{\mte{C}}=\mte{C}$ and thus
\[%
{\mathbf N}_{GME_2/\Q}\mte{C}= {\mathbf N}_{E_2M/\Q}({\mathbf N}_{GME_2/E_2M}(\mte{C})) ={\mathbf
N}_{E_2M/\Q}(\mte{C}\bar{\mte{C}})={\mathbf N}_{E_2M/\Q} \left (\prod_{\pp \in \calY_0}\pp^{2r_0(\pp)}\prod_{\pp
\in \calY_1}\pp^{2r_1(\pp)}\prod_{\pp \in \calY_2}\pp^{2r_2(\pp)} \right ),
\]%
Now it is clear that ${\mathbf N}_{GME_2/\Q}(\mte{C}) \Big{|} {\mathbf N}_{GME_2/\Q}(\mte{Y})^2$.

Since on the one hand $x$ does not have positive order at any prime of $\overline{\calZ}_{GME_2}$, as in Lemma
\ref{le:order}, we can conclude that ${\mathbf N}_{GME_2/\Q}(\mte{X})$ and ${\mathbf N}_{GME_2/\Q}(\mte{Y})$ are
relatively prime as $\Q$-divisors. Thus, ${\mathbf N}_{GME_2/\Q}(\mte{Y})$ is the divisor of $Y$. On the other
hand, if we let $\frac{B}{C}= {\mathbf N}_{GME_2/\Q}(2\alpha y_1)$, where $B, C \in \Z$ and are relatively prime
in $\Z$, then certainly $C \Big{|} {\mathbf N}_{GME_2/\Q}(\mte{C})$. Thus, $\frac{Y^2}{Z}{\mathbf
N}_{GME_2/\Q}(2\alpha y_1)$ will have no rational primes lying below the primes of $\calZ_{GME_2}$ in the
denominator. Further, by assumption, as in Lemma \ref{le:order} we have that for all primes $\mte{t}$ such that
$\ord_{\mte{t}}\mte{Z}>0$ it is the case  that $\ord_{\mte{t}}x \geq\ord_{\mte{t}}\mte{Z}$ and
$\ord_{\mte{t}}\bar{x} \geq\ord_{\mte{t}}\mte{Z}$. Therefore, $\ord_{\mte{t}}(x-\bar{x}) \geq
\ord_{\mte{t}}\mte{Z}$ and consequently,
\[%
\ord_{\mte{t}}2\alpha y_1=\ord_{\mte{t}}(x-\bar{x})  \geq \ord_{\mte{t}}\mte{Z}.
\]%

\end{proof}%

The next lemma follows from the fact that $Q(X) \in \Z[X]$ has a positive leading coefficient.
\begin{lemma}%
\label{le:frombelow}%
 For any positive integer $k$ there exists a  number $A(k)>1$ such that for any real $x >A(k)$, we have that $Q(x) > k$.
\end{lemma}%
\begin{notation}%
\label{not:B}
For future use we introduce the following notation.
\begin{itemize}
\item Let $B=A(2)+1$.
\end{itemize}
\end{notation}%
The following lemma will allow us to establish some bounds on coordinates in a degree 2 extension of a totally real
field.
\begin{lemma}%
\label{le:boundtwo}%
Let $y \in E_2GM$ and let $x = y_0 + y_1\alpha, y_0,y_1 \in E_2M$.    Let $z \in E_2GM$ and suppose that for every
$\sigma_1,\ldots,\sigma_{pr_{GM}}:E_2GM \longrightarrow \tilde{\Q} \cap \R$, we have that
\begin{equation}%
\label{eq:y>1}
1\leq |\sigma_i(x)| < |\sigma_i(z)|,
\end{equation}%
while for all $\tau_1,\ldots,\tau_{ps_{GM}} : E_2GM\longrightarrow \tilde{\Q} $ with $\tau_i(E_2GM) \not \subseteq
\R$, we have that
\begin{equation}%
\label{eq:z>1}
|\tau_i(z) | \geq 1
\end{equation}
 Then $|{\mathbf N}_{GME_2/\Q}(y_1) | \leq |{\mathbf N}_{GME_2/\Q}(z) {\mathbf N}_{GME_2/\Q}(x)|$.
\end{lemma}%
\begin{proof}%
First of all observe that for all non-real embeddings $\tau_i : E_2GM \rightarrow \tilde{\Q}$, we have that
$\overline{\tau_i(x)}=\tau_i(y_0) - \tau_i(\alpha)\tau_i(y_1)$ and therefore,
\begin{equation}
\label{eq:nonreal}
|\tau_i(y_1)| =\left |\frac{\tau_i(x)-\overline{\tau_i(x)}}{2\tau_i(\alpha)}\right|\leq \left
|\frac{\tau_i(x)}{\tau_i(\alpha)}\right |< |\tau_i(x)|,
\end{equation}%
since by assumption $|\tau_i(\alpha)| >1$. On the other hand, for any real embedding $\sigma_i : E_2GM \rightarrow
\tilde{\Q}$, we have that $\sigma_i(x) = \sigma_i(y_0) +\sigma_i(\alpha)\sigma_i(y_1)$, while for some $i' \in
\{1,\ldots, pr_{GM}\}$ we also have that $\sigma_{i'}(x) = \sigma_i(y_0) -\sigma_i(\alpha)\sigma_i(y_1)$.  Thus,
\begin{equation}%
\label{eq:real}%
|\sigma_i(y_1)| =\left |\frac{\sigma_i(x)-\sigma_{i'}(x)}{2\sigma_i(\alpha)}\right|\leq \left
|\frac{2\sigma_i(z)}{2\sigma_i(\alpha)}\right |< |\sigma_i(z)|,
\end{equation}%
Putting together (\ref{eq:nonreal}) and (\ref{eq:real}) we obtain,
\begin{equation}%
\label{eq:together}%
|{\mathbf N}_{E_2GM/\Q}(y_1) | \leq \prod_{\tau_i}|\tau_i(x)|\prod_{\sigma_j}|\sigma_j(z)|\leq |{\mathbf
N}_{GME_2/\Q}(xz)|,
\end{equation}%
where the last  inequality follows from (\ref{eq:y>1}) and (\ref{eq:z>1}).
\end{proof}%

\begin{lemma}%
\label{le:n>1}%
Let $x \in GME_2$. Let $l >2$ be an integer such that $l >\max(\{|\beta_j|, j=1,\ldots,p_2\})$. Let
$x_k=Q(x-8(k+1)l)$. Then for some value of $k \in \{0,\ldots, 2p_2n\}$ we have that for any embedding $\phi: GME_2
\rightarrow \tilde \Q$, it is the case that $|\phi(x_k)| > 2$.  (Here we remind the reader that from Notation and
Assumptions \ref{not:norms} we have that $n=[M:\Q]$, $[E_2M:M]=p_2$, and $[E_2GM:E_2M]=2$.)
\end{lemma}%
\begin{proof}%
 In $GME_2$ we can factor
\[%
Q(x-8(k+1)l)=\prod_{j=1}^{p_2}(N_{\calE}x-8(k+1)l- \beta_j).
\]%
Let $B_k \subset \C$ be the closed ball of radius $2l$ centered at $8(k+1)l$ (i.e. $B_k = \{ z \in \C:
|z-(8(k+1)l)| \leq 2l\}$).  We claim that for all $j, k$ it is the case that $\beta_j + 8(k+1)l \in B_k$.
Further,  the distance between any point of $B_k$ and any  point of  $B_{k'}$ for $k\not = k'$ is at least $4l\geq
4$ and so $B_k \cap B_{k'}=\emptyset$ for $k \not = k'$.

Let $\phi_1,\ldots, \phi_{2p_2n}$ be all the embeddings  of $E_2GM$ over $\Q$. Then each $\phi_i(N_{\calE}x)$ can
be located in at most one $B_k$.  Further, if $\phi_i(N_{\calE}x) \not \in B_k$ then for any $j$ we have that
$|\phi_i(N_{\calE}x) - \beta_j-8(k+1)l| >l > 2$.  Since there is at least one $B_k$ without any
$\phi_i(N_{\calE}x)$'s, for some $k=\{0,\ldots, 2p_2n\}$, for all $i,j$, we have that
\[%
|\phi_i(N_{\calE}x)-8l(k+1)-\beta_j| > 2
\]%
implying that for all $i$ we have that
\[%
\phi_i(Q(x-8(k+1)l)) =\phi_i(P(N_{\calE}x-8(k+1)l))=\prod_{j=1}^{p}(\phi_i(N_{\calE}x)-8(k+1)l-\beta_j) \geq 2.
\]%
\end{proof}%


\section{\bf Diophantine Definability for Extensions of Degree 2.}%
\label{sec:deg2}
\setcounter{equation}{0}%

In this section we consider two versions of Diophantine definability for extensions of degree 2 of totally real
fields. In the first version we will restrict ourselves to the ring $O_{GME_2,\calU_{GME_2}}$, but the definition
will not use the degree of $GM$ over $\Q$. In the second  version we will use explicitly the degree of $GM$ over
$\Q$ but allow any prime of $\calZ_{E_2GM}$ in the denominator.
\begin{notation}%
\label{not:degree2}
We add the following to our notation and assumption list.
\begin{itemize}%
\item Let $l$ be as in Lemma \ref{le:n>1}, i.e. let $l \in \Z_{>0}$ be an upper bound for the absolute values of
all the roots of $P(X)$.
\end{itemize}%
\end{notation}
We start with a technical lemma.
\begin{lemma}%
\label{le:closeto1} Let $\varepsilon \in O_{HGME_2}$ be a solution to (\ref{sys:norm}).  Then for any positive
integer $k$ and any $\lambda >0$ there exists a positive integer $r$ such that for all $\tau : GME_2H  \rightarrow
\tilde{\Q}$ with $\tau(E_2HM) \not \subset \R$ we have that
\[%
\left|\frac{\tau(\varepsilon^{rk}-1)}{\tau(\varepsilon^r-1)}-k\right|< \lambda.
\]%
\end{lemma}%
\begin{proof}%
We start with an elementary observation:
\[%
T^k-1=(T-1+1)^k-1= \sum_{i=1}^{k} \frac{k!}{i!(k-i)!}(T-1)^{i}=\sum_{i=1}^{k} \left (\begin{array}{c} k\\i
\end{array} \right )(T-1)^i,
\]%
and so
\[%
\frac{T^k-1}{T-1}=\sum_{i=1}^{k} \frac{k!}{i!(k-i)!}(T-1)^{i-1}=   \sum_{i=1}^{k} \left (\begin{array}{c} k\\i \end{array}
\right )(T-1)^{i-1}
\]%
Assume that $k$ and $0<\lambda<1$ are fixed and suppose $z \in \C$ is such that $|z-1| <2^{-k}\lambda$.  Then
\[%
\left |\frac{z^k-1}{z-1}- k\right | = \left |  \sum_{i=2}^{k} \left (\begin{array}{c} k\\i \end{array}\right
)(z-1)^{i-1} \right | \leq 2^{-k}\nu\sum_{i=0}^{k} \left (\begin{array}{c} k\\i \end{array}\right ) =\lambda.
\]%

Now let $\varepsilon \in O_{GME_2H}$ be a solution to (\ref{sys:norm}).  Let $m$ be defined as in  Corollary
\ref{cor:m} and deduce that   $\varepsilon^m \in E_2MH$ with $|\tau(\varepsilon^m)| =1$ for all $\tau$, non-real
embeddings of $E_2MH$ into $\tilde{\Q}$.  Thus, given $\nu>0$ and $k$ the problem reduces to showing that for some
power $r \cong 0 \mod m$ of $\varepsilon$ we will have $|\tau(\varepsilon)^r-1| < 2^{-k}\nu$ for all non-real
embeddings of $E_2HM$ into $\tilde{\Q}$.  The proof of this fact is completely analogous of the proof of Lemma 12
of \cite{Sh2}.
\end{proof}%

\begin{lemma}%
\label{le:inf2} %
Let
\[%
x_0, x_1 \in O_{GM,\calU_{GM}},
\]%
\[%
a_1, a_2, b_1, b_2, c, d, u, v \in O_{GM,\calU_{GM}}[\mu] \subset O_{E_2GM, \calU_{GME_2}},
\]%
\[%
\varepsilon_i, \gamma_i \in O_{GM,\calU_{GM}}[\mu, \delta] \subset O_{E_2GMH, \calU_{E_2GMH}}, i=1,\ldots,4.
\]%
Assume also that the following conditions and equations are satisfied.
\begin{equation}
\label{eq:supplement}
\forall \pp \in \overline{\calE}_{GM}, \,\, \ord_{\pp}x_0 \geq 0
\end{equation}
\begin{equation}%
\label{eqdeg2:1}
x_1 = Q(x_0),
\end{equation}

\begin{equation}%
\label{eqdeg2:2}
\left \{%
\begin{array}{c}
{\mathbf N}_{HGME_2/EGM}(\varepsilon_i)=1, i=1,\ldots,4,\\%
{\mathbf N}_{GMHE_2/HGM}(\varepsilon_i)=1, i=1,\ldots,4,%
\end{array} \right .%
\end{equation}%

\begin{equation}%
\label{eqdeg2:3}%
\gamma_i = \varepsilon_i^m, \gamma_i \not =1, i=1, \ldots,4,
\end{equation}%
\begin{equation}%
\label{eqdeg2:4}%
\frac{\gamma_{2j}-1}{\gamma_{2j-1}-1} =a_j-\delta b_j, j=1,2
\end{equation}%
\begin{equation}
\label{eqdeg2:4.2}%
\gamma_3=c+\delta d,
\end{equation}%

\begin{equation}%
\label{eqdeg2:7}%
1 \leq |\sigma(x_1)| \leq Q(B + \sigma(a_1^2-db_1^2)^2),
\end{equation}%
where $\sigma$ ranges over all real embedding of $E_2GM$ into $\tilde{\Q}$,

\begin{equation}%
\label{eqdeg2:9}%
x_1-(a_2-\delta b_2) =(c-1+\delta d)(u+v\delta),
\end{equation}%

\begin{equation}%
\label{eqdeg2:10}%
Px_1Q(B + (a_1^2-db_1^2)^2) \big{|} (c-1+\delta d),
\end{equation}%
where $P$ is a rational prime without any factors in $\calU_{GMHE_2}$.   (For example, $P$ can be any prime
splitting completely in the extension $\Q(\mu)/\Q)$.) Then  $x_1 \in M$.

Conversely, if $x_0 \in \N$, the conditions and equations above can be satisfied in variables ranging over the
prescribed sets.
\end{lemma}%

\begin{proof}%
From (\ref{eqdeg2:2}),  (\ref{eqdeg2:3}) and Corollary \ref{cor:m} we conclude that for all $i=1,\ldots,4$, we have
that $\gamma_i \in O_{HE_2M}$.  Therefore, $c, d \in O_{E_2GM,\calU_{E_2GM}}$. Since $\delta$ generates $HE_2M$
over $E_2M$ and $E_2GMH$ over $E_2GM$, we conclude that $c, d \in O_{E_2M,\calU_{E_2M}}$. A similar argument tells
us that $a_1,b_1, a_2, b_2 \in O_{E_2M,\calU_{E_2M}}$.

Next from (\ref{eqdeg2:1}) and Lemma \ref{le:allprimes} we conclude that for all $\mte{p} \in
\overline{\calU}_{GME_2}$ we have that $\ord_{\pp}x_1\leq 0$ and
\[%
\ord_{\pp}Q(B+ (a_1^2-db_1^2)^2)\leq 0.
\]%
From definition of $B$ (see Notation \ref{not:B}) and the fact that $a_1, b_1 \in E_2M$ - a totally real field, we
have that
\begin{equation}%
\label{eqdeg2:8}%
1 \leq  Q(B + \tau(a_1^2-db_1^2)^2),
\end{equation}%
where $\tau$ ranges over all non-real embeddings of $GME_2$ into $\tilde{\Q}$. Combining the bound equations
(\ref{eqdeg2:7}) and (\ref{eqdeg2:8}), and writing $x_1 = y_0 + y_1\alpha$, where $y_0, y_1\in O_{M,\calU_{M}}$,
we conclude by Lemma \ref{le:boundtwo}
\begin{equation}%
\label{eqdeg2:11}%
|{\mathbf N}_{E_2GM/\Q}(2\alpha y_1)| \leq      |{\mathbf N}_{E_2GM/\Q}(x_1){\mathbf N}_{E_2GM/\Q}(Q(B+
(a_1^2-db_1^2)^2))|.
\end{equation}%

Next consider the divisor $\mte{D}$ of $c-1+\delta d$.  We can write as $\mte{D}=\mte{D}_1\mte{D}_2$, where%
\[%
\mte{D}_1=\prod_{\qq \not \in \calU_{GME_2H}}\mte{q}^{\ord_{\qq}(c-1+\delta d)}
\]%
is an integral divisor, and $\mte{D}_2$ is comprised of primes of $\calU_{GME_2H}$ only. Observe that from
(\ref{eqdeg2:9}), we have that
\[%
x_1-(a_2-b_2\delta)\cong 0 \mod \mte{D}_1 \mbox{ in } O_{GMHE_2,\calU_{GME_2H}}
\]%
Let
\[%
D_1= |{\mathbf N}_{E_2GMH/\Q}(\mte{D}_1)| \in \Z_{>0},
\]%
let
\[%
|{\mathbf N}_{E_2GMH/\Q}(x_1)|=\frac{X}{Y},
\]%
and let
\[%
|{\mathbf N}_{E_2GMH/\Q}(B+Q(a_1^2-db_1^2))|=\frac{U}{V},
\]%
where $X,Y,U,V \in \Z_{>0}, (X,Y)=1, (U,V)=1$, and $X, U$ are not divisible by any rational primes with factors in
$\calU_{GME_2H}$. Then from (\ref{eqdeg2:10}) we have that
\begin{equation}%
\label{eq:XU}
XU < D_1.
\end{equation}%
By Lemma \ref{le:order}, on the one hand we have that
\[%
\frac{Y{\mathbf N}_{E_2GMH/\Q}(2\alpha y_1)}{D_1} \in \Z,
\]%
and therefore
\[%
|Y{\mathbf N}_{E_2GMH/\Q}(2\alpha y_1)| \geq D_1\mbox{ or } y_1=0.
\]%
On the other hand, combining (\ref{eqdeg2:11}) and (\ref{eq:XU}), we have that $\displaystyle |Y{\mathbf
N}_{E_2GMH/\Q}(2\alpha y_1)| \leq XU < D_1$. Thus $y_1$ is 0 and $x_1 \in M$.

We will now show that assuming that $x_0 > 0$ is a natural number, we can satisfy all the equations and conditions
(\ref{eqdeg2:1})--(\ref{eqdeg2:10}). Observe that by (\ref{eqdeg2:1}), we have that $x_1$ is also a natural
number. Let $\nu \in U_{E_2HM} \cap O_{M}[\delta, \mu]$ be a solution to (\ref{sys:norm}) such that it is not a
root of unity. Such a solution exists by Lemma \ref{le:rank1}, Corollary \ref{cor:m} and by Section 2.1.1 of
\cite{Sh16}. Let $\{\phi_1,\ldots,\phi_{s_{E_2HM}}\}$ be a set containing a representative from every
complex-conjugate pair of non-real conjugates of $\nu$. By Lemma \ref{le:closeto1}, we can find a positive integer
$r \cong 0 \mod m$ such that  for all $i=1,\ldots, s_{E_2HM}$ we have that
\[%
\left |\frac{\phi_i^{rA}-1}{\phi_i^r-1}-A\right| < \frac{1}{2},
\]%
where $A=A(x_1) +1$ (see Lemma \ref{le:frombelow}), and thus,
\[%
\left |\frac{\phi_i^{rA}-1}{\phi_i^{r}-1}\right| > A-\frac{1}{2}>A(x_1).
\]%

 So we set $\varepsilon_1=\nu^{r/m}, \gamma_1=\varepsilon^r, \varepsilon_2=\varepsilon^{rA/m},
\gamma_2=\varepsilon^{rA}$. Then for $i=1,2$ the system (\ref{eqdeg2:2}) is satisfied. We also satisfy
(\ref{eqdeg2:3}) for these values of $i$. Next we define $a_1$ and $b_1$ so that (\ref{eqdeg2:4}) is satisfied for
$j=1$. Next let $\sigma$ be an embedding of $M$ into $\tilde{\Q}$ extending to a real embedding of $GM$ and
therefore to a real embedding of $GME_2$. Then by assumption on $H$, we have that $\sigma$ extends to a non-real
embedding $\hat{\sigma}$ on $E_2MH$.  Thus, without loss of generality, for some $i =1,\ldots,s_{E_2HM}$ we have
that
\[%
\hat\sigma(a_1-\delta b_1) = \hat \sigma\left (\frac{\varepsilon^{rA}-1}{\varepsilon^r-1}\right )=
\frac{\phi_i^{rA}-1}{\phi_i^r-1},
\]%
and therefore
\[%
\sigma(a_1^2-db_1^2)=\left | \frac{\phi_i^{rA}-1}{\phi_i^r-1}\right |^2 > A(x_1)^2 >A(x_1),
\]%
leading to
\[%
Q(B+ \sigma(a_1^2-db_1^2)^2) > x_1 = \sigma(x_1) >1.
\]%
Thus we can satisfy (\ref{eqdeg2:7}).

Let $\varepsilon_3$ to be a solution to (\ref{eqdeg2:2}) in $O_{GM}[\delta, \mu]$ such that $\gamma_3=
\varepsilon_3^m \in U_{E_2MH} \cap O_{M}[\delta, \mu]$, (\ref{eqdeg2:3}), (\ref{eqdeg2:4.2}) for $i=3$, and
(\ref{eqdeg2:10}) are satisfied. Again this can be done by Lemma \ref{le:rank1}, Corollary \ref{cor:m} and by
Section 2.1.1 of \cite{Sh16}. Finally, set $\varepsilon_4=\varepsilon_3^{x_1}, \gamma_4=\gamma_3^{x_1}$. In this
case we can satisfy (\ref{eqdeg2:2}), (\ref{eqdeg2:4.2}) for $i=4$.

We now  observe that
\[%
a_2-\delta b_2=\frac{\gamma_4-1}{\gamma_3-1} =x_1 + (\gamma_3-1)(u+\delta v)=x_1 +(c-1-\delta
d)(u+v\delta),
\]%
where $u, v \in O_{GM,\calU_{GM}}[\mu]$.  Thus (\ref{eqdeg2:9}) will also be satisfied.
\end{proof}%

Next we prove a slightly different version of the result above.  We will explicitly use the degree of $M$ over $\Q$.
\begin{lemma}%
\label{le:ext2nf} Let $x, x_0, \ldots,x_{2p_2n} \in O_{GM,\calZ_{GM}}, a_0, b_0,\ldots, a_{2p_2n}, b_{2p_2n}, v, u
\in O_{GM,\calZ_{GM}}[\mu]\subset O_{GME_2H,\calZ_{GME_2H}}$, $\varepsilon_i, \gamma_i \in
O_{GM,\calZ_{GM}}[\mu,\delta] \subset O_{GME_2H, \calZ_{GME_2H}}, i=0,\ldots,2p_2n$.  Assume also that  the
following equations hold.
\begin{equation}
\label{eq:supplement1}
\forall \pp \in \overline{\calE}_{GM}, \,\, \ord_{\pp}x_0 \geq 0
\end{equation}
\begin{equation}%
\label{eqdeg2:1f}
x_k = Q(x + 8l(k+1)), k=0,\ldots, 2p_2n
\end{equation}
\begin{equation}%
\label{eqdeg2:2f}
\left \{%
\begin{array}{c}
{\mathbf N}_{HGME_2/E_2GM}(\varepsilon_i)=1, i=0,\ldots,2p_2n,\\%
{\mathbf N}_{GMHE_2/HGM}(\varepsilon_i)=1, i=0,\ldots,2p_2n,%
\end{array} \right .%
\end{equation}%
\begin{equation}%
\label{eqdeg2:3f}%
\gamma_i = \varepsilon_i^m, \gamma_i \not = 1, i=0, \ldots,2p_2n,
\end{equation}%
\begin{equation}%
\label{eqdeg2:4f}%
\frac{\gamma_{j+1}-1}{\gamma_{0}-1} =a_j-\delta b_j,   j=0,\ldots,2p_2n
\end{equation}%
\begin{equation}
\label{eqdeg2:4.2f}%
\gamma_0=c+\delta d,
\end{equation}%

\begin{equation}%
\label{eqdeg2:9f}%
x_j-(a_j-\delta b_j) =(c+\delta d)(u+v\delta),    j=0,\ldots,2p_2n
\end{equation}%

\begin{equation}%
\label{eqdeg2:10f}%
\left (\prod_{j=0}^{2p_2n} x^2_j\right) \Big{|} (c_0-1+\delta d_0) \mbox{ in } O_{E_2GMH,\calZ_{E_2GMH}}.
\end{equation}%

Then for some $j \in \{0,\ldots,2p_2n\}$ we have $x_j \in M$. Conversely, if $x_0 \in \Z_{>0}$, then equations
(\ref{eqdeg2:1f}) -- (\ref{eqdeg2:10f}) can be satisfied with all the variables in the prescribed sets.
\end{lemma}%
\begin{proof}%
We start as in Lemma \ref{le:inf2} with concluding that $\gamma_j \in E_2HM$ for all $j=0,\ldots, 2p_2n$, and
therefore $a_j, b_j \in O_{E_2M,\calW_{E_2M}}$.  Also as in Lemma \ref{le:inf2}, we note that $\ord_{\pp}x_k \leq
0$ for all $\pp \in \calZ_{GME_2H}$.  By Lemma \ref{le:n>1}, for some $j$ we have that that all the
$\Q$-conjugates of $x_j$ have absolute value greater than 2.  Further, if $x_j = y_{0,j}+ y_{1,j}\alpha$,where
$y_{0,j}, y_{1,j} \in M$,  $\bar{x}_j$ is the conjugate of $x_j$ over $M$, and $\rho : GM \rightarrow \tilde{\Q}$
is an embedding of $GM$ into its algebraic closure,  then
\[%
|2 \rho(\alpha y_{1,j})|=|\rho(x_j)-\rho(\bar{x}_j) | \leq 2\max\{|\rho(x_j)|, |\rho(\bar{x}_j)|\} <
|\rho(x_j)\rho(\bar{x}_j)|.
\]%
Thus,
\begin{equation}%
\label{eqdeg2:11f}%
|{\mathbf N}_{GMHE_2/\Q}(2\alpha y_{1,j})| < |{\mathbf N}_{GMHE_2/\Q}(x_j\bar{x}_j)|= |{\mathbf
N}_{GMHE_2/\Q}(x_j){\mathbf N}_{GMHE_2/\Q}(\bar{x}_j)|= {\mathbf N}_{GMHE_2/\Q}(x_j^2).
\end{equation}%
Next consider the divisor $\mte{D}$ of $c-1+\delta d$.  We can write as $\mte{D}_1\mte{D}_2$, where%
\[%
\mte{D}_1=\prod_{\qq \not \in \calZ_{GME_2H}}\qq^{\ord_{\qq}(c-1+\delta d)},
\]%
is an integral divisor and  $\mte{D}_2$ is divisible by primes of $\calZ_{GME_2H}$ only. Observe that from
(\ref{eqdeg2:9f}), we have that
\[%
x_j-(a_j-b_j\delta)\cong 0 \mod \mte{D}_1
\]%
in $O_{E_2GMH, \calW_{E_2GMH}}$.  Let
\[%
D_1= |{\mathbf N}_{E_2GMH/\Q}(\mte{D}_1)|,
\]%
and let
\[%
|{\mathbf N}_{E_2GMH/\Q}(x_j)|=\frac{X_j}{Y_j}.
\]%
Then by Lemma \ref{le:order2} we conclude that
\[%
\frac{Y^2_j{\mathbf N}_{GME_2H/\Q}(2\alpha y_{1,j})}{D} \in \Z \Rightarrow |Y^2_j{\mathbf N}_{GME_2H/\Q}(2\alpha
y_{1,j})| \geq D_1,
\]%
unless $y_{1,j}=0$.  At the same time,  we also have from (\ref{eqdeg2:10f}) that $X_j^{2} \leq D_1$, and further
from  (\ref{eqdeg2:11f}) we deduce that
\[%
|Y^2_j{\mathbf N}_{GME_2H/\Q}(2\alpha y_{1,j})| <  Y^2_j {\mathbf N}_{GMHE_2/\Q}(x_j^2) = X_j^2 \leq D_1.
\]%
Hence we must conclude that $y_{1,j}=0$.

The argument that the equations above can be satisfied if $x$ is a positive integer is analogous to the argument
used in Lemma \ref{le:inf2}.
\end{proof}%

\section{\bf Diophantine Definability and Decidability in Big  Subrings of Extensions of Degree 2 of Totally Real
Number Fields.}
In this section we will use the technical results from Sections \ref{sec:norms}, \ref{sec:bound2}, and
\ref{sec:deg2} to show that in \emph{any} extension of a degree 2 of a totally real number field, the elements of
$\Q$ contained in some big rings have a Diophantine definition over these rings. Given this definability result,
by now well-explored technique will immediately produce a Diophantine definition of $\Z$ in smaller (but still
big) subrings, as well as a counter examples for the archimedean and non-archimedian versions of a Mazur's
conjecture over these rings.

We start with observing that we have done most of the work in proving the following definability result.

\begin{proposition}%
\label{prop:nf}
$O_{GM,\calZ_{GM}} \cap M$ has a Diophantine definition over $O_{GM,\calZ_{GM}}$.
\end{proposition}%
\begin{proof}%
Lemma \ref{le:ext2nf} will serve as the basis of our proof. First, we define recursively several constants. Let
$N_1$ be a positive integer such that for a any $k, k' \in \{0,\ldots,2p_2n\}$, we have that polynomials
$Q(X+8l(k+1))$ and $Q(X+N_1+8l(k'+1))$ are linearly independent over $\C$. Such a $N_1$ exists by Lemma
\ref{le:linind}. Assume, $N_1,\ldots,N_s, s < p_2$ have been defined recursively, and define $N_{s+1}$ to be a
natural number such that for any $k_0,\ldots,k_s, k_{s+1} \in \{0,\ldots,2p_2n\}$ we have that the set of polynomials
\[%
\{Q(X+8l(k_0+1)),Q(X+N_1+8l(k_1+1)),\ldots,Q(X+N_s+8l(k_s+1)),Q(X+N_{s+1}+8l(k_{s+1}+1))\}
\]%
is linearly independent over $\C$. As above, $N_{s+1}$ exists by Lemma \ref{le:linind}.

Let $N_0=0$ and suppose now that Equations (\ref{eqdeg2:1f})--(\ref{eqdeg2:10f}) are satisfied for
$x=y+N_0,y+N_1,\ldots,y+N_{p_2}$, where $y \in O_{GM,\calZ_{GM}}$. Then by Lemma \ref{le:ext2nf}, for some $k_0,
\ldots,k_{p_2} \in \{0,\ldots,2p_2n\} $ we have that%
\be%
\item $Q(y + N_s+8l(k_s+1)) \in O_{M,\calV_M}$ for $s=0,\ldots, p_2$,\\%
and%
\item the set of polynomials $\{ Q(X + N_s+8l(k_s+1)), s=0,\ldots, p_2 \}$ is linearly independent of $\C$.%
\ee%
Therefore, by Lemma 5.1 of \cite{Sh1} we have that $y \in O_{M,\calZ_M}$. We also know by Lemma \ref{eqdeg2:1f}
that if $y$ is a positive integer then all the equations can be satisfied with variables taking values in the
prescribed sets. To get all the other elements of $O_{M,\calZ_M}$ we can use any integral basis of $M$ over $\Q$.
Thus the only remaining task is making sure that all the Equations (\ref{eqdeg2:1f})--(\ref{eqdeg2:10f}) can be
rewritten in polynomial form with variables ranging over $O_{GM,\calZ_{GM}}$. We can rewrite all the equations
with coefficients and variable in $O_{GM, \calZ_{GM}}$ instead of $O_{GME,\calZ_{GME_2}}$ and
$O_{GMEH,\calZ_{GME_2H}}$ by Proposition \ref{prop:rewriting3} and Proposition \ref{prop:rewriting2}.
\end{proof}%

We can summarize the discussion of the degree 2 extensions of totally real number fields in the following theorem.
\begin{theorem}%
\label{thm:main1}%
Let $K$ be a totally real number field. Let $G$ be any extension of $K$ of degree 2. Let $T$ be any totally real
cyclic extension of $\Q$ of odd prime degree $p >0$ such that $p$ does not divide the degree of the Galois closure
of $K$ over $\Q$. Let $E=KT$. Let $\calX_{G}$ be a set of primes of $G$ such that all but finitely primes of
$\calX_{G}$ are not splitting in the extension $GE/G$. Then $O_{G,\calX_{G}}\cap K$ has a Diophantine definition
over $O_{G,\calX_{G}}$.
\end{theorem}%
\begin{proof}%
Let $K_G$ be the Galois closure of $K$ over $\Q$.  Given our assumption on $p$, we have that
\[%
[EK_G:K_G]=[EG:G]=[E:K]=[T:\Q]=p
\]%
and a rational prime $\mte{P}$ does not split in the extension $T/\Q$ if and only if all of its factors in $K$ and
$G$ do not split in the extensions $E/K$ and $EG/G$ respectively by Lemmas \ref{le:notsplit} and \ref{le:split}.
Further any generator $\mu$ of $T$ over $\Q$ will also generate $E$ over $K$. Thus if $P(X)$ is the monic
irreducible polynomial of $\mu$ over $K$ or over $GK$, it will have rational integer coefficients. Since $p \geq
3$, by Dirichlet Unit Theorem we have that $T$ has units which are not roots of unity. We can set $\mu$ to be
such a unit and satisfy Notation and Assumptions \ref{not:bounds}.

Given that we can define integrality at finitely many primes over number fields (see Proposition
\ref{prop:finmany}), we can restrict all the variables to the values in $O_{G,\calX_G}$ integral at all the primes
splitting in the extension $EG/G$ or dividing the discriminant of $P(X)$ (this set of ``inconvenient'' primes was
denoted by $\calE_G$). Note that we can reconstruct all the values in the ring $O_{G,\calX_{G}}$ by taking the
ratios of the variables whose values are restricted. This is so because we have an existential definition of all
the non-zero values from Proposition \ref{prop:non-zero}. Then by Proposition \ref{prop:nf} we conclude that
$O_{G,\calX_{G}} \cap K$ has a Diophantine definition over $O_{G,\calX_{G}}$.
\end{proof}%

We should note next that the theorem above is a (stronger) analog of Theorem 3.6 and Corollary 3.7 of \cite{Sh3}
where a similar result was proved for totally complex extensions of degree 2 of totally real fields.  Now using
almost exactly the same method as in \cite{Sh3} we can derive analogs of Theorems 3.8, 3.10, 3.11, 3.12, and 3.14
of \cite{Sh3}.  Further using the natural version of the Tchebotarev density theorem (see \cite{Serre3}), we can
replace Dirichlet density by natural density in the statements of all the propositions.  We list the statements of
these theorems below.

\begin{theorem}%
\label{thm:main2}%
Let $K, G, E,  \calX_G$ be as in Theorem \ref{thm:main1}.  Then there exists a set of $G$-primes
$\calN_{G}$ such that $\calX_{G} \subseteq \calN_{G}$, $ \calN_{G} \setminus \calX_{G}$ is
a finite set,  and $O_{G, \calN_{G}} \cap \Q$ has a Diophantine definition over $O_{G,\calN_{G}}$.
\end{theorem}%

\begin{theorem}%
\label{thm:main3}%
Let $G$ be any extension of degree 2 of a totally real field.   Let $\calY_{G}$ be any set of primes of $G$. Then for
any $\varepsilon >0$ there exists a set $\calD_{G}$ such that $\calY_{G} \setminus \calD_{G}$ is
contained in a set of natural density less than $\varepsilon$, $\calD_{G} \setminus \calY_{G}$ is finite, and
 $O_{G, \calD_{G}} \cap \Q$ has a Diophantine definition over $O_{G, \calD_{G}}$.
\end{theorem}%

\begin{theorem}%
\label{thm:main4}%
Let $\calY_{\Q}$ be any set of rational primes. Then for any $\varepsilon > 0$ and any  degree 2 extension $G$ of a
totally real number field, there exists a set of rational primes $\calD_{\Q}$ such that $\calD_{\Q}
\setminus \calY_{\Q}$ is finite, $\calY_{\Q} \setminus \calD_{\Q}$ is contained in a set of primes of
natural density less than $\varepsilon$, and $O_{\Q,\calD_{\Q}}$ has a Diophantine definition in its integral
closure in $G$.
\end{theorem}%

\begin{theorem}%
\label{thm:main5}%
Let $G$ be any  extension of degree 2 of a totally real number field. Let $\chi_G$ be the density of the set of
rational primes splitting completely in $G$. Then for any $\varepsilon >0$ there exists a recursive set
$\calY_{G}$ of primes of $G$ whose natural density is bigger than $1-\chi_{G} - \varepsilon$ and such that $\Z$
has a Diophantine definition over $O_{G,\calY_{G}}$. (Thus, Hilbert's Tenth Problem is undecidable in
$O_{G,\calY_{G}}$.)
\end{theorem}%

\begin{corollary}%
\label{cor:main5}%
Let $G$ be any  extension of degree 2 of a totally real number field. Then for any $\varepsilon >0$ there exists a
recursive set $\calY_{G}$ of primes of $G$ whose natural density is bigger than $1-1/[G:\Q] - \varepsilon$ and
such that $\Z$ has a Diophantine definition over $O_{G,\calY_{G}}$.\\
\end{corollary}%

\begin{theorem}%
\label{thm:main6}%
Let $G$ be any extension of degree 2 of a  totally real number field and let $\varepsilon >0$ be given.
Let $\calY_{\Q}$ be the set of all rational primes splitting in $G$. (If the extension is Galois but not cyclic,
$\calY_{\Q}$ contains all the primes.) Then there exists a set of $G$-primes $\calD_{G}$ such that the set of
rational primes $\calD_{\Q}$ below $\calD_{G}$ differs from $\calY_{\Q}$ by a set contained in a set of
natural density less than $\varepsilon$ and such that $\Z$ has a Diophantine definition over $O_{G,\calD_{G}}$.
\end{theorem}%

As we discussed in the introduction, given Theorem \ref{thm:main2}, we can also reproduce results concerning
existential definability of discrete sets in the archimedean and non-archimedean topologies and a ring version of
Mazur's conjecture on topology of rational points. The proof of these results depends on the analogs of Theorem
\ref{thm:main2} only and therefore can be lifted almost verbatim from the proofs of Theorem 3.6 of \cite{Sh21} and
Theorem 1.8 of \cite{PS}. We state these two results below with Dirichlet density again replaced by natural
density.

\begin{theorem}%
\label{thm:main7}%
Let $G$ be an  extension of degree 2 of a totally real number field. Then for any $\varepsilon > 0$, there exists a
recursive set of $G$-primes $\calY_{G}$ such that the natural density of $\calY_{G}$ is greater $1-\varepsilon$
and there exists an  affine algebraic set $V$ defined over $G$ such that its intersection with
$O_{G,\calY_{G}}$ is infinite and discrete in the usual archimedean topology, and therefore
$\overline{V(O_{G,\calY_{G}})}$, the topological closure of the set of points of $V$ which happen to be in
$O_{G,\calY_{G}}$ in $\C$ if $G$ is not real, and in $\R$, if $G$ is real, has infinitely many connected
components.
\end{theorem}%

\begin{theorem}%
\label{thm:mainnorm}%
 Let $G$ be a  degree-$2$ extension of a totally real number field. Let $\pp$ be any prime of $K$ and let $p_{\Q}$
be the rational prime below it. Then for any $\varepsilon >0$ there exists a recursive set of $G$-primes
$\calY_{G} \ni \pp$ of natural density $> 1-\varepsilon$ such that there exists an infinite Diophantine subset of
$O_{G,\calY_{G}}$ that is discrete and closed when viewed as a subset of the completion $G_\pp$. In fact, such a
subset can be found inside $\Z[1/p_{\Q}]$. \end{theorem}

\section{\bf Diophantine Decidability and Definability over Totally Real Infinite Extensions of $\Q$: an Update.}
\setcounter{equation}{0}%
In this section using the updated version of the norm equations, we update some definability and decidability
results for totally real infinite extensions of $\Q$.  The main difference from our earlier results is in the fact
that we will be able to include factors of any finite set of $K$ primes in the allowed denominators for the rings
under consideration, assuming these primes do not split in the extension $K_{\infty}/K$.
\begin{notationassumption}%
\label{not:2}%
 In this section we will use the following notation and assumptions together with Notation and Assumptions
\ref{not:new}, \ref{not:norms},  \ref{not:bounds}, \ref{not:Q}, and \ref{not:B} which are now assumed to hold for
any field $M$ such that $M$ is contained in a field $K_{\infty}$ described below and $K \subset M$.
\begin{itemize}%
\item Let $K_{\infty}$ be a totally real normal algebraic extension of $\Q$ with $K \subset K_{\infty}$. \item
Assume that only finitely many rational primes are ramified in $K_{\infty}$.%
\item There are only finitely many primes $p$ dividing  $[M:K]$ for any number field $M$ such that $K\subset M
\subset
K_{\infty}$.  %
\item Let $A$  be a positive constant.%
\item Assume that the extension $K_{\infty}/K$ satisfies the following conditions. For any
number field $M$ with $K \subset M \subset K_{\infty}$ besides assumptions described above,  we also  have that
\begin{itemize}%
\item There exists a subfield $\bar M \subset M$ such that $K \subset \bar M$ and $[M:\bar M] \leq A$. %
\item There exists a basis $\Omega=\{\omega_1=1, \omega_2, \ldots, \omega_{n_M}\} \subset O_M$ of $M$ over $\bar M$
such that for all embeddings $\sigma$ of $K_{\infty}$ into its algebraic closure, $|\sigma(\omega_j) | <A$.%
\end{itemize}%
\item Let $D \in \Z_{>0}$ satisfy the following conditions.%
\begin{itemize}%
\item For all $\pp \in \calW_K$ we have that $\ord_{\pp}D=0$.%
\item $D$ is greater than any conjugate of the discriminant of ${\tt D}_{M/\bar M}(\Omega)$ of
$\Omega$ over $\Q$ for any $\Omega$, $M$ and $\bar{M}$ as above.%
\end{itemize}%
\item Let $I_{\calS_K/K_{\infty}}(x,t_1,\ldots,t_k)$ and $u_{\calS_K/K_{\infty}}$ be as in Notation \ref{not:I}.
(Such a polynomial and a rational constant exist by Corollary \ref{cor:deg}
given our assumptions on primes dividing the degrees of subextensions of $K_{\infty}$.)%
\item Let $O_{K_{\infty},\calW_{K_\infty}},O_{K_{\infty},\calS_{K_\infty}}$ be the integral closures of
$O_{K,\calW_K}$ and $O_{K,\calS_K}$ in $K_{\infty}$ respectively (or alternatively one can think of
$\calW_{K_{\infty}}$ and $\calS_{K_{\infty}}$ as being the set of prime ideals of the ring of integers of
$K_{\infty}$ containing all the prime ideals $\pp$
such that $\pp \cap K \in \calW_K$ or $\pp \cap K \in \calS_K$ respectively). %

\item Let $B< l_0< l_1 <\ldots < l_{h_{LE_1}p_2} \in \Z_{>0}$ be a set of positive integers such that the set of
polynomials $\{Q(X+l_i), i=0,\ldots, h_{LE_1}p_2\}$ is linearly independent over $\C$. (Such a set of positive
integers exists by Lemma \ref{le:linind} and the constant $B$ is defined in Notation \ref{not:B}.)%
\item Let $\gamma_{E_1}, \gamma_L$ generate $E_1$ and $L$ over $\Q$.%
 \item Let $C$ be a constant defined in Lemma 4.1 of \cite{Sh26}.

\end{itemize}%

\end{notationassumption}%
The following proposition contains the technical core of this section and is a slightly modification of Proposition
6.2 of \cite{Sh26}.
\begin{proposition}%
\label{prop:totreal}%
Suppose  the following set of equations is satisfied for all $i \in \{0,1,\ldots,h_{KLE_1}p_2\}$;  some
$t_1,\ldots,t_k \in K_{\infty}, y, x_i, y_i \in O_{K_{\infty},\calW_{K_{\infty}}}; \bar{\nu}, \nu,
\bar{\lambda}_i, \lambda_i, \bar{\varepsilon}, \varepsilon_i, w_i, z_i, a_i, Z_i, W_i \in
O_{K_{\infty},\calW_{K_{\infty}}}[\gamma_L,\gamma_{E_1}]$.
\begin{equation}%
\label{eq:1}%
\left \{
\begin{array}{c}
{\mathbf N}_{K_{\infty}E_1L/LK_{\infty}}(\bar{\nu}) = 1,\\%
{\mathbf N}_{K_{\infty}E_1L/EK_{\infty}}(\bar{\nu}) = 1,\\
\bar{\nu} \not = \pm 1,\\
\nu = \bar{\nu}^{2h_{LE_1}}
\end{array}%
\right .
\end{equation}%
\begin{equation}%
\label{eq:2}%
\left \{
\begin{array}{c}%
{\mathbf N}_{K_{\infty}E_1L/LK_{\infty}}(\bar{\lambda}_i) = 1,\\%
{\mathbf N}_{K_{\infty}E_1L/EK_{\infty}}(\bar{\lambda}_i) = 1, \\%
\bar{\lambda}_i \not = \pm 1,\\
\lambda_i = \bar{\lambda}_i^{2h_{LE_1}},
\end{array}%
\right .
\end{equation}%
\begin{equation}
\label{eq:3}
\left \{
\begin{array}{c}
{\mathbf N}_{K_{\infty}E_1L/LK_{\infty}}(\bar{\varepsilon}_i) = 1,  \\%
{\mathbf N}_{K_{\infty}E_1L/E_1K_{\infty}}(\bar{\varepsilon}_i) = 1,    \\%
\bar{\varepsilon}_i \not= \pm 1, \\%
\varepsilon_i = \bar{\varepsilon}_i^{2h_{LE_1}}, \\%
\end{array}
\right .
\end{equation}
\begin{equation}
\label{eq:4}%
\lambda_i - 1=(\nu- 1)z_i,
\end{equation}%
\begin{equation}%
\label{eq:5}
\varepsilon_i - 1=(\nu - 1)w_i,
\end{equation}%
\begin{equation}%
\label{eq:6}
x_i - z_i =(\nu -1)Z_i
\end{equation}%
\begin{equation}%
\label{eq:7}%
\nu-1 = x_i^{2}a_i,%
\end{equation}%
\begin{equation}
\label{eq:8}
x_i = Q(u_{\calS_K/K_{\infty}}y+l_i),
\end{equation}%
\begin{equation}%
\label{eq:8.1}
y_i = x_i^{h_{LE_1}},
\end{equation}
\begin{equation}%
\label{eq:8.2}
 I_{\calS_K/K_{\infty}}(y,t_1,\ldots,t_k)=0,
\end{equation}%
\begin{equation}%
\label{eq:9}
|\sigma(y_i)| > 1, \forall \sigma: K_{\infty} \rightarrow \C,
\end{equation}
\begin{equation}%
\label{eq:10}
 y_i - w_i = y_i^{2A}CDW_i.
\end{equation}%

Then $y \in O_{K,\calW_K}$.

Conversely, if $y \in \Z_{>0}$, then these equations can be satisfied for all $i \in \{0,1,\ldots,h_{KLE_1}p_2\}$;
some $t_1,\ldots,t_k \in K, y, x_i, y_i \in O_{K,\calW_K}; \bar{\nu}, \nu, \bar{\lambda}_i, \lambda_i,
\bar{\varepsilon}, \varepsilon_i, w_i, z_i, a_i, Z_i, W_i \in O_{K,\calW_K}[\gamma_L,\gamma_{E_1}]$.
\end{proposition}%

\begin{proof}%
 Suppose all the equations are satisfied with variables as indicated in the statement of the proposition. Let
$\hat{M}\subset K_{\infty}$ be the smallest overfield of $K$ such that for all $i \in \{0,1,\ldots,h_{LE_1}p_2\}$
we have that $t_1,\ldots,t_k \in \hat{M}, y, x_i, y_i \in O_{\hat{M}, \calW_{\hat{M}}},\bar{\nu}, \nu,
\bar{\lambda}_i, \lambda_i, \bar{\varepsilon}, \varepsilon_i, w_i, z_i \in
O_{\hat{M},\calW_{\hat{M}}}[\gamma_{E_1}, \gamma_L]$.   If $K \not = K(y)=M\subseteq \hat{M}$, then let ${\bar M}$
be a proper subfield of $M$ satisfying the conditions in the Notation and Assumptions \ref{not:2}.

Since, by assumption for any subfield $M$ of $K_{\infty}$ such that $M$ contains $K$ we have that $[ME_1:M]=p_1$
and $[ML:M]=2$ we conclude that $\gamma_{E_1}$ and $\gamma_L$ have the same conjugates over $LK_{\infty}$  and
$E_1K_{\infty}$ respectively as over $L\hat{M}$ and $E_1\hat{M}$ respectively, and therefore we can rewrite the
equations (\ref{eq:1})-- (\ref{eq:3}) as
\begin{equation}%
\label{eq:1.1}%
\left \{
\begin{array}{c}
{\mathbf N}_{\hat{M}E_1L/L\hat{M}}(\bar{\nu}) = 1,\\%
{\mathbf N}_{\hat{M}E_1L/E_1\hat{M}}(\bar{\nu}) = 1,\\
\bar{\nu} \not = \pm 1,\\
\nu = \bar{\nu}^{2h_{LE_1}}
\end{array}%
\right .
\end{equation}%
\begin{equation}%
\label{eq:2.1}%
\left \{
\begin{array}{c}%
{\mathbf N}_{\hat{M}E_1L/L\hat{M}}(\bar{\lambda}_i) = 1,\\%
{\mathbf N}_{\hat{M}E_1L/E_1\hat{M}}(\bar{\lambda}_i) = 1, \\%
\bar{\lambda}_i \not = \pm 1,\\
\lambda_i = \bar{\lambda}_i^{2h_{LE_1}},
\end{array}%
\right .
\end{equation}%
\begin{equation}
\label{eq:3.1}
\left \{
\begin{array}{c}
{\mathbf N}_{\hat{M}E_1L/L\hat{M}}(\bar{\varepsilon}_i) = 1,  \\%
{\mathbf N}_{\hat{M}E_1L/E_1\hat{M}}(\bar{\varepsilon}_i) = 1,    \\%
\bar{\varepsilon}_i \not= \pm 1, \\%
\varepsilon_i = \bar{\varepsilon}_i^{2h_{LE_1}}, \\%
\end{array}%
\right .%
\end{equation}%
Now from Lemma \ref{le:norms} we conclude that $\nu, \lambda_i$ and $\varepsilon_i$ for all $i \in
\{0,1,\ldots,h_{LE_1}p_2\}$ are in $O_{E_1L,\calW_{E_1L}}$. Thus, $z_i, w_i \in O_{E_1L,\calW_{E_1L}}$ for all $i
\in \{0,1,\ldots,h_{LE_1}p_2\}$ as well. Note also that since $y_i \in M$, and $\nu, z_i, w_i \in E_1L$ we have
that $Z_i, W_i \in MLE_1$ by equations (\ref{eq:6}) and (\ref{eq:10}) respectively.  In other words,
(\ref{eq:6})-- (\ref{eq:8.1}) hold over $ME_1L$.  Further, while (\ref{eq:8.2}) a priori  might hold over a bigger
field, its implication about lack of factors of primes in $\calW_{M}$ in the numerator of the divisor of  $x_i$
hold over $ME_1L$. Similarly, (\ref{eq:9}) also holds over $M$.

Given the discussion above, by Lemma 5.1 of \cite{Sh26}, using equations (\ref{eq:6})--(\ref{eq:8.2}), we now
conclude that $y_i= x_i^{h_{LE_1}} = u_i\delta_i^{-1}$, where $u_i \in O_{KLE}$, $\delta_i \in O_{MLE}$ and all
the primes in the divisor of $\delta_i$ are in $\calW_{MLE_1}$. Proceeding further, by a slight modification of
Lemma 5.2 of \cite{Sh26} we obtain from equation (\ref{eq:10}) that $y_i \in \bar{M}LE_1$. (In Lemma 5.2 of
\cite{Sh26} we had that $w_i, W_i \in M$ and $w_i \in K$ instead of $MLE_1$ and $LE_1$ respectively, and the
conclusion is that $y_i \in \bar{M}$. However, the argument is the same for our case.) Since $y_i \in M$ and $M$
and $LE_1$ are linearly disjoint over $\bar{M}$, we conclude that $y_i \in M \cap \bar{M}LE_1=\bar{M}$. The final
step is to note that for all $i$ we have that $y_i=(Q(u_{\calS_K}y+l_i))^{h_{LE_1}} \in \Z[y]$, where and
$(Q(u_{\calS_K}y+l_i))^{h_{LE_1}}$ is of degree $h_{LE_1p_2}$. Thus, if $y_0,\ldots, y_{h_{LE_1}p_2} \in \bar{M}$,
then by Lemma 5.1 of \cite{Sh1}, it is the case that $y \in \bar{M}$ and therefore actually $y \in K$.

Now the satisfiability assertion can be shown in exactly same fashion as it was done in Proposition 6.2 of
\cite{Sh26}.  We should note only that given our choice of $l_i$'s, we can satisfy the equations above for any
positive integer $y$.
\end{proof}%
Before we state the main result of this section, we need to revisit some old number field results.
\begin{theorem}%
\label{thm:old}%
There exists a set $\hat{\calW}_K$ of primes of $K$ such that
$\hat{\calW}_K\setminus\calW_K=\{\ttt_1,\ldots,\ttt_r\}$ is a finite set, no $\ttt_i$ lies above a rational prime
ramified in $K_{\infty}$, and $O_{K,\hat{\calW}_K}\cap \Q$ has a Diophantine definition over
$O_{K,\hat{\calW}_K}$.
\end{theorem} %
\begin{proof}%
Let $K_G$ be the Galois closure of $K$ over $\Q$. Note that due to Notation and Assumption \ref{not:new} we have
that $([K_G:K],p_2)=1$ and therefore $\mu$ is of degree $p_2$ over $K_G$. Let $F_1,\ldots, F_k$ be all the cyclic
subextensions of $K_G$. By Lemma \ref{le:notsplitabove}, for each $i$ there are infinitely many $K_G$-primes
${\mathfrak T}_i$ such that each  ${\mathfrak T}_i$ lies above a $F_i$-prime not splitting in the extension
$K_G/T_i$ and each  ${\mathfrak T}_i$ splits completely in $E_2/F_i$. We claim that by Theorem 2.2 of \cite{Sh3} we
have that $O_{K_G,\calW_{K_G}\cup \{{\mathfrak T}_1, \ldots, {\mathfrak T}_r\}} \cap \Q$ has a Diophantine
definition over $O_{K_G,\calW_{K_G}\cup \{{\mathfrak T}_1, \ldots, {\mathfrak T}_r\}}$. If we compare our data to
the data in Theorem 2.2 of \cite{Sh3}, we will see that we seem to be out of compliance on two points. First of
all we need an element $\gamma$ with $\gamma^2 \in K_G$ such that $K_G(\gamma)$ is totally complex and all
${\mathfrak T}_i$ split in the extension $K_G(\gamma)/K_G$. By the Weak Approximation Theorem we can find $b \in
K_G$ such that all the conjugates of $b$ over $\Q$ are negative and $b \equiv 1 \mod \prod_{i=1}^r {\mathfrak
T}_i$. If we choose a complex number $\gamma$ satisfying $\gamma^2 = b$, then $K_G(\gamma)$ will satisfy the
requirements by Proposition 25 of Section 8, Chapter I and Proposition 16, Section 3, Chapter III of \cite{L}. The
other part out of compliance is the potential presence of finitely many primes in $\calW_{K_G}$ dividing the
discriminant of the power basis of $\mu$ over $K_G$. We take care of this problem by using Proposition
\ref{prop:finmany}.

Now let $\{\ttt_1,\ldots,\ttt_r\}$ be the set  of $K$-primes lying below $\{{\mathfrak T}_1, \ldots, {\mathfrak
T}_r\}$. Let $\hat{\calW}_K = \calW_K \cup \{\ttt_1,\ldots,\ttt_r\}$ and let $\hat{W}_{K_G}$ be the set of all
$K_G$-primes lying above $\hat{W}_K$ primes. Observe that $O_{K_G,\hat{\calW}_{K_G}}$ is the integral closure of
$O_{K,\hat{\calW}_{K}}$ in $K_G$ and $\hat{\calW}_{K_G}\setminus (\calW_{K_G}\cup \{{\mathfrak T}_1, \ldots,
{\mathfrak T}_r\})$ is a finite set. (The extra primes are other factors of $\ttt_i$'s in $K_G$.) Using what
we know about $O_{K_G,\calW_{K_G}\cup \{{\mathfrak T}_1, \ldots, {\mathfrak T}_r\}}$ and Proposition
\ref{prop:finmany} again we can assert that $O_{K_G,\hat{\calW}_{K_G}}\cap \Q$ has a Diophantine definition over
$O_{K_G,\hat{\calW}_{K_G}}$. Finally, by Proposition \ref{prop:rewriting}, we finally conclude that
$O_{K,\hat{\calW}_{K}}\cap \Q$ has a Diophantine definition over $O_{K,\hat{\calW}_{K}}$.
\end{proof}%
We are now ready for the main theorem of this section.
\begin{theorem}%
\label{main:real}
\be%
\item There exist a positive integer $n$ and a polynomial $F(t,\bar{x}) \in K[t,\bar{x}]$ satisfying the following
conditions. For any $t \in O_{K_{\infty},\calW_{K_{\infty}}}$, if there exists $\bar{x} \in
(O_{K_{\infty},\calW_{K_{\infty}}})^n$ such that $F(t,\bar{x})=0$, then $t \in O_{K,\calW_K}$. Further, if $t
\in O_{K,\calW_K}$, there exist $\bar{x} \in (O_{K,\calW_K})^n$ such that $F(t,\bar{x})=0$. Thus,
$O_{K,\calW_K}$ is existentially definable over $O_{K_{\infty},\calW_{K_{\infty}}}$.%

\item There exists a set of $K$-primes $\hat{\calW}_K$, a positive integer $n$, and a polynomial $F(t,\bar{x}) \in
K[t,\bar{x}] $ satisfying the following conditions.
\be%
\item   $\hat{\calW}_K \setminus \calW_K$ is a finite set.%
\item For any $t \in O_{K_{\infty},\hat{\calW}_{K_{\infty}}}$, where $O_{K_{\infty},\hat{\calW}_{K_{\infty}}}$ is
the integral closure of $O_{K,\hat{\calW}_K}$ in $K_{\infty}$, if there exists $\bar{x} \in
(O_{K_{\infty},\hat{\calW}_{K_{\infty}}})^n$ such that $F(t,\bar{x})=0$, then $t \in
O_{K_{\infty},\hat{\calW}_{K_{\infty}}} \cap \Q$. Further, if $t \in O_{K_{\infty},\hat{\calW}_{K_{\infty}}} \cap
\Q$,
there exist $\bar{x} \in (O_{K,\hat{\calW}_K})^n$ such that $F(t,\bar{x})=0$.%
\item $O_{K_{\infty},\hat{\calW}_{K_{\infty}}} \cap \Q$ is existentially definable over
$O_{K_{\infty},\calW_{K_{\infty}}}$.%
\ee%

\item There exists a positive integer $n$ and a polynomial $F(t,\bar{x}) \in K[t,\bar{x}]$ satisfying the following
conditions. For any $t \in O_{K_{\infty},\calS_{K_{\infty}}}$, if there exists $\bar{x} \in
(O_{K_{\infty},\calS_{K_{\infty}}})^n$ such that $F(t,\bar{x})=0$, then $t \in O_{K_{\infty},\calS_{K_{\infty}}}
\cap \Q$. Further, if $t \in O_{K_{\infty},\calS_{K_{\infty}}} \cap \Q$, there exist $\bar{x} \in
(O_{K,\calS_K})^n$ such that $F(t,\bar{x})=0$. Thus, $O_{K_{\infty},\calS_{K_{\infty}}} \cap \Q$ is existentially
definable over $O_{K_{\infty},\calS_{K_{\infty}}}$.
\ee%
\end{theorem}%
\begin{proof}%
Most of the work for the proof of the first assertion has already been done in Proposition \ref{prop:totreal}. We
just have to note that by the discussion in Section \ref{sec:prelim}, all the equations and conditions
(\ref{eq:1})--(\ref{eq:10}) can be rewritten as polynomial equations with coefficients in $K$ and with the
variables ranging in $O_{K_{\infty}, \calW_{K_{\infty}}}$.

To show that the second assertion holds we need to show that $O_{K,\hat{\calW}_K}$ is existentially definable over
$O_{K_{\infty},\hat{\calW}_{K_{\infty}}}$, where following the notational scheme used so far,
$O_{K_{\infty},\hat{\calW}_{K_{\infty}}}$ is the integral closure of $O_{K,\hat{\calW}_K}$ in $K_{\infty}$. Now by
Theorem \ref{thm:old}, we can assume that the new primes allowed in the denominators of divisors are not ramified
in $K_{\infty}$. Thus, we can use the fact that we can define integrality at such primes to obtain the requisite
existential definition. More precisely, let $F(t,\bar{x})$ be the polynomial from the first assertion of the
theorem. Let $\calT_K = \{\ttt_1,\ldots,\ttt_k\}$ and let $I_{\calT_K/K_{\infty}}(x,t_1,\ldots,t_k)$ be defined
as in Notation \ref{not:I}.  Given the choice of primes in $\calT_K$, we can take $u_{\calT_K/K_{\infty}}=1$.
Next consider the following system of equations
\begin{equation}%
\label{eq:extra}%
\left \{
\begin{array}{c}%
F(t,\bar{x})=0\\
I_{\calT_K/K_{\infty}}(t,w_1,\ldots,w_k)=0\\
I_{\calT_K/K_{\infty}}(x_1,w_{1,1},\ldots,w_{1,k})=0\\
\ldots
\end{array}%
\right .
\end{equation}%
Suppose this system has solutions in $O_{K_{\infty},\hat{\calW}_{K_{\infty}}}$.  Then $   F(t,\bar{x})=0$ has
solutions in   $O_{K_{\infty},\calW_{K_{\infty}}}$ and $t \in O_{K,\calW_K}$.  Conversely, if $t \in
O_{K,\calW_K}$, then  $ F(t,\bar{x})=0$ has solutions in  $O_{K,\calW_K}$ and we can find solutions to
$I_{\calT_K/K_{\infty}}(t,w_1,\ldots,w_k)=0, I_{\calT_K/K_{\infty}}(x_1,w_{1,1},\ldots,w_{1,k})=0, \ldots$ in
$K$ also.

To show that the third assertion is true recollect that $\Z$ is existentially definable over $O_{K, \calS_K}$ by
results of Denef and Proposition \ref{prop:finmany}, and since the primes of $\calS_K$ are the only primes which
\emph{have} to be contained in $\calW_K$ in order for the arguments to go through (i.e. to have solutions to the
norm equations in the rings under consideration), the first assertion of the theorem implies the third one.

\end{proof}%

\section{\bf Diophantine Definability and Decidability in the Integral Closure of  Big  and Small Subrings of
Extensions of Degree 2 of Totally Real Algebraic Extensions of $\Q$.}%
\setcounter{equation}{0}%
In this section we consider the definability for the extensions of degree 2 when the underlying totally real field
possibly has an infinite degree over $\Q$.
\begin{notationassumption}%
\label{not:1}%
We start  again with adding to notation and assumptions we have used so far. We continue to think of $M$ as
ranging over all subextensions of $K_{\infty}$ containing $K$ with all the preceding assumptions (i.e assumption
in Assumptions and Notation \ref{not:new}, \ref{not:norms}, \ref{not:bounds}, \ref{not:degree2}, \ref{not:2})
holding for any such $M$.
\begin{itemize}%
\item Let $\calM_K$ be the set of primes of $K$ not splitting in the extension $G/K$. \item Let $\calN_K $ be the
set of $K$-primes not splitting in the extension $E_1E_2G/K$. We will assume that $\calN_K \subset \calM_K \cap
\calV_K \subset \calZ_K$. Let $\calA_K = \calN_K \cup \calS_K$. Given our assumptions we also have that  $\calA_K
\subset \calU_K \cap \calW_K$. %
\item Let $\calL_K$ be formed by removing the highest degree prime of $K$ from every complete set of
$\Q$-conjugates in $\calN_K$.%
\item Let $\calR_K=\calL_K \cup \calS_K$.%
 \item Let $G_{\infty} = K_{\infty}(\alpha)=GK_{\infty}, H_{\infty}=K_{\infty}(\delta)=HK_{\infty}$.%
\item Let $O_{G_{\infty},\calR_{G_{\infty}}}$ be the integral closure of $O_{K,\calR_K}$ in $G_{\infty}$.%
\item Assume that either any root of unity in $E_2GHK_{\infty}$ is already in $G_{\infty}$ or the
group of roots of unity of $E_2GHK_{\infty}$ is finite. In the first case set $m=2$, in the second
case let $m$  be a multiple of the size of the group of roots of unity of $EGHK_{\infty}$.%
 \item Let $O_{K_{\infty}, \calU_{K_{\infty}}}$, $O_{G_{\infty}, \calU_{G_{\infty}}}$ be the integral closures of
$O_{K,\calU_K}$ and $O_{G, \calU_{G}}$ respectively.%
\item  Let $O_{G_{\infty},\calR_{G_{\infty}}}$ be the integral closure of $O_{K,\calR_K}$ in $G_{\infty}$.%
\item Let $O_{K_{\infty},\calS_{K{\infty}}}, O_{G_{\infty},\calS_{G_{\infty}}}$ be the integral closures of
$O_{K,\calS_K}$ in $K_{\infty}$ and $G_{\infty}$ respectively.
\end{itemize}%
\end{notationassumption}%

We will separate the following assumption from the rest, since we will not be using it all the time.  We will
specify explicitly where this assumption is used.

\begin{assumption}%
\label{assum:odd}%
\begin{itemize}%
\item For all number fields $M$ as above we have that $[M:K]$ is odd.%
\end{itemize}%
\end{assumption}%

\begin{lemma}%
\label{le:M_M}%
Suppose Assumption \ref{assum:odd} holds. Then the primes in $\calM_M$ do not split in the extension $GE_2M/M$,
i.e. this notation use is consistent with Notation and Assumption \ref{not:new} and \ref{not:bounds}.
\end{lemma}%

\begin{proof}%
Given Assumption \ref{assum:odd}, we have that  $([M:K],[GE_2:K])=1$ and therefore we can apply Lemmas
\ref{le:notsplit} and \ref{le:split} of the Appendix to reach the desired conclusion.
\end{proof}%
The following lemma also follows from the consideration of the degrees of the extension.
\begin{lemma}
\label{le:Z_G}
The primes of $\calZ_G$ do not split in the extension $E_2G/G$.
\end{lemma}

\begin{proposition}%
\label{prop:infty2}%
Under Assumption \ref{assum:odd} we have that $O_{G_{\infty},\calU_{G_{\infty}}}$ contains a Diophantine subset
$\tt B$
satisfying the following conditions:%
\be%
\item If $x \in \tt B$, then $x \in O_{K_{\infty},\calU_{K_{\infty}}}$.%
\item If $x \in O_{K,\calU_K}$, then $x \in \tt B$.%
\ee%
Alternatively, we can say that there exists a polynomial $P(t,X_1,\ldots, X_l)$ with coefficients in $K$, such that
\begin{equation}
\label{eq:def3} \forall t, X_1, \ldots, X_n  \in O_{G_{\infty}, \calU_{G_{\infty}}}: T_{\calU}(t,X_1,\ldots,X_n)=0
\Longrightarrow t \in K_{\infty} \cap O_{G_{\infty}, \calU_{G_{\infty}}},
\end{equation}%
and
\begin{equation}%
\label{eq:def4}
\forall t \in K, \exists X_1, \ldots, X_n  \in O_{K, \calU_K}: T_{\calU}(t,X_1,\ldots,X_n)=0,
\end{equation}%
\end{proposition}%
\begin{proof}%
The proof will use Proposition \ref{le:inf2} as its foundation.   However we have to adjust somewhat the equations
used in that proposition.  First of all we change the initial range of values for the variables.  Let $x_0, x_1
\in O_{G_{\infty},\calU_{G_{\infty}}}, a_1, a_2, b_1, b_2, u, v \in O_{G_{\infty},\calU_{G_{\infty}}}[\mu]$,
$\varepsilon_i \in O_{G_{\infty},\calU_{G_{\infty}}}[\mu, \delta], i=1,\ldots,4, t_1,\ldots, t_k \in G_{\infty}$
and assume the following conditions and equations are satisfied.
\begin{equation}%
\label{eq:0.1}
 I_{\calE_K/G_{\infty}}(x_0,t_1,\ldots,t_k)=0,
\end{equation}%
\begin{equation}%
\label{eqdeg2:1i}
x_1 = Q(u_{\calE_K/G_{\infty}}x_0),
\end{equation}

\begin{equation}%
\label{eqdeg2:2i}
\left \{%
\begin{array}{c}
{\mathbf N}_{E_2HG_{\infty}/E_2G_{\infty}}(\varepsilon_i)=1, i=1,\ldots,4,\\%
{\mathbf N}_{E_2HG_{\infty}/HG_{\infty}}(\varepsilon_i)=1, i=1,\ldots,4,%
\end{array} \right .%
\end{equation}%

\begin{equation}%
\label{eqdeg2:3i}%
\gamma_i = \varepsilon_i^m, i=1, \ldots,4,
\end{equation}%
\begin{equation}%
\label{eqdeg2:4i}%
\frac{\gamma_{2j}-1}{\gamma_{2j-1}-1} =a_j-\delta b_j, j=1,2
\end{equation}%
\begin{equation}
\label{eqdeg2:4.2i}%
\gamma_i=c_i+\delta d_i, i=1,\ldots, 4,
\end{equation}%
\begin{equation}%
\label{eqdeg2:7i}%
1 \leq |\sigma(x_1)| \leq Q(B+\sigma(a_1^2-db_1^2)^2),
\end{equation}%
where $\sigma$ ranges over all real embedding of $E_2G_{\infty}$ into $\tilde{\Q}$,
\begin{equation}%
\label{eqdeg2:9i}%
x_1-(a_2-\delta b_2) =(c_3+\delta d_3)(u+v\delta),
\end{equation}%
\begin{equation}%
\label{eqdeg2:10i}%
Pax_1Q(B+(a_1^2-db_1^2)^2) \big{|} (c_3-1+\delta d_3),
\end{equation}%
where $P$ is defined as in Proposition \ref{le:inf2} and $I_{\calE_K/G_{\infty}}, u_{\calE_K/G_{\infty}}$ as in
Notation \ref{not:I}.  (We remind the reader that the polynomial $I_{\calE_K/G_{\infty}}$ exists by  Corollary
\ref{cor:deg} and we can choose $u_{\calE_K/G_{\infty}} \in \Z_{>0}$.)  Then, we claim,  $x_1 \in K_{\infty}$.

Conversely, we claim that if $x_0 \in \Z_{>0}$ , the conditions and equations above can be satisfied with $a_1,
a_2, b_1, b_2, u, w \in O_{K,\calU_K}[\mu]$, $\varepsilon_i \in O_{K,\calU_K}[\mu, \delta], i=1,\ldots,4$, and
$t_1,\ldots,t_k \in K$.

To prove the first claim, observe the following. Let $M$  such that $GHE_2M$ contains $\alpha, \delta,\mu,
x_0,a_1,a_2,b_1,b_2, \newline u, v, \varepsilon_i, t_1,\ldots,t_k$. Then given our assumptions on the fields under
consideration, in the equations above we can replace $E_2HG_{\infty}$ by $GE_2HM$, $G_{\infty}E_2$ by $GE_2M$, and
finally $K_{\infty}$ by $M$, while the equalities and other conditions will continue to be true, assuming we
modify the prime sets by choosing the primes in the finite extensions which are below $\calU_{G_{\infty}}$. Then
we can use Proposition \ref{le:inf2} to reach the conclusion that $x_1 \in M \subset K_{\infty}$. The converse
claim follows directly from Proposition \ref{le:inf2}.

The only remaining issue is being able to rewrite all the equations and conditions as polynomial equations with
variables taking values in $G_{\infty}$, and also observe that we can require $x_0 + 1,\ldots,x_0 + p_2$ to satisfy
the equations above.  Here we can proceed exactly as in Proposition \ref{prop:nf}.
\end{proof}%

\begin{remark}%
Let $\hat{\calU}_K$ be a subset of primes of $K$ such that $\hat{\calU}_K \setminus \calU_K$ is a finite set
containing no factors of rational primes ramified in $G_{\infty}$.  Then by Corollary 3.8 the statement of the
proposition above will apply to $O_{G_{\infty},\hat{\calU}_{G_{\infty}}}$ -- the integral closure of
$O_{K,\hat{\calU}_K}$.
\end{remark}%
We now specialize the proposition above for ``small'' rings.  Please note that we {\it do not need} Assumption
\ref{assum:odd} below.

\begin{corollary}%
\label{cor:S}%
$O_{G_{\infty},\calS_{G_{\infty}}}$ contains a Diophantine subset $\tt B$ satisfying the following conditions:%
\be%
\item If $x \in \tt B$, then $x \in O_{K_{\infty},\calS_{\infty}}$.%
\item If $x \in O_{K,\calS_K}$, then $x \in \tt B$.%
\ee%
or alternatively, we can say that there exists a polynomial $T_{\calS}(t,X_1,\ldots, X_l)$ with coefficients in
$K$, such that
\begin{equation}%
\label{eq:def3.2} \forall t, X_1, \ldots, X_n  \in O_{G_{\infty}, \calS_{\infty}}: T_{\calS}(t,X_1,\ldots,X_n)=0
\Longrightarrow t \in K_{\infty} \cap O_{G_{\infty}, \calS_{\infty}},
\end{equation}%
and
\begin{equation}%
\label{eq:def4.2}
\forall t \in K, \exists X_1, \ldots, X_n  \in O_{K, \calS_K}: T_{\calS}(t,X_1,\ldots,X_n)=0,
\end{equation}%
\end{corollary} We can now combine the results above with Theorem \ref{main:real} to obtain the results below.
Observe that we now need the ``majority'' of $K$-primes  not splitting in $E_1/K$  or $E_2/K$.  The analogous
requirements should also hold for primes in the extensions of $K$.  This leads us to use $\calA_K$ as the
set of the allowed denominators.
\begin{theorem}%
\label{main:deg2}
\be%
\item Assume Assumption \ref{assum:odd} holds.  Then there exist a positive integer $n$ and a polynomial
$F(t,\bar{x}) \in K[t,\bar{x}]$ satisfying the following conditions. For any $t \in
O_{G_{\infty},\calA_{G_{\infty}}}$, if there exists $\bar{x} \in (O_{G_{\infty},\calA_{G_{\infty}}})^n$ such that
$F(t,\bar{x})=0$, then $t \in O_{K,\calA_K}$. Further, if $t \in O_{K,\calA_K}$, there exist $\bar{x} \in
(O_{K,\calA_K})^n$ such that $F(t,\bar{x})=0$. Thus,
$O_{K,\calA_K}$ is existentially definable over $O_{G_{\infty},\calA_{G_{\infty}}}$.%

\item Assume Assumption \ref{assum:odd} holds. Then there exists a set of $K$-primes $\hat{\calA}_K$, a positive
integer $n$, and a polynomial $F(t,\bar{x}) \in K[t,\bar{x}] $ satisfying the following conditions.
\be%
\item   $\hat{\calA}_K \setminus \calA_K$ is a finite set.%
\item For any $t \in O_{G_{\infty},\hat{\calA}_{G_{\infty}}}$, where $O_{G_{\infty},\hat{\calA}_{G_{\infty}}}$ is
the integral closure of $O_{K,\hat{\calA}_K}$ in $G_{\infty}$, if there exists $\bar{x} \in
(O_{G_{\infty},\hat{\calA}_{G_{\infty}}})^n$ such that $F(t,\bar{x})=0$, then $t \in
O_{K_{\infty},\hat{\calA}_{K_{\infty}}} \cap \Q$. Further, if $t \in O_{G_{\infty},\hat{\calA}_{G_{\infty}}} \cap
\Q$, there exist $\bar{x} \in (O_{K,\hat{\calA}_K})^n$ such that $F(t,\bar{x})=0$.%
\item $O_{G_{\infty},\hat{\calA}_{G_{\infty}}} \cap \Q$ is existentially definable over
$O_{G_{\infty},\calA_{G_{\infty}}}$.%
\ee%
\item Assume Assumption \ref{assum:odd} holds. Then there exists a set of $K$-primes $\hat{\calR}_K$, a positive
integer $n$, and a polynomial $F(t,\bar{x}) \in K[t,\bar{x}] $ satisfying the following conditions.
\be%
\item   $\hat{\calR}_K \setminus \calR_K$ is a finite set.%
\item For any $t \in O_{G_{\infty},\hat{\calR}_{G_{\infty}}}$, where $O_{G_{\infty},\hat{\calR}_{G_{\infty}}}$ is
the integral closure of $O_{K,\hat{\calR}_K}$ in $G_{\infty}$, if there exists $\bar{x} \in
(O_{G_{\infty},\hat{\calR}_{G_{\infty}}})^n$ such that $F(t,\bar{x})=0$, then $t \in
O_{G_{\infty},\hat{\calR}_{G_{\infty}}} \cap \Q$. Further, if $t \in O_{G_{\infty},\hat{\calR}_{G_{\infty}}} \cap
\Q$, there exist $\bar{x} \in (O_{K,\hat{\calR}_K})^n$ such that $F(t,\bar{x})=0$.%
\item $O_{G_{\infty},\hat{\calR}_{G_{\infty}}} \cap \Q$ is existentially definable over
$O_{G_{\infty},\calR_{G_{\infty}}}$.%
\item $\Z$ is existentially definable over $O_{G_{\infty},\hat{\calR}_{G_{\infty}}}$ and therefore HTP is
undecidable over $O_{G_{\infty},\hat{\calR}_{G_{\infty}}}$.
\ee%
\item There exists a positive integer $n$ and a polynomial $F(t,\bar{x}) \in K[t,\bar{x}]$ satisfying the following
conditions. For any $t \in O_{G_{\infty},\calS_{G_{\infty}}}$, if there exists $\bar{x} \in
(O_{G_{\infty},\calS_{G_{\infty}}})^n$ such that $F(t,\bar{x})=0$, then $t \in O_{G_{\infty},\calS_{G_{\infty}}}
\cap \Q$. Further, if $t \in O_{G_{\infty},\calS_{K_{\infty}}} \cap \Q$, there exist $\bar{x} \in
(O_{K,\calS_K})^n$ such that $F(t,\bar{x})=0$. Thus, $O_{G_{\infty},\calS_{G_{\infty}}} \cap \Q$ and $\Z$ are is
existentially definable over $O_{G_{\infty},\calS_{G_{\infty}}}$, and HTP is not decidable over this ring.
\ee%
\end{theorem}%
\begin{proof}%
The only point which requires clarification  is the definability of $\Z$ over
$O_{G_{\infty},\calR_{G_{\infty}}}$. Here we just point out that by construction,
$O_{G_{\infty},\calR_{G_{\infty}}} \cap \Q$ is a ``small'' subring  of $\Q$, and by Proposition \ref{prop:finmany}
we know that  $\Z$ is definable in ``small'' subrings of $\Q$.
\end{proof}%

The main drawback of the results above is that the numerous conditions make it unclear if any ``nice''  (or for
that matter any) class of infinite algebraic extensions of $\Q$ is covered by the theorems.  In the final section
of the paper we will show that infinite cyclotomics with finitely many ramified rational primes, and consequently
all the abelian extensions embedded in them, satisfy the assumptions of our propositions.

\section{\bf Infinite Cyclotomic and Abelian Extensions.}
To begin with we revisit the issue we have investigated in \cite{Sh26} and \cite{Sh17}. This issue concerns the
number of factors a rational prime can have in a in infinite cyclotomic extension. This matter was investigated in
\cite{Sh26} using an elementary argument. Unfortunately, it was done for odd primes only which is not sufficient
for our current purposes.

\begin{notation}%
We add the following notation.%
\begin{itemize}%
\item Let $\{q_1,\ldots,q_k\}$ be a finite set of rational primes. %
\item For $i =1,\ldots,k$ and any positive integer $j$, let $\xi_{i,j}$ be a primitive $q_i^{j}$-th root of unity.%
\item We specialize $K_{\infty}$ to be the largest totally real subfield of $G_{\infty}=\Q(\{\xi_{i,j},
i=1,\ldots,k, j \in \Z_{>0}\})$.%
\item We now let $M \subset K_{\infty}$ range over number fields contained in $K_{\infty}$. We also vary $K$ across
subfields of $K_{\infty}$ while preserving the assumption that $K \subset M$.%
\item For a rational number $p$, let $g_p(M)$ be the number of factors $p$ has in $M$. Let $g_p(K_{\infty}) =
\max_M\{g_p(M)\}$.%
\item For a rational prime $p \not \in \{q_1,\ldots,q_k\}$, let $n_i$ be the order of $p$ modulo $q_i^{a_i}$, where
$a_i=1$ if $q_i$ is odd, and $a_i=2$ if $q_i=2$. In other words, $n_i$ is the smallest positive integer such that
$p^{n_i} \equiv 1 \mod q_i^{a_i}$.  Also, let $r_i = \ord_{q_i}(p^{n_i}-1)$.%

\end{itemize}%
\end{notation}%
\begin{lemma}%
\label{le:ordmodp}%
Let $x \in \Z$ and let $q$ be a rational prime. Assume further that $\ord_q(x - 1)=n$, where $n$ is a positive
integer. If $q=2$, we will assume that $n\geq2$. Let $l \not = q$ be a prime number. Then $\ord_q(x^l-1)=n$ while
$\ord_q(x^q- 1)= n+1$.
\end{lemma}%
\begin{proof}%
Let $x,q, n,l$ be as in the statement of the lemma and consider the factorization of $x^l-1$ over $\Q(\mu_l)$
where $\mu_l$ is a primitive $l$-th root of unity,
\[%
x^l-1=(x-1)(x-\mu_l)(x-\mu_l^2)\ldots(x-\mu^{l-1})
\]%
Let $\qq$ be a factor of $q$ in $\Q(\mu_l)$, and observe that since $l \not= q$, we have that $\qq$ is not ramified
over $q$.  Thus, on the one hand,
\[%
\ord_q(x^l-1) = \ord_{\qq}(x^l-1) = \ord_{\qq}(x-1) + \sum_{j=1}^{l-1}\ord_{\qq}(x-\mu_l^j).
\]%
On the other hand, $\ord_{\qq}(x- \mu_l^j) = \min( \ord_{\qq}(x- 1),  \ord_{\qq}(1- \mu_l^j)) =0$, since the only
factor occurring in the divisor of $1-\mu_q^j$ is the factor of $l$.

Next consider the factorization of $x^q-1$ in $\Q(\mu_q)$, where $\mu_q$ is a primitive $q$-th root of unity.
\[%
x^q-1=(x-1)(x-\mu_q) \ldots (x-\mu_q^{q-1}),
\]%
Let $\qq$ be the ramified factor of $q$ in this extension.  Then $\ord_{\qq}(1-\mu_q^j) =1$ for $j=1,\ldots,q-1$
and, given our assumptions for the case of $q=2$, we have that  $\ord_{\qq}(1-\mu_q^j) < \ord_{\qq}(x-1)$.
\[%
\ord_{\qq}(x^q-1) = n(q-1) +
\sum_{j=1}^{q-1}\ord_q(x-\mu_q^j)=n(q-1)+\sum_{j=1}^{q-1}\min(\ord_{\qq}(x-1),\ord_{\qq}(1-\mu_p^j))
\]%
\[%
=n(q-1) + (q-1) = (n+1)(q-1).
\]%
Hence,
\[%
\ord_q(x^q-1) = \frac{\ord_{\qq}(x^q-1)}{q-1}= n+1.
\]%
\end{proof}%
From this lemma it immediately follows that the following statement is true.
\begin{corollary}%
\label{cor:iff}%
For any $i = 1,\ldots,k$ and any integers $l \geq 0$, $s >0$, we have that $\ord_{q_i}(p^s-1) \geq r_i +l$ if and
only if $s \equiv 0 \mod n_ip^l$.
\end{corollary}%
Using this corollary we can prove another consequence of Lemma \ref{le:ordmodp}.
\begin{corollary}%
\label{cor:degroot}
Let $\{l_1,\ldots,l_k\}$ be a set of non-negative numbers. Let $m=\prod_{i=1}^kq_i^{l_i+r_i}$. Let
\[%
n=\mbox{LCM}(n_1q_1^{l_1},\ldots,n_kq_k^{l_k}).
\]%
Then $\mu_m$, a primitive $m$-th root of unity, is of degree $n$ over $\F_p$ -- a finite field of $p$ elements.
\end{corollary}%
\begin{proof}%
Let $F$ be a finite field of characteristic $p$. Let $p^s=|F|$. Then $\mu_m \in F$ if and only if $p^s-1 \equiv 0
\mod m$. At the same time, by Corollary \ref{cor:iff}, we have that $p^s-1\equiv0 \mod m$ if and  only if $s \equiv
0 \mod n$, and thus the assertion of the corollary is true.
\end{proof}%

\begin{proposition}%
\label{prop:factors}%
Let $p$ be a rational prime. Then $g_p(K_{\infty}) <\infty$.
\end{proposition}%
\begin{proof}%
It is enough to show that the proposition holds for $\Q(\xi_{i,j}, i=1,\ldots,k, j \in \Z_{>0})$.   We first
consider the extension  $\Q(\xi_i^{r_i+l_i}, i=1,\ldots,k)/\Q$, where, as above, $\{l_1,\ldots,l_k\}$ is a set of
non-negative integers.  Let $\pp$ be a factor of $p$ in this extension and let $f$ be its relative degree.  Since
a power basis of a root of unity is always an integral basis over $\Q$, by Proposition 25 of Section 8, Chapter I
of \cite{L},  to determine $f$, it is enough to determine the degree of $\xi_m$ over $\F_p$, where $m$ is as in
Corollary \ref{cor:degroot}.  By Corollary \ref{cor:degroot}, this degree is equal
$LCM(q_1^{l_1}n_1,\ldots,q_k^{l_k}n_k)$.  Since $p$ is not ramified in the extension $\Q(\xi_i^{r_i+l_i},
i=1,\ldots,k)/\Q$, we can conclude that the number of factors of $p$ in $\Q(\xi_i^{r_i+l_i}, i=1,\ldots,k)$ is
equal to $\displaystyle \frac{q_1^{r_1+l_1-1}(q_1-1)\ldots
q_k^{r_k+l_k-1}(q_k-1)}{LCM(q_1^{l_1}n_1,\ldots,q_k^{l_k}n_k)}\leq q_1^{r_1}\ldots q_k^{r_k}$.
\end{proof}%

\begin{lemma}%
Assuming that ramification degree of 2 is finite,  there exists a number field $K \subset K_{\infty}$ such that
for all number fields $M$ with $K \subset M \subset K_{\infty}$ we have that $[M:K]$ is odd.  Further the same
assertion will be true for any finite extension of $K$ in $K_{\infty}$.
\end{lemma}%
\begin{proof}%
Given our assumptions, without loss of generality,  we can assume that
\[%
G_{\infty}=\Q(\xi_{1,r}, \{\xi_{i,j}, i=2,\ldots,k,j \in \Z_{>0}\}),
\]%
where $q_1=2$ and $r \in \Z_{>0}$. Let $G=\Q(\xi_{1,r}, \xi_{2,1}, \ldots,\xi_{k,1})$. Then for any number field
$R$ with $G \subset R \subset G_{\infty}$, we have that $[R:G]=\prod_{i=2}^kq_i^{a_i}, a_i \in \Z_{\geq 0}$.  Set
$K$ to be the largest totally real subfield of $G$ and consider a number field  $M \subset K_{\infty}$ with $K
\subset M$.  For some number field $R$, as above, we have that $[R:M]=2$, and therefore by comparing degrees we can
conclude that $[M:K]=[R:G]$ is an odd number.
\end{proof}%

 We now consider various other assumptions on $K_{\infty}$ and
$G_{\infty}$ and prove the following proposition.
\begin{proposition}%
\label{prop:satisfy}%
Let  $K_0 \subset K_{\infty}$, and let $\calS_{K_0}$ be a finite set of primes of $K_0$.  Then for some infinite
set $\calA_{K_0}$ of $K_0$-primes containing $\calS_{K_0}$ there exists a finite extension $K$ of $K_0$ and finite
extensions $E_1, E_2$, and $L$ of $K$ so that all the assumptions in Notation and Assumptions \ref{not:new},
\ref{not:norms}, \ref{not:bounds}, \ref{not:degree2}, and \ref{not:2} hold with respect to $\calS_K$ and
$\calA_{K}$ -- the sets of $K$-primes above $\calS_{K_0}$ and $\calA_{K_0}$ respectively.  Further if $K_0 \not =
\Q$ we can arrange that all the primes of $\calN_{K_0} = \calA_{K_0}\setminus \calS_{K_0}$ lie above rational
primes splitting completely in the extension $K_0/\Q$ so that if we remove one prime from $\calN_{K_0}$ per every
complete set of $\Q$-conjugates, the remaining set $\calL_{K_0}$  will still be infinite.
\end{proposition}%
\begin{proof}%
If  we set $K$ to be any number field containing $\Q(\cos(2\pi/q_1^{a_1}\ldots q_k^{a_k}))$, where $a_i=1$ if $q_i$
is odd and $a_i=2$ if $q_i=2$, then any number field $M$ with $K \subset M \subset K_{\infty}$ will have a subfield
$\bar{M}$ such that $[M: \bar M]= q_i$ for some $i=1,\ldots,k$ and   $M=\bar M(\cos(2\pi/q_i^{b_i}))$, where $b_i$
is a positive integer and $2\cos(2\pi/q_i^{b_i})=\xi_{i,b_i} + \xi^{-1}_{i,b_i}$ is an algebraic integer.  Thus
the assumptions that only finitely many primes divide the degrees of subextensions and the integral basis elements
and their conjugates  are bounded in absolute value hold.   Further the condition on finite number of rational
primes ramified in $K_{\infty}$ also holds by our choice of $K_{\infty}$.

Our next job is to make sure that primes of $\calS_K$ do not split in the extensions $K_{\infty}/K$.  Since every
prime can have only finitely many factors in $K_{\infty}$ we can certainly choose a number field $K$ contained in
$K_{\infty}$ so that  it  contained $\Q(\cos(2\pi/q_1^{a_1}\ldots q_k^{a_k}))$, where $a_i=1$ if $q_i$ is odd and
$a_i=2$ if $q_i=2$ and the maximum possible number of factors for each prime in $\calS_{K_0}$.  Then the primes of
$\calS_K$ will remain prime in the extension $K_{\infty}/K$.

We now produce cyclic extensions $E_1$ and $E_2$ with the required properties. Choose two distinct odd rational
prime numbers $p_1$ and $p_2$ such that each $p_i$ is prime to $\prod (q_i-1)q_i$.  By Lemma \ref{le:extensions},
there exists a cyclic degree $p_1$ extension $\hat{E}_1$ of $\Q$ such that all the prime below $\pp_1,\ldots,
\pp_s$ split completely in the extension $\hat{E}_1/\Q$. Also by Lemma \ref{le:extensions}, there exists a cyclic
degree $p_2$ extension $\hat{E}_2$ of $\Q$ such that all the prime below $\pp_1,\ldots, \pp_s$ do not split in the
extension $\hat{E}_2/\Q$. Now set $E_1=K\hat{E}_1, E_2=K\hat{E}_2$. Then, given the degrees of the extensions
involved,  by Lemmas \ref{le:notsplitabove} and \ref{le:split} we have that $\pp_1,\ldots,\pp_s$ split completely
in the extension $E_1/K$ and do not split in the extension $E_2/K$.  We also not here that since $\hat{E}_1$ and
$\hat{E}_2$ are  Galois extensions of $\Q$ of odd degree, they must be totally real.

We still have to construct $L$ so that $\pp_i$'s split in $L/K$, choose a generator for $H$ with the correct sign
for the conjugates and make sure that the requirements for roots of unity are satisfied: $LE_1K_{\infty}$ should
have no roots of unity beyond $\pm1$ and $G_{\infty}E_2H$ should not have any roots of unity which are not already
in $G_{\infty}$.

To make sure that $LE_1K_{\infty}$ has no non-real roots of unity, it is enough, by Lemma 2.4 of \cite{Sh17} to
make sure that in the extension $LE_1/E_1$ we have ramification of at least two $K$-primes lying above two
different rational primes. We can do this by choosing $c \in K$ such that besides satisfying the inequality
$\sigma(c) <0$ for all embeddings $\sigma$ of $K$ into its algebraic closure and  equivalencies $c\equiv 1 \mod
\pp_i$, it is the case that $c$ also satisfies the condition that it has order 1 at two $K$-primes lying above two
different rational primes. Such an $c$ can be found by the Weak Approximation Theorem. Now we can set
$L=K(\sqrt{c})$.

Now note that by Lemma \ref{le:extensions} we can make sure that no factor of $2,3, q_1,\ldots, q_k$ is ramified in
the extension $E_2/K$. Let $\ell $ be the prime ramified in this extension and note that the ramification degree
of $\ell$ is $p_2 <\ell -1$. Next let's consider $G_{\infty}E_2$ and note that the only primes which are ramified
in any finite subextension of this field are $q_1,\ldots,q_k, \ell$.

Now choose $d \in O_K$  ($\sqrt{d}$ will generate $H$) so that $\sigma(d) >0$ for all embeddings $\sigma$ of $K$
into its algebraic closure, $d \equiv 1 \mod 4$ and $d$ is a unit at all the factors of $2, 3, q_1,\ldots, q_k$
and $\ell$.   Now the only rational primes ramified in any finite subextension of $G_{\infty}HE_2$ are
$q_1,\ldots,q_k,\ell,$ and factors of $d$.  Thus, the only ``extraneous'' roots of unity which can occur are
$\xi_{\ell}$ and $\xi_t$, where a factor of $t$ divides $d$.  But if $\xi_{\ell} \in G_{\infty}HE_2$, then for some
$GHE_2$ we must have ramification of $\ell$ in the extension equal to $\ell -1$, which is not the case by the
argument above.  Similarly, if $\xi_t \in G_{\infty}HE_2$, then $t$ must have ramification at least $t-1$ in some
$GHE_2$.  However, by construction, the ramification degree of $t$ can be at most 2 and $t >3$.

Now we are ready to choose $\calA_{K_0}$.  We will describe $\calN_{K_0}$ and add $\calS_{K_0}$ to the set.  By
Lemma \ref{le:infset} there exists an infinite set $\calN_{\Q}$ of rational primes $P$ such that
\be%
\item $P$ splits completely in the extension $K/\Q$ and therefore also in the extension $K_0/\Q$.
\item No factor of $P$ in $K$ splits in the extension $E_1E_2G/K$.
\ee
Now let $\calN_{K_0}$ consist of all the $K_0$-primes lying above the primes of $\calN_{\Q}$.

Finally we note that since only finitely many primes divide the degrees of subextensions, i.e divisors of
$\prod_{i=1}^qq_i(q-1)$, and all the subextensions are Galois.  Therefore, by Corollary \ref{cor:deg}  we will be
able to use  Theorems \ref{thm:ifmp} and \ref{thm:ifmpu}
\end{proof}%

  We are now ready for the following theorem.

\begin{theorem}%
\label{thm:cyclo} %
Let $G_{\infty}$ be an infinite cyclotomic extension of $\Q$ with finitely many ramified rational primes and finite
ramification degree for 2. Let $K_0 \not = \Q$ be a totally real number field contained in $G_{\infty}$.  Then for
some large subring $O_{K_0,\calR_{K_0}}$ of $K_0$ its integral closure $O_{G_{\infty}, \calR_{G_{\infty}}}$ in
$G_{\infty}$ satisfies the following conditions:
\be%
\item   $O_{G_{\infty}, \calR_{G_{\infty}}} \cap \Q= O_{\Q,\calS_{\Q}}$, where $\calS_{\Q}$ is finite.%
\item There exists a positive integer $n$ and a polynomial $F(t,\bar{x}) \in K_0[t,\bar{x}]$ satisfying the
following conditions. For any $t \in O_{G_{\infty},\calR_{G_{\infty}}}$, if there exists $\bar{x} \in
(O_{K_{\infty},\calR_{G_{\infty}}})^n$ such that $F(t,\bar{x})=0$, then $t \in O_{G_{\infty},\calR_{G_{\infty}}}
\cap \Q$.  Further, if $t \in O_{G_{\infty},\calR_{G_{\infty}}} \cap \Q$, there exist $\bar{x} \in
(O_{K_0,\calR_{K_0}})^n$ such that $F(t,\bar{x})=0$.  Thus $\Z$ is definable over  $O_{G_{\infty},
\calR_{G_{\infty}}}$.
\ee%
 It is also possible to find a totally real number field $K \subset K_{\infty}$ such that the
Dirichlet  (or natural) density of $\calR_K$ can be made arbitrarily close to 1/2.
\end{theorem}%
\begin{proof}%
This theorem follows almost immediately from Theorem \ref{main:deg2} and Proposition \ref{prop:satisfy}, if we
let $K$, $\calA_{K_0}, \calL_{K_0}$ be constructed as in  Proposition \ref{prop:satisfy} and set $\calR_{K_0} =
\calL_{K_0}\cup \calS_{K_0}$, since bounded ramification for 2 implies that Assumption \ref{assum:odd} holds for
some finite extension $K$ of $K_0$. There is only one point which requires explanation:  the question of density.
We now show how  the density of $\calA_K$ can be arranged to be arbitrarily close to $1/2$.  Here we can
assume that $K=K_0$ satisfies Assumption \ref{assum:odd} and review the definition of $\calR_K$. It can be formed
in several steps. We start with $\calN_K$-- a set of $K$ primes not splitting in the extensions $E_1/K$, $E_2/K$
and $G/K$. Next out of $\calN_K$ we form a set of $K$-primes $\calL_K$ by removing the highest degree prime out of
every complete set of $\Q$-conjugates in $\calL_K$. Finally, we let $\calR_K = \calL_K \cup
\calS_K$. Here we note that if a number field $K$ satisfies the assumptions of
Theorem \ref{main:deg2}, then so does any finite extension of $K$ inside $K_{\infty}$.  Hence when required we can
make the degree of $K$ arbitrarily large.  We start with the fact that, by Tchebotarev Density Theorem (the classic
or natural versions), the density (Dirichlet or natural) of primes of $K$ not splitting in the extension $G/K$ is
$1/2$.  However, out of this set of $K$-primes we have to remove the primes splitting either in $E_1$ (density
$1/p_1$) or $E_2$ (density $1/p_2$) and primes of the highest relative degree in complete sets of $\Q$-conjugates.
Since the only primes which contribute to density are primes of relative degree  1 over $\Q$, we should worry only
about complete conjugates sets lying above rational primes splitting completely in the extension $K/\Q$.  The
density of the set containing exactly one prime for each complete set of conjugates lying above a completely
splitting rational prime is $\frac{1}{[K:\Q]}$.  Using Tchebotarev Density Theorem and a Galois extension
$GE_1E_2/\Q$ one can deduce that the set of removed primes has density (natural and Dirichlet), and by making the
degree of $K$ over $\Q$, and the degrees of $E_1$, and $E_2$ over $K$  high enough we can make this density
arbitrarily small.
 
\end{proof}%
We now extend results above to complex number fields contained in $G_{\infty}$.
\begin{corollary}%
\label{cor:cyclo}
Let $G_0$ be a number field contained in $G_{\infty}$.  Then for
some large subring $O_{G_0,\calR_{G_0}}$ of $G_0$, its integral closure $O_{G_{\infty}, \calR_{G_{\infty}}}$ in
$G_{\infty}$ satisfies the following conditions:
\be%
\item   $O_{G_{\infty}, \calR_{G_{\infty}}} \cap \Q= O_{\Q,\calS_{\Q}}$, where $\calS_{\Q}$ is finite.%
\item There exists a positive integer $n$ and a polynomial $F(t,\bar{x}) \in G_0[t,\bar{x}]$ satisfying the
following conditions. For any $t \in O_{G_{\infty},\calR_{G_{\infty}}}$, if there exists $\bar{x} \in
(O_{K_{\infty},\calR_{G_{\infty}}})^n$ such that $F(t,\bar{x})=0$, then $t \in O_{G_{\infty},\calR_{G_{\infty}}}
\cap \Q$.  Further, if $t \in O_{G_{\infty},\calR_{G_{\infty}}} \cap \Q$, there exist $\bar{x} \in
(O_{G_0,\calR_{G_0}})^n$ such that $F(t,\bar{x})=0$.  Thus $\Z$ is definable over  $O_{G_{\infty},
\calR_{G_{\infty}}}$.
\ee%
\end{corollary}%

\begin{proof}%
Either $G_0$ is totally real and we are done, or $G_0$ is an extension of degree 2 of some totally real number
field $K_0$.  In the latter case construct $O_{K_0,\calR_{K_0}}$ as in Theorem \ref{thm:cyclo} and let
$O_{G_0,\calR_{G_0}}$ be the integral closure   $O_{K_0,\calR_{K_0}}$ in $G_0$.
\end{proof}%
Our next goal is to apply results above to small rings -- rings of $\calS$-integers with finitely many primes in
$\calS$.  To obtain the most general results we need a technical lemma.
\begin{lemma}%
\label{le:dioph}
As above let $G/K$ to be an extension of degree 2 generated by $\alpha \in O_G$, and let $\calC_G$ be a set of
$G$-primes such that no prime of $\calC_G$ has its $K$-conjugate in $\calC_G$ (thus $O_{G,\calC_G} \cap K = O_K$).
 Let $\calS_K$ be the set of all the  $K$-primes of lying below the primes of $\calC_G$.  Let
$O_{K_{\infty},\calS_{K_{\infty}}}, O_{G_{\infty},\calC_{G_{\infty}}}$  be the integral closures of $O_{K,\calS_K}$
and  $O_{G,\calC_G}$ respectively.  Next let ${\mathbf A} \subset O_{K_{\infty}}$ and be such that $O_K \subset$
${\mathbf A}$.  Now consider the following set
\[%
{\mathbf B} =\{\frac{x}{y}: x, y \in {\mathbf A} \land (\exists a, b, c, d, e, f \in
O_{K_{\infty}})(cde \not =0 \land (e,f) =1 \land   \frac{a}{c}+\alpha\frac{b}{d} \in
O_{G_{\infty},\calC_{G_{\infty}}} \land 2\frac{a}{c}=\frac{e}{f} \land\frac{fx}{y} \in {\mathbf A}\}.
\]%
We claim that $O_{K,\calS_K} \subset {\mathbf B} \subset O_{K_{\infty},\calS_{K_{\infty}}}$.
\end{lemma}%
\begin{proof}%
Let $\pp_K \in \calS_K$, let $\pp_G \in \calC_G$ be a prime above $\pp_K$ in $G$,  and let $m \in \Z_{>0}$ be a
multiple of the class numbers of $K$ and $G$ .  Then there exists $a, c, b, d \in O_K, cd \not = 0$ such that
\[%
\frac{a}{c}+\alpha\frac{b}{d} \in O_{G,\calC_G}
\]%
and
\[%
\ord_{\pp_G} \left (\frac{a}{c}+\alpha\frac{b}{d}\right ) = -m,
\]%
while being integral at all the other primes.  Since $\bar{\pp}_G$ - the conjugate of $\pp_G$ over $K$, is not in
$\calC_G$, we conclude that
\[%
\ord_{\bar{\pp}_G} \left (\frac{a}{c}+\alpha\frac{b}{d} \right ) =\ord_{\pp_G} \left (\frac{a}{c}-\alpha\frac{b}{d}
\right) =0.
\]%
Consequently,
\[%
\ord_{\pp_K}\left (\frac{2a}{c}\right )=\ord_{\pp_G} \left (\frac{a}{c}+\alpha\frac{b}{d} \right )= -m.
\]%
Using the finiteness of the $K$-class number again, we can find $e, f \in O_K$ such that $(e,f)=1$ and
$\displaystyle \frac{e}{f} = \frac{2a}{b}$.  Then we conclude that
\[%
\ord_{\pp_K} f = -\ord_{\pp_K}\left (\frac{2a}{c}\right )  =-\ord_{\pp_G}\left
(\frac{a}{c}+\alpha\frac{b}{c}\right) = m.
\]%
Hence, if  $\displaystyle\frac{x}{y} \in O_{K, \calS_K}$, then there exists $f$ as above (essentially the
$K$-``denominator'' of an element of $O_{G,\calC_G}$)  so that $\displaystyle\frac{fx}{y} \in {\mathbf A}$.

Conversely, suppose $\displaystyle\frac{a}{c}+\alpha\frac{b}{d} \in GM \cap O_{G_{\infty},\calC_{G_{\infty}}}$ for
some finite extension $M$ of $K$ in $K_{\infty}$ and $\ord_{\pp_M}f < 0$ for some prime of $M$, where $a,b,c,d,e,f
$ are as described in the statement of the lemma. Then
$\displaystyle\ord_{\pp_M}\frac{e}{f}=\ord_{\pp_M}(2\frac{a}{c}) <0$, and consequently $\pp_M \in \calS_M$. So if
$\displaystyle\frac{fx}{y} \in {\mathbf A} \subset O_{K_{\infty}}$, then $\displaystyle\frac{x}{y} \in
O_{K_{\infty}, \calS_{\infty}}$.
\end{proof}%
\begin{remark}%
A more natural way to state the lemma above would be to assert that for some set ${\mathbf B}$ such that
$O_{K,\calS_K} \subset {\mathbf B} \subset O_{K,\calS_{K_{\infty}}}$ we have that ${\mathbf B} \dg
{\mathbf A} \dg O_{G_{\infty},\calC_{G_{\infty}}}$.  (For a discussion of Diophantine generation see
either \cite{Sh5} or \cite{Sh34}.)
 
\end{remark}%
We are now ready to deal with rings of $\calS$-integers where $\calS$ is finite, also known as ``small'' rings.
\begin{theorem}%
\label{thm:sint}

Let $G_{\infty}$ be a cyclotomic extension of $\Q$ with finitely many ramified rational primes. Let $R$ be any
number field contained in $G_{\infty}$ and let $\calS_R$ be a non-empty finite set of primes of $R$.  Then  there
exists a positive integer $n$ and a polynomial $F(t,\bar{x}) \in K[t,\bar{x}] $ satisfying the following
conditions. For any $t \in O_{G_{\infty},\calS_{G_{\infty}}}$, if there exists $\bar{x} \in
(O_{K_{\infty},\calS_{G_{\infty}}})^n$ such that $F(t,\bar{x})=0$, then $t \in O_{G_{\infty},\calS_{G_{\infty}}}
\cap \Q$. Further, if $t \in O_{G_{\infty},\calS_{G_{\infty}}} \cap \Q$, there exist $\bar{x} \in
(O_{R,\calS_R})^n$ such that $F(t,\bar{x})=0$. Thus, $O_{G_{\infty},\calS_{G_{\infty}}} \cap \Q$ is definable over
$O_{G_{\infty},\calS_{G_{\infty}}}$. Consequently $\Z$ is existentially definable in the integral closure of
$O_{R,\calS_R}$ in $G_{\infty}$ and Hilbert's Tenth Problem is not decidable over
$O_{G_{\infty},\calS_{G_{\infty}}}$.

\end{theorem}%
\begin{proof}%
If the number field where we select the ring of $\calS$-integers is totally real (a field $K$ in our notation),
then the assertion of the theorem follows directly from Theorem \ref{thm:cyclo}. If, however, the field is totally
complex (a field $G$ in our notation), we have to be more carefully. Let $\calC_G$ and
$O_{G_{\infty},\calC_{G_{\infty}}}$ be defined as above. Since the construction of a Diophantine definition of
$O_{K_{\infty},\calS_{K_{\infty}}}$ over $O_{G_{\infty},\calS_{G_{\infty}}}$, where $\calS_K$ is defined as usual
to be a finite set of primes of some totally real number field $K$ lying below primes of $\calC_G$ with $[G:K]=2$,
was carried out over the ring of integers, while we ``neutralized'' the ``denominators'' by using a polynomial
$Q(X)$ and the fact that all the primes allowed in the denominator of the divisors of elements of the rings in
question did not split in the extension $E_2G_{\infty}/K_{\infty}$ (see Section \ref{sec:norms} and Lemma
\ref{le:inf2}), we can replicate this process no matter what finite set of primes of $G$ we select. However, the
difficulty can arise when we find ourselves in $K_{\infty}$. In order to carry out the totally real part of the
construction we need at least one ``prime in the denominator'', i.e. if $O_{G_{\infty}, \calC_{G_{\infty}}} \cap
K_{\infty} = O_{K_{\infty}}$ we will not be able to proceed directly. We need somehow to construct $O_{K,
\calS_K}$ so that we have solutions to norm equations (\ref{norm:1}).   To show that this can be in fact be done,
we use Lemma \ref{le:dioph}.  We note that a set ${\mathbf A}$, as described in the statement of the lemma is
indeed Diophantine over $O_{G_{\infty}, \calC_{G_{\infty}}}$ and Lemma \ref{le:dioph} then tells us that using
polynomial equations we can represent elements of a set ${\mathbf B}$ containing $O_{K,\calS_K}$.  Consequently,
the totally real part of the construction can be carried out.

Finally we note that being non-zero and relatively prime in a ring
of integers are Diophantine conditions by Propositions \ref{prop:non-zero} and \ref{prop:relprime}.
\end{proof}%
We now can make use of Kronecker-Weber Theorem and Lemma \ref{le:KW} to assert the following.
\begin{theorem}%
Let $A_{\infty}$ be an abelian extension of $\Q$ with finitely many ramified rational primes. Then the following
statements are true.
\begin{itemize}%
\item If the ramification degree of 2 is finite, then for any number field $X \subset A_{\infty}$ there exists an
infinite set of $X$-primes $\calW_X$ such that $\Z$ is existentially definable in the integral closure of
$O_{X,\calW_X}$ of $A_{\infty}$.%
\item For any number field $X \subset A_{\infty}$ and any finite non-empty set $\calS_X$ of its primes we have that
$\Z$ is existentially definable in the integral closure of $O_{X,\calS_X}$ in $A_{\infty}$.
\end{itemize}%
\end{theorem}%
\begin{proof}%

Let $A_{\infty}$ be an abelian equation with finitely many ramified primes, and assume that the ramification degree
of 2 is finite (possibly 1). Then by Lemma \ref{le:KW} we have that $A_{\infty} \subset G_{\infty} = \Q(\xi_{1,r},
\{\xi_{i,j}, i=2,\ldots,k,j \in \Z_{>0}\})$, where $p_1=2$ and $r$ is a positive integer. Now the theorem follows
from Corollary \ref{cor:cyclo} and Theorem \ref{thm:sint} via polynomials $F(t,\bar{x})$.
\end{proof}%

\begin{remark}%
While infinite abelian and infinite cyclotomic extensions with finitely many ramified rational primes are probably
the ``nicest'' examples of the fields to which our results apply, they are certainly not the only ones.   One can
produce more examples by starting with a totally real subfield of an infinite cyclotomic with finitely many
ramified rational primes, attaching it to an arbitrary totally real number field and then adding an arbitrary
extension of degree 2.
\end{remark}%

\section{Appendix}
This section contains some technical  results used in the paper. This first  lemma is a modification of Lemma 8.1
of \cite{Sh26}.%
\begin{lemma}
\label{le:linind}%
Let $K$ be a real number field. Let $F(T) \in K[T]$ be a polynomial of degree $n> 0$. Suppose that for some
positive numbers $l_0, l_1,\ldots, l_k, k < n$, we have that polynomials $F(T+l_0), F(T+l_1),\ldots, F(T+l_k)$ are
linearly independent over $\C$. Then there exist a positive constant $C$ such that for any real $l >C$, polynomials
$F(T+l_0), F(T+l_1),\ldots, F(T+l_k), F(T+l)$ are also linearly independent over $\C$.
\end{lemma}%
\begin{proof}%
 Let $F(T) = a_0 + a_1T + \ldots a_nT^n$. Then for $l \in \N$ we have that
 \[%
 F(T+l) = a_0 + a_1(T+l) + \ldots + a_n(T+l)^n =
\]%
\[%
a_0 + a_1(T+l) + \ldots + a_i\left(\sum_{j=0}^i \left ( \begin{array}{c} i\\j\end{array}\right ) T^jl^{i-j}\right ) + \ldots +
a_n\left (\sum_{j=0}^n\left ( \begin{array}{c}n\\j\end{array}\right ) T^jl^{n-j}\right ) =
\]%
\[%
\sum_{j=0}^na_jl^j + \ldots + \left(\sum_{j=k}^n\left ( \begin{array}{c} j\\k\end{array}\right )a_jl^{j-k}\right ) T^k + \ldots
+ a_nT^n=
\]%
\begin{equation}%
\label{sys:1}%
 P_n(l) + P_{n-1}(l)T + \ldots P_0(l)T^n,%
\end{equation}%
where $P_i(l)\in K[l]$ is a polynomial of degree $i$ in $l$.  (The coefficient of $l^i$ in $P_i(l)$ is $a_n \not =
0$ by assumption on the degree of $F(T)$.)  Let $F_k(T+l)= \sum_{j=0}^kP_j(l)T^{n-j}$. Suppose now that we found
$l_0,\ldots,l_k, k < n$ such that $F_k(T), F_k(T+l_1),\ldots, F_k(T+l_k)$ are linearly independent over $\R$.  Let
 $l \in \N$ be such that  $F_{k+1}(T+l) = \sum_{i=0}^k A_i F_{k+1}(T+l_i), A_i(l)=A_i \in \C$. Then, we have a
linear system
\begin{equation}%
\label{sys:poly}%
 P_j(l) = \sum_{i=0}^k A_iP_j(l_i), j = 0,\ldots,k+1.
\end{equation}%
We can solve the first $k+1$ equations simultaneously for $A_i$ using Cramer's rule. Thus,%
\[%
A_i=\frac{\sum_{j=0}^kb_jP_j(l)}{\det(P_j(l_i))},%
\]%
where $\det(P_j(l_i)), j=0,\ldots,k, i=0,\ldots,k$ is not zero by induction hypothesis and $b_j \in \C$. Therefore,
for each $i=0,\ldots,k$, we have that $A_i=A_i(l)$ is a fixed polynomial in $l$ of degree at most $k$. Next
consider the equation of system (\ref{sys:poly}) number $k+2$.
\[%
P_{k+1}(l) = \sum_{i=0}^k A_i(l)P_{k+1}(l_i).
\]%
Note that on the left we have a polynomial in $l$ of degree $k+1$ and on the right a polynomial of degree at most
$k$. Thus, the equality will not hold for all sufficiently large $l$.  Finally, note that  for any  non-negative
integer $k \leq n$, any $l_0,\ldots,l_k \in \R$, we have that the set $\{F(T+l_0),\ldots, F(T+l_k)\}$ is
linearly dependent only if the set $\{F_k(T+l_0),\ldots, F_k(T+l_k)\}$ is linearly dependent.
\end{proof}%
The next lemma deals with degrees of certain extensions used to define integrality at finite sets of primes.

\begin{lemma}%
\label{le:extq}%
Let $K$ be a field. Let $q$ be a rational prime. Let $b \in K$ be such that $b$ is not a $q$-th power in $K$. Let
$\beta$, an element of the algebraic closure of $K$, be a root of $X^q -b$. Then $[K(\beta):K]=q$.
\end{lemma}%
\begin{proof}%
It is obvious that $[K(\beta):K] \leq q$. So suppose $[K(\beta):K] =m < q$. Let $\beta_1=\beta,\ldots,\beta_m$ be
all the conjugates of $\beta$ over $K$. Observe that $\beta_i= \xi_i\beta$, where $\xi_i$ is a $q$-th root of
unity. Further let
\[%
c={\mathbf N}_{K(\beta)/K}=\prod_{i=1}^m\xi_i\beta_i=\xi\beta^m \in K,
\]%
where $\xi$ is again a $q$-th root of unity. Now let $x,y \in \Z$ be such that $xm + yq=1$.
Then
\[%
c^xb^y=\xi^x\beta^{xm}\beta^{qy}=\xi'\beta \in K,
\]%
where $\xi'$ is another $q$-th root of unity. But $(\xi'\beta)^q= b$ in contradiction of our assumption on $b$.
\end{proof}%
The next lemma uses the Kronecker-Weber Theorem to determine how to embed an abelian extension into the smallest
possible cyclotomic one.
\begin{lemma}%
\label{le:KW}
Let $A_{\infty}$ be an abelian extension of $\Q$ with finitely many ramified rational primes $p_1,\ldots, p_k$.
Then $A_{\infty}$ is contained in the $F=\Q(\xi_{1,l},\ldots,\xi_{k,l}, l \in \Z_{>0})$, where for $i=1,\ldots,k, j
\in \Z_{>0}$ we have that $\xi_{i,j}$ is $p_i^j$-th primitive root of unity.
\end{lemma}%
\begin{proof}%
Suppose the assertion of the lemma does not hold. By Kronecker-Weber theorem, $A_{\infty}$ must be contained in a
cyclotomic extension.  Then for some $\alpha \in A_{\infty}$ we have that $\alpha \not \in F$ but $\alpha \in
F(\xi_n)$, where $\xi_n$ is an $n$-th primitive root of unity and $(n,p_i)=1$ for $i=1,\ldots,k$.    Observe that
the only rational primes ramified in the extension $F(\alpha)/\Q$ are $p_1,\ldots,p_k$ (see Proposition 8, Section
4, Chapter II of \cite{L}).  Next consider $F(\xi_n)/F$ and let $\tau \in \Gal(F(\xi_n)/F)$ be such that
$F(\alpha)$ is the fixed field of the subgroup generated by $\tau$.  Since $F$ and $\Q(\xi_n)$ are linearly
disjoint over $\Q$, we have that  $\Gal(F(\xi_n)/F) \stackrel{\rightarrow}{\cong} \Gal(\Q(\xi_n)/\Q)$ with the
isomorphism realized by restriction. Therefore, restriction of $\tau$ to $\Q(\xi_n)$ will not generate all of
$\Gal(\Q(\xi_n)/\Q)$.   Let $\mu \in \Q(\xi_n)$ generate the fixed field of the subgroup of $\Gal(\Q(\xi_n)/\Q)$
generated by restriction of $\tau$ to $\Q(\xi_n)$.  Then $\mu \not \in \Q$ and $\mu \in F(\alpha)$.  Hence
$\Q(\mu) \subset F(\alpha)$.  But one of the rational divisors of $n$ is ramified in the extension $\Q(\mu)/\Q$
contradicting our assumption on $A_{\infty}$.
\end{proof}%

The following lemmas will all deal with some technical aspects of prime splitting in number fields.%
\begin{lemma}%
\label{le:notsplit}%
 Let $T/K$ be a cyclic extension of number fields, and let $M$ be an  extension of $K$ such that
it is  Galois and $([M:K],[T:K])=1$.  Let $\pp_K$ be a prime of $K$ not splitting in the extension $T/K$.  Let
$\pp_M$ be an $M$-prime above $\pp_K$.  Then $\pp_M$ does not split in the extension $MT/M$.
\end{lemma}%
\begin{proof}%
First of all observe that the extension $MT/K$ is Galois and every factor of $\pp_K$ in $M$ has the same number of
factors in $MT$.  If this number is not 1 then it is a non trivial divisor of $[MT:M]=[T:K]$.  Thus, if the
factors of $\pp_K$ do not stay prime in the extension $MT/M$, the number of $MT$-factors of $\pp_K$ has a
non-trivial common divisor with $[T:K]$.   On the other hand, since $\pp_K$ does not split in the extension $T/K$,
the number of $MT$-factors of $\pp_K$ is a divisor of $[MT:T]=[M:K]$.  Thus, if some factor of $\pp_K$ splits in
the extension $MT/M$, we have a contradiction of our assumption on the degrees of extensions $M/K$ and $T/K$.
\end{proof}%

\begin{lemma}%
\label{le:notsplitabove}%
Let $K/B$ be a Galois extension and let $T/B$ be a cyclic extension. Assume further that $K$ and $T$ are linearly
disjoint over $B$. Let $\pp_B$ be a $B$-prime not splitting in the extension $T/B$ and splitting completely in the
extension $K/B$. Then the following statements are true.
\be%
\item There are infinitely many primes of $B$ satisfying the two requirements for $\pp_B$.%
\item Let $\pp_K$ be a $K$-prime lying above a $B$-prime $\pp_B$ as above. Then $\pp_K$ does not split in the
extension $TK/K$.
\ee%
\end{lemma}%
\begin{proof}%
The linear disjointness implies that
\[%
\Gal(KT/B) \cong \Gal(KT/T)\times \Gal(KT/K) \stackrel{\rightarrow}{\cong} \Gal(K/B) \times \Gal(T/B),
\]%
where the last isomorphism is realized by restriction.  So consider an element $(\id_K,\sigma_T) \in \Gal(KT/B) $,
where $\id_K$ is the identity element of  $\Gal(K/B)$ and $\sigma_T$ is a generator of $\Gal(T/B)$.  Let
$\pp_{KT}$ be a $KT$-prime whose Frobenius isomorphism is $(\id_K,\sigma_T)$.  By Tchebotarev density theorem there
are infinitely many such  primes.  Next let $\pp_K$ be a prime $K$ prime below it.  The decomposition group of
$\pp_{KT}$ over $K$ is the intersection of the decomposition group of $\pp_{KT}$ over $B$ and $\Gal(KT/K)$.  This
intersection is all of $\Gal(KT/K)$ and therefore $\pp_K$ will not split in the extension $KT/K$.  Finally, since
the decomposition groups of $\pp_{KT}$ over $K$ and over $B$ are the same, we conclude that the decomposition group
of $\pp_K$ over $B$ is trivial.  Thus, if $\pp_B$ is the $B$-prime below $\pp_{KT}$, then $\pp_B$ splits completely
in the extension $K/B$.

Next let $\pp_B$ be as above and let $\pp_{KT}$ be its factor  in $KT$.  Then, since the relative degree of $\pp_T$
over $\pp_B$,  $f(\pp_T/\pp_B)=[T:B]$, we have that $[T:B] \geq f(\pp_{KT}/\pp_B) =
f(\pp_{KT}/\pp_{K})f(\pp_K/\pp_B)=f(\pp_{KT}/\pp_K) \leq [T:B]$, where the last equality holds because $\pp_B$
splits completely in $K$.  Thus, $f(\pp_{KT}/\pp_K) = [T:B]$ and $\pp_K$ remains prime in $KT$.
\end{proof}%
\begin{lemma}%
\label{le:split}%
Let $K/B$ be a finite extension of number fields and let $T/B$ be a Galois extension. Assume further that $K$ and
$T$ are linearly disjoint over $B$. Let $\qq_B$ be a $B$-prime splitting completely (into distinct factors) in the
extension $T/B$. Let $\qq_K$ be a $K$-prime lying above $\qq_B$. Then $\qq_K$ splits in the extension $TK/K$.
\end{lemma}%
\begin{proof}%
In this case the linear disjointness implies that
\[%
\Gal(KT/K) \stackrel{\rightarrow}{\cong} \Gal(T/B),
\]%
where the isomorphism, as above, is realized by restriction. Let $\sigma \in \Gal(KT/K)$ be the Frobenius
isomorphism of some $KT$-factor of $\qq_{K}$. Then $\sigma$ restricted to elements of $T$ should be an element
of the decomposition group of $\qq_T$, a factor of $\qq_B$ in $T$. But the decomposition group of any factor of
$\qq_B$ in $T$ is trivial. Thus, since the restriction to $B$ is an isomorphism, we conclude that the
decomposition group of any factor of $\qq_K$ in $KT$ is trivial. Thus $\qq_K$ splits completely.
\end{proof}%

\begin{lemma}%
\label{le:existprime} Let $U/K$ be a Galois extension of number fields.  Let $F_i/U, i=1,\ldots,k$ be a cyclic
number field extension such that each $F_j$ is  linearly disjoint from $\prod_{i\not=j}F_i$ and the extension
$\left(\prod_{i=1}^kF_i\right )/K$ is Galois.   Then there are infinitely many primes of $U$ not splitting in the
extension $F_i/U$ for any $i$ and lying above a prime of $K$ splitting completely in $U$.
\end{lemma}%
\begin{proof}%
The linear disjointness condition implies that $\Gal(\prod_{i=1}^kF_i/U)\cong \prod_{i=1}^k\Gal(F_i/U)$.  Let
$\sigma_i$ be a generator of $\Gal(F_i/U)$.  Then any prime of $\prod_{i=1}^k F_i$ whose Frobenius is
$(\sigma_1,\ldots,\sigma_n) \in \Gal(\prod_{i=1}^kF_i/U) \subset  \Gal(\prod_{i=1}^kF_i/K) $ will have the desired
property.  Now Tchebotarev Density Theorem tells us that there are infinitely many such primes.
\end{proof}%
The following two lemmas are slight generalizations of Lemma 2.6 of \cite{PS}.
\begin{lemma}%
\label{le:frob}%
Let $F/U$ be a cyclic extension such that for some rational prime $q$ we have that $[F:U]=\ell\equiv 0 \mod q$. Let
$N$ be the unique $q$-th degree extension of $U$ contained in $F$. Let $\pp_F$ be a prime of $F$ and let $\pp_U$ be
the $U$-prime below it.  Let  $\sigma$ be the Frobenius automorphism of $\pp_F$.  Then $\pp_U$ splits in $N$
if and only if $\sigma$ is a $q$-th power in $\Gal(F/U)$.
\end{lemma}%
\begin{proof}%
The unique index $q$ subgroup $H$ of $\Gal(F/U)$ consists of $q$-th powers. Further, $N$ is the fixed field of $H$.
Suppose now that $\sigma \in H$.  Then the decomposition group of $\pp_F$ over $U$ (denoted by $G_{\pp_F}(F/U)$)
and $N$ (denoted by $G_{\pp_F}(F/N)$) are the same.  Let $\pp_N$ be the $N$-prime below $\pp_F$.  In this case the
decomposition group of $\pp_N$ over $U$, equal to $G_{\pp_F}(F/U)/G_{\pp_F}(F/N)$, is trivial and $\pp_U$ splits
completely in $N$. Conversely, suppose $\sigma \not \in H$.  Then       $G_{\pp_F}(F/N) \not =G_{\pp_F}(F/U) $ and
$G_{\pp_F}(F/U)/G_{\pp_F}(F/N)$ is not trivial.  Thus, $\pp_U$ does not split completely in the extension $N/U$.
Since the degree of the extension is prime, however, for an unramified prime not splitting completely is equivalent
to staying prime.
\end{proof}%
\begin{lemma}%
\label{le:extensions}%
Let $p_1, \ldots,p_k, t$ be a finite set of distinct rational primes.  Then there exists a totally real cyclic
extension of $\Q$ of degree $t$ where none of $p_i$'s splits and there exists a totally real cyclic extension of
$\Q$ of degree $t$ where all of $p_i$ split.  Further, we can arrange for any given finite subset of primes not to
ramify in these extensions.
\end{lemma}%
\begin{proof}%
Let $\ell$ be a prime splitting completely into distinct factors in the extension
$\Q(\xi_t,\sqrt[t]{p_1},\ldots,\sqrt[t]{p_k})$, where $\xi_t$ is a primitive $t$-th root of unity. Then
$\ell\equiv 1\mod t$ and $\mod \ell$ we have that $p_i$ is a $t$-th power. Now consider the extension
$\Q(\xi_{\ell})/\Q$, where $\xi_{\ell}$ a primitive $\ell$-th root of unity. Let $\tau_i$ be the Frobenius of
$p_i$. Then $\tau_i(\xi_{\ell})=\xi_{\ell}^{p_i}$ and $\tau_i$ is a $t$-th power in $\Gal(\Q(\xi_{\ell})/\Q)$. Let
$G$ be the unique degree $t$ extension of $\Q$ inside of $\Q(\xi_{\ell})$. Then by Lemma \ref{le:frob} we have
that $p_i$ splits completely in this extension, and the first assertion of the lemma holds.

Now let $\ell$ satisfy the following conditions:
\be%
\item $\ell$ splits completely in    $\Q(\xi_t)/\Q$. %
\item Factors of $\ell$ in $\Q(\xi_t)$ do not split in any of the extensions $\Q(\xi_t, \sqrt[t]{p_i})/\Q(\xi_t)$.
(By Lemma \ref{le:existprime} there are infinitely many such $\ell$'s.)
\ee%
Then we conclude that as above $\ell\equiv 1\mod t$, but $p_i$ is not a $t$-th power mod $\ell$. Now considering
the extension $\Q(\xi_{\ell})/\Q$ as above, by Lemma \ref{le:frob}, we conclude that none of $p_i$ will split in
the unique degree $t$ extension of $\Q$ contained in $\Q(\xi_{\ell})$.

Finally we observe that $\ell$ would be the only prime ramifying in either extension, and in choosing $\ell$ we can
always avoid any finite set of primes.
\end{proof}%
\begin{lemma}%
\label{le:infset}%
Let $Z$ be a number field and let $K/Z$ be a finite extension.  Let $\beta \in \tilde{\Q}$ be such that
$\beta^2\in K$, $\beta \not \in K$ and $K(\beta)/Z$ is Galois.  Let $E$ be a cyclic extension of $Z$ of odd
degree $l$ with $(l,[K:Z])=1$. Then there exists an infinite set $\calB_{Z}$ of primes of $Z$ such that every
$K$-prime above a prime of $\calB_{Z}$ does not split in the extension $KE(\beta)/K$ and every prime in
$\calB_{Z}$ splits completely in the extension $K/Z$.
\end{lemma}%
\begin{proof}%
Given our assumptions on the degrees of the extensions, $\Gal(KE(\beta)/K) = \Gal(K(\beta)/K)) \times \Gal(E/Z)$.
Let $\sigma$ be the generator of $\Gal(K(\beta)/K)$ and let $\tau$ be a generator of $\Gal(E/Z)$. Let
$\pp_{EK(\beta)}$ be a prime of $EK(\beta)$ whose Frobenius is $(\sigma,\tau) \in \Gal(EK(\beta)/Z)$. Note that
$(\sigma,\tau)$ generates $\Gal(EK(\beta)/K)$ and therefore, if we let $\pp_K$ be the prime below
$\pp_{EK(\beta)}$ in $K$ we will have that $\pp_K$ does not split in the extension $KE(\beta)/K$. Next let $\qq_K$
be the conjugate of $\pp_K$ over $Z$. Since $KE(\beta)/Z$ is Galois, $\qq_K$ remains prime in the extension
$KE(\beta)/K$ also. Let $\pp_Z$ be the $Z$-prime below $\pp_K$ and $\qq_K$. Then no factor of $\pp_{Z}$ in $K$
splits in the extension $KE(\beta)/K$. Finally we note that the decomposition group of $\pp_{EK(\beta)}$ over $K$
and over $Z$ are the same. Thus, $\pp_{Z}$ must split completely in the extension $K/Z$.
\end{proof}%
\bibliography{mybib}%

\end{document}